\documentclass[10pt]{article}
\usepackage{amsfonts}
\usepackage{amsfonts,latexsym,amsmath,mathrsfs,amscd,geometry}
\usepackage{amssymb}
\usepackage{latexsym}
\usepackage{color}
\usepackage{upgreek}
\usepackage{amsmath}
\numberwithin{equation}{section}

\newcommand \nc{\newcommand}
\newtheorem{theorem}{Theorem}[section]
\newtheorem{lemma}[theorem]{Lemma}
\newtheorem{proposition}[theorem]{Proposition}

\newtheorem{remark}[theorem]{Remark}

\nc{\ba}{\begin{array}}\nc{\ea}{\end{array}}
\nc{\be}{\begin{eqnarray}}\nc{\ee}{\end{eqnarray}}
\nc{\beq}{\begin{equation}}\nc{\eeq}{\end{equation}}
\nc{\bex}{\begin{eqnarray*}}\nc{\eex}{\end{eqnarray*}}
\nc{\btm}{\begin{theorem}} \nc{\etm}{\end{theorem}}
\nc{\blm}{\begin{lemma}} \nc{\elm}{\end{lemma}}
\nc{\R}{\mathbb{R}} \nc{\va}{\varepsilon} \nc{\ls}{\limits}

\def\pf{\noindent{\bf Proof.\quad}}\def\endpf{\hfill$\Box$}

\allowdisplaybreaks

\usepackage[colorlinks,linkcolor=blue]{hyperref}

\begin{document}
\title{The Cauchy problem for an inviscid Oldroyd-B model in $\mathbb{R}^3$}
\author{Sili Liu\footnote{School of Mathematics, South China University of Technology, Guangzhou 510631, China. E-mail: maslliu@mail.scut.edu.cn},\, Wenjun Wang\footnote{College of Science,
University of Shanghai for Science and Technology,
Shanghai 200093, China. E-mail:
wwj001373@hotmail.com},\, Huanyao Wen
\footnote{School of Mathematics, South China University of Technology, Guangzhou 510631, China. Email: mahywen@scut.edu.cn.}
}

\maketitle

\begin{abstract}
In this paper, we consider the Cauchy problem for an inviscid compressible Oldroyd-B model in three dimensions. The global well posedness of
strong solutions and the associated time-decay estimates in Sobolev spaces are established near an equilibrium state. The vanishing of viscosity is the main challenge compared with our previous work \cite{W1W} where the viscosity coefficients are included and the decay rates for the highest-order derivatives of the solutions seem not optimal. One of the main objectives of this paper is to develop some new dissipative estimates such that the smallness of the initial data and decay rates are independent of the viscosity. In addition, it proves that the decay rates for the highest-order derivatives of the solutions are optimal. Our proof relies on Fourier theory and delicate energy method. This work can be viewed as an extension of \cite{W1W}.
\end{abstract}

\bigbreak \textbf{{\bf Key Words}:}  Inviscid compressible Oldroyd-B model; global well posedness; decay rates.

\bigbreak  {\textbf{AMS Subject Classification 2020:} 76N10, 76N17, 74H40.

\section{Introduction}
The Oldroyd-B model is a widely used constitutive model to describe the motion of viscoelastic fluids.
One of the known derivations is that it can be derived as a macroscopic closure of Navier-Stokes-Fokker-Planck system which is a micro-macro
model describing dilute polymeric fluids in dumbbell Hookean setting, see \cite{Bris-Lelievre} and \cite{BLS} for the incompressible case and the compressible case, respectively. The compressible Oldroyd-B model in the space-time cylinder $Q_T=\mathbb{R}^3\times(0,T]$ is stated as follows:
\beq\label{clc-1}
\left\{
\begin{split}
&\rho_t+\mathrm{div} (\rho u)=0,\\
&(\rho u)_t+\mathrm{div} (\rho u\otimes u)+\nabla P(\rho)-\mu\Delta u-(\mu+\nu)\nabla \mathrm{div}u
=\mathrm{div} \big(\mathbb{T}-(kL\eta+\mathfrak{z}\eta^{2})\mathbb{I}\big),\\
&\eta_t+\mathrm{div} (\eta u)=\varepsilon\Delta \eta,\\
&\mathbb{T}_t+\mathrm{div} (u\mathbb{T})-(\nabla u
\mathbb{T}+\mathbb{T}\nabla^Tu)
=\varepsilon\Delta \mathbb{T}+\frac{kA_0}{2\lambda}\eta\mathbb{I}-\frac{A_0}{2\lambda}\mathbb{T},
\end{split}
\right.
\eeq
where the pressure $P(\rho)$ and the density $\rho(x,t)\geq0$ of the fluid are supposed to be related by the typical power law relation for simplicity:
\[P(\rho)=a\rho^{\gamma}\]
for some known constants $a>0$, $\gamma>1$; $u(x,t)\in\mathbb{R}^{3}$ denotes the velocity field of the fluid. $\mu$ and $\nu$ are viscosity coefficients satisfying $\mu\geq0,\; 2\mu+3\nu\geq0$. The polymer number density $\eta(x,t)\geq0$ represents the integral of
the probability density function $\psi$ which is a microscopic variable in the modelling of dilute polymer chains, i.e.,
\[\eta(x,t)=\int_{\mathbb{R}^3}\psi(x,t,q)\,\mathrm{d}q,\]
where $\psi$ is governed by the Fokker-Planck equation. The extra stress tensor $\mathbb{T}(x,t)= (\mathbb{T}_{i,j})(x,t)\in\mathbb{R}^{3\times 3},\ 1\leq i,j\leq 3$ is a positive definite symmetric matrix defined on $Q_T$, and the notation $\mathrm{div} (u\mathbb{T})$ is understood as
\begin{equation*}
\left(\mathrm{div} (u\mathbb{T})\right)_{i,j} = \mathrm{div} (u\mathbb{T}_{i,j}), \quad 1\leq i,j \leq 3.
\end{equation*}
The constant parameter $\varepsilon$ is the centre-of-mass diffusion coefficient and other parameters $k, L, \mathfrak{z}, A_0, \lambda$ are all positive and known numbers, whose meanings were explained in \cite{BLS}. (\ref{clc-1}) is known as diffusive Oldroyd-B model when the diffusion coefficient $\varepsilon>0$. The corresponding micro-macro version of (\ref{clc-1}) can be referred for instance to \cite{Barrett-Suli} and references therein.

Note that the centre-of-mass diffusion term is usually smaller than other effects (\cite{Bhave1}). For such a reason, in early mathematical studies of macroscopic Oldroyd--B model, the stress diffusion is omitted, see \cite{Oldroyd}. In this context, \cite{Renardy90} established the local existence theory for Dirichlet problem. Guillop\'e and Saut \cite{GS1} obtained the existence and uniqueness of global strong solution in the Sobolev space $H^s(\Omega)$ for bounded domains $\Omega\in \mathbb{R}^3$ with a small initial data. Some other related results can be referred to \cite{FGO,MT}. In exterior domains, Hieber, Naito and Shibata \cite{HNS} obtained a global existence and uniqueness of the solution provided the initial data and the coupling constant are sufficiently small. Fang, Hieber and Zi \cite{FHZ} extended the work \cite{HNS} to the case without any smallness assumption on the coupling constant. The existence of a global-large-data weak solution was established by Lions and Masmoudi (\cite{LM}) in the corotational derivative setting. For long-time behavior of the solution, please refer to \cite{HWZ, HWWZ}. There are also some interesting results on other macroscopic model of Oldroyd type concerning viscoelastic flow introduced by Lin, Liu and Zhang (\cite{Lin-Liu-Zhang}), for example \cite{Hu-Lin, Hu-Wu1, Lai} and others.

However, the center-of-mass diffusion can be physically justified
 to model the shear and vorticity banding phenomena (\cite{Bhave2,Cates_2006, Chupin, Dhont_2008,El-Kareh, Liu, Malek2018}), although it is small. In this case, some interesting works have been achieved. More specifically, the global-in-time existence of large-data solutions in two dimensional setting was obtained by Barrett-Boyaval (\cite{Barrett-Boyaval}) for weak solutions and by  Constantin and Kliegl (\cite{CK}) for strong solutions. In three-dimensional setting, Bathory, Bul\'{i}\v{c}ek and M\'{a}lek (\cite{Bathory_etal}) proved the global existence of weak solutions for a generalized rate-type viscoelastic fluids in bounded domains. For the inviscid case, Elgindi and Rousset (\cite{ER}) obtained the global existence and uniqueness of regular solutions in two dimensions with arbitrarily large initial data when $Q=Q(\nabla u,\mathbb{T})$ is omitted and with small initial data when $Q\neq0$. We refer to \cite{EL} for the three-dimensional case with small initial data. For the case of fractional Laplace, please refer to \cite{Constantin-3}. Very recently, the second author, the third author and their collaborators (\cite{HWWZ}) studied the long-time behavior of the solutions and obtained some decay estimates. These results are concerned with homogeneous fluids, i.e., the density is constant.

For the compressible case, there are a lot of fundamental problems which are still open. We recall some mathematical results for compressible viscoelastic models, which have been the subject of related fields in recent years. The well posedness in local time and global well posedness near an equilibrium for macroscopic models of three-dimensional compressible viscoelastic fluids were considered in \cite{FZ, Hu-Wang3, Hu-Wu, Lei, Qian-Zhang} (see \cite{BuFeMa} for global existence of weak solutions). In particular, Fang and Zi (\cite{FZ}) proved the local well posedness of strong solutions to a compressible Oldroyd--B model and established a blow-up criterion. Soon afterwards, the authors (\cite{Z}) obtained the global well posedness in critical spaces. Lei (\cite{Lei}), Fang and Zi (\cite{FZ1}), and Guillop\'e, Salloum and Talhouk \cite{GST} investigated the incompressible limit problem in torus, the whole space and bounded domain, respectively. Very recently, Zhou, Zhu and Zi (\cite{ZZZ}) obtained some time-decay estimates of strong solutions. Zhu \cite{Zhu} obtained the global well posedness of small classical solutions to a generalized inviscid compressible Oldroyd--B model in Sobolev space $H^s$ for $s\geq5$. In \cite{BLS}, Barrett, Lu and S\"uli not only showed the derivation of the compressible viscous Oldroyd--B model with stress diffusion \eqref{clc-1} via a macroscopic closure of a micro-macro model, but also proved the existence of global-in-time finite-energy weak solutions with arbitrarily large initial data in two dimensions. The global-in-time existence of solutions strong or weak with arbitrarily large initial data is unknown in three dimensions either with stress diffusion or not. In two and three dimensional setting, Lu and Zhang (\cite{LZ}) obtained the local-in-time well posedness of strong solutions together with a blow-up criterion and weak-strong uniqueness. Very recently, the second author and the third author (\cite{W1W}) showed the global well posedness and optimal time-decay rates of strong solutions for Cauchy problem in three dimensions. In critical Besov spaces, one can refer to \cite{Zhai20}. Less is known concerning the vanishing of centre-of-mass diffusion and the inviscid case in \eqref{clc-1} either for global well posedness or for long time behavior, until very recently the first author, the third author and their collaborator investigated the first case (i.e., $\varepsilon=0$) in \cite{LLW}. This work is devoted to the latter one which is more challenging.



 More precisely, we consider the case that $\mu=\nu=0$ in \eqref{clc-1}, i.e.,
\beq\label{clc-3}
\left\{
\begin{split}
&\rho_t+\mathrm{div} (\rho u)=0,\\
&(\rho u)_t+\mathrm{div} (\rho u\otimes u)+\nabla  P(\rho)
=\mathrm{div} (\mathbb{T}-(kL\eta+\mathfrak{z}\eta^2)\mathbb{I}),\\
&\eta_t+\mathrm{div} (\eta u)=\varepsilon\Delta \eta,\\
&\mathbb{T}_t+\mathrm{div} (u\mathbb{T})-(\nabla u\mathbb{T}+\mathbb{T}\nabla^Tu)
=\varepsilon\Delta \mathbb{T}+\frac{kA_0}{2\lambda}\eta\mathbb{I}-\frac{A_0}{2\lambda}\mathbb{T}.
\end{split}
\right.
\eeq
\subsection{Reformulation of the problem}
In this section, we give a reformulation of (\ref{clc-3}) to make the analysis more convenient right behind. In fact, when $\varepsilon=0$ and the viscosity coefficients are fixed, a similar reformulation was given in our previous work \cite{LLW}. Thus this section is a slight modification of the corresponding part in \cite{LLW}. More specifically, multiplying $(\ref{clc-3})_3$ by $k\mathbb{I}_{ij}$, we have
\beq\label{clc-2}
(k\eta \mathbb{I}_{ij})_t+\mathrm{div} (k\eta \mathbb{I}_{ij} u)=\varepsilon\Delta(k\eta\mathbb{I}_{ij}).
\eeq
Then subtracting $(\ref{clc-2})$ from $(\ref{clc-3})_4$ yields that
\beq\label{clc-4}
\partial_t(\mathbb{T}_{ij}-k\eta \mathbb{I}_{ij})+\partial_l \left((\mathbb{T}_{ij}-k\eta \mathbb{I}_{ij})u_l\right)
-(\partial_l u_i\mathbb{T}_{lj}+\mathbb{T}_{il}\partial_l u_j)
=\varepsilon\Delta(\mathbb{T}_{ij}-k\eta \mathbb{I}_{ij})-\frac{A_0}{2\lambda}(\mathbb{T}_{ij}-k\eta \mathbb{I}_{ij}).
\eeq
Further, denoting $\tau_{ij}=\mathbb{T}_{ij}-k\eta \mathbb{I}_{ij}$, and conducting direct calculations, we can get
\beq\label{clc-5}
\partial_l u_i\mathbb{T}_{lj}=\partial_l u_i\tau_{lj}+k\partial_l u_i\eta \mathbb{I}_{lj}=\partial_l u_i\tau_{lj}+k\partial_j u_i\eta,
\eeq
and
\beq\label{clc-6}
\mathbb{T}_{il}\partial_l u_j=\tau_{il}\partial_l u_j+k\eta \mathbb{I}_{il}\partial_l u_j=\tau_{il}\partial_l u_j+k\eta\partial_i u_j.
\eeq
Putting $(\ref{clc-5})$ and $(\ref{clc-6})$ into $(\ref{clc-4})$ yields
\begin{equation*}
\partial_t\tau_{ij}+\partial_l (\tau_{ij}u_l)
-(\partial_l u_i\tau_{lj}+\tau_{il}\partial_l u_j)-k\eta(\partial_j u_i+\partial_i u_j)
=\varepsilon\Delta\tau_{ij}-\frac{A_0}{2\lambda}\tau_{ij},
\end{equation*}
which is
\beq\label{clc-7}
\partial_t\tau+\mathrm{div} (u\,\tau)
-(\nabla u\tau+\tau\nabla^T u)-k\eta(\nabla u+\nabla^T u)
=\varepsilon\Delta\tau-\frac{A_0}{2\lambda}\tau.
\eeq
Next, the term of the right-hand side in $(\ref{clc-3})_2$ can be transformed into the following form
\begin{equation*}
\partial_j \big(\mathbb{T}_{ij}-(kL\eta+\mathfrak{z}\eta^2)\mathbb{I}_{ij}\big)
=\partial_j \big(\tau_{ij}+k\eta\mathbb{I}_{ij}-(kL\eta+\mathfrak{z}\eta^2)\mathbb{I}_{ij}\big)
=\partial_j \tau_{ij}-\partial_i\left(k(L-1)\eta+\mathfrak{z}\eta^2\right),
\end{equation*}
which together with $(\ref{clc-3})_2$ implies that
\beq\label{clc-8}
(\rho u)_t+\mathrm{div} (\rho u\otimes u)+\nabla \left(P(\rho)+k(L-1)\eta+\mathfrak{z}\eta^2\right)
=\mathrm{div} \tau.
\eeq
Finally, combining $(\ref{clc-3})_1$, $(\ref{clc-8})$, $(\ref{clc-3})_3$ and $(\ref{clc-7})$ yields
\beq\label{1}
\left\{
\begin{split}
&\rho_t+\mathrm{div} (\rho u)=0,\\
&(\rho u)_t+\mathrm{div} (\rho u\otimes u)+\nabla \big( P(\rho)+k(L-1)\eta+\mathfrak{z}\eta^{2}\big)
=\mathrm{div} \tau,\\
&\eta_t+\mathrm{div} (\eta u)=\varepsilon\Delta \eta,\\
&\tau_t+\mathrm{div} (u\tau)-(\nabla u
\tau+\tau\nabla^Tu)-k\eta(\nabla u+\nabla^Tu)
=\varepsilon\Delta \tau-\frac{A_0}{2\lambda}\tau,
\end{split}
\right.
\eeq
which is equipped with the following initial condition:
\begin{equation}\label{initial-condition}
(\rho, u, \eta, \tau)(x,t)|_{t=0}=(\rho_0, u_0, \eta_0, \tau_0)(x)\rightarrow(\tilde{\rho},0,\tilde{\eta},0),\quad {\rm as}\,\, |x|\rightarrow \infty.
\end{equation}
Note that (\ref{1}) is equivalent to the system (\ref{clc-3}) with the regularity of the solution in the present paper and that it seems more convenient to consider (\ref{1}) in the proof. Therefore we will state the main results afterwards for the reformulated system (\ref{1}) only.

\subsection{Main Results}
Our main results are stated as follows.
\begin{theorem}\label{Theorem1} Let $L\geq 1,\,\mathfrak{z}\geq0$. Assume that $(\rho_0-\tilde{\rho}, u_0, \eta_0-\tilde{\eta}, \tau_0)\in H^3(\mathbb{R}^3)$  for constants $\tilde{\rho}, \tilde{\eta}>0$. Then there exists a positive constant $\theta$ sufficiently small such that if
\begin{align}
\|(\rho_0-\tilde{\rho}, u_0, \eta_0-\tilde{\eta}, \tau_0)\|_{H^3(\mathbb{R}^3)}\leq \theta,\label{T1}
\end{align}
the initial-value problem $(\ref{1})$-$(\ref{initial-condition})$ admits a unique global strong solution $(\rho, u, \eta, \tau)$ which satisfies
\begin{equation*}
\begin{split}
&(\rho-\tilde{\rho},u)\in
\mathcal{C}([0,\infty);H^3(\mathbb{R}^3)),\ (\rho_t,u_t)\in \mathcal{C}([0,\infty); H^{2}(\mathbb{R}^3)),\,\rho>0,\ \eta> 0,\\
& (\eta-\tilde{\eta},\tau)\in \mathcal{C}([0,\infty);H^3(\mathbb{R}^3))\cap L^2(0,\infty;H^{4}(\mathbb{R}^3)),\, (\eta_t,\tau_t)\in  \mathcal{C}([0,\infty);H^1(\mathbb{R}^3))\cap L^2(0,\infty;H^{2}(\mathbb{R}^3)).
\end{split}
\end{equation*}
\end{theorem}

\begin{theorem}\label{Theorem2} In addition to the conditions of Theorem \ref{Theorem1}, we assume that $(\rho_0-\tilde{\rho}, u_0, \eta_0-\tilde{\eta}, \mathrm{div}\tau_0)\in L^1(\mathbb{R}^3)$.
Then there exists a positive constant $C$ independent of $t$ such that the solution $(\rho,u,\eta,\tau)$ satisfies the following time-decay estimates:
\begin{align*}
&\|\nabla^m\tau(t)\|_{L^2(\mathbb{R}^3)}\leq C(1+t)^{-\frac{5}{4}-\frac{m}{2}},\;\;\;m=0,1,2,\\
&\|\nabla^m(\rho-\tilde{\rho}, u, \eta-\tilde{\eta})(t)\|_{L^2(\mathbb{R}^3)}\leq C(1+t)^{-\frac{3}{4}-\frac{m}{2}},\;\;\;m=0,1,2,3,\\
&\|\nabla^3\tau(t)\|_{L^2(\mathbb{R}^3)}\leq C(1+t)^{-\frac{9}{4}},
\end{align*}for any $t\geq0$.

\end{theorem}

\begin{remark}
From (\ref{1})$_4$ and the conclusion in Theorem \ref{Theorem2}, one can observe that the decay rate of $\|\nabla^l\tau\|_{L^2(\mathbb{R}^3)}$ is the same as that of $\|\nabla^{l+1}u\|_{L^2(\mathbb{R}^3)}$ where the maximum of $l$ is 2 according to the regularity of the solutions. Therefore the decay rate of $\|\nabla^3\tau(t)\|_{L^2(\mathbb{R}^3)}$ is not expected to be sharper.
\end{remark}

\medskip

We would like to introduce the main novelty of this work. Compared with \cite{W1W} where the global well posedness of strong solution for \eqref{clc-1} with positive shear viscosity $\mu$ is established subject to some smallness assumptions, the vanishing of viscosity in the present paper will bring new challenges such as the loss of regularity for the velocity. We introduce a good unknown $\tau_{ij}=\mathbb{T}_{ij}-k\eta \mathbb{I}_{ij}$ for $i,j=1,2,3$ inspired by \cite{LLW, Lu-Pokorny20} and derive some new dissipative estimates of velocity from the equation of $\tau_{ij}$ such that the smallness of the initial data does not depend on the viscosity. This demonstrates that the coupling yields new dissipative effect of the velocity satisfying the Euler equation only. Besides, the construction of the basic energy estimates in $H^3$-norm relies strongly on the dissipative estimate of $\nabla^4u$ due to the second term of the continuity equation and the pressure term of the momentum equation. It seems impossible to get the dissipative estimate of $\nabla^4u$ in the non-viscous case. To handle the issue, inspired by \cite{Zhu}, we use the variation of the continuity equation $\mathrm{div}u=-\frac{\rho_t+\beta u\cdot\nabla\rho}{r_1+\beta\rho}$ and integration by parts to transfer the derivative to other term. Concerning the optimal time-decay estimates, the loss of dissipation of velocity due to the vanishing of viscosity is the main difficulty compared with the viscus case in \cite{W1W}. Delicate energy method and low-high-frequency decomposition is the main tool in the proof. In this context, we observe that the reformulated equation of $\tau_{ij}$ can produce the dissipation mechanism of velocity such that the decay estimates do not rely on the viscosity, which is the key to obtain the optimal time-decay estimates of the solution except for its highest order. It is unusual that the optimal decay rate for the highest-order derivatives of the solution to some hyperbolic-parabolic systems even as (\ref{clc-1}) with viscosity (see \cite{W1W}) could be obtained. To get the dissipative estimate for the hyperbolic quantities $\nabla^k\rho$ and $\nabla^k\tau$ where $k=3$ is the maximal one, the usual energy method is to construct the interaction energy functional between the second-order and the third-order. Therefore it implies the decay rate for the third-order will be the same as that for the second-order. Here we use the low-high-frequency decomposition and employ the high-frequency part of velocity at $2^{th}$ order as a test function of the equation of $\nabla^2\mathrm{div}\tau_{ij}$. The high-frequency quantity will make the damping term in the equation of $\nabla^2\mathrm{div}\tau_{ij}$ keep the desired order, see Section \ref{sec4.4} for more details. This is different from our previous work \cite{WW} for compressible Navier-Stokes equations with reaction diffusion where a new observation for cancellation of a low-medium-frequency quantity was adopted to get the optimal time-decay estimate at the highest order, see also \cite{WZZ} for a two-phase fluid model. In addition, to get the decay estimates of the low-frequency part to the linearized system \eqref{A1}, we apply the Hodge decomposition to $u$ and $\tau$, and transfer the linearized system into two system \eqref{A2} and \eqref{A3}. We introduce some corrected modes different from \cite{W1W} to overcome the difficulties caused by the lack of dissipation of $u$. With the help of these estimates, the decay properties for the low-frequency part of the solutions to the nonlinear system are obtained by using the Duhamel principle. Combining the delicate energy estimates with the decay estimates of low-frequency part, we obtain the same decay rates of the solution to \eqref{1} up to the second-order as those for viscous case, see \cite{W1W}. Moreover the decay rate for the third-order in the present paper is sharper.

The rest of the paper is organized as follows. In Section \ref{sec2}, we linearize the reformulated system which will make the following analysis more convenient. In Section \ref{sec3}, the proof of the global well posedness of the solutions will be given by using delicate energy method combined with the continuity technique. In Section \ref{sec4}, we establish some optimal time-decay estimates and finish the proof of Theorem \ref{Theorem2}.

\section{Linearization of the reformulated system}\label{sec2}
To simplify the proof of the main theorems, we linearize the reformulated problem (\ref{1})-(\ref{initial-condition}) as follows.
Taking change of variables by $(\rho, u, \eta, \tau)\rightarrow \left(\rho'+\tilde{\rho}, \beta u', \eta'+\tilde{\eta}, \tau\right)$ with $\beta>0$ to be determined, the initial-value problem $(\ref{1})$-$(\ref{initial-condition})$ is written as below
$$
\left\{
\begin{split}
&\rho'_t+\beta\tilde{\rho}\,\mathrm{div}u'=S_1',\\
&\beta u'_t+\frac{P'(\tilde{\rho})}{\tilde{\rho}}\nabla\rho'
+\frac{k(L-1)+2\mathfrak{z}\tilde{\eta}}{\tilde{\rho}}\nabla\eta'
-\frac{\mathrm{div}\tau}{\tilde{\rho}}=S_2',\\
&\eta'_t+\beta\tilde{\eta}\,\mathrm{div}u'-\varepsilon\Delta\eta'=S_3',\\
&\tau_t+\frac{A_0}{2\lambda}\tau-\varepsilon\Delta\tau-\beta k\tilde{\eta}(\nabla u'+\nabla^Tu')=S_4',
\end{split}
\right.
$$
where
$$
\left\{
\begin{split}
&S_1'=-\beta\mathrm{div}(\rho'u'),\\
&S_2'=-\beta^2u'\cdot\nabla u'-(\frac{P'(\rho'+\tilde{\rho})} {\rho'+\tilde{\rho}}-\frac{P'(\tilde{\rho})}{\tilde{\rho}})\nabla \rho'-(\frac{k(L-1)+2\mathfrak{z}(\eta'+\tilde{\eta})} {\rho'+\tilde{\rho}}-\frac{k(L-1)+2\mathfrak{z}\tilde{\eta}}{\tilde{\rho}})\nabla \eta'\\
&\qquad\;+(\frac{1}{\rho'+\tilde{\rho}}-\frac{1}{\tilde{\rho}})\mathrm{div} \tau,\\
&S_3'=-\beta\mathrm{div}(\eta'u'),\\
&S_4'=-\beta\mathrm{div} (u'\tau)+\beta(\nabla u'
\tau+\tau\nabla^Tu')+\beta k \eta'(\nabla u'+\nabla^Tu'),
\end{split}
\right.
$$
with initial data
$$
(\rho', u', \eta', \tau)(x,0)=(\rho_0',  u_0', \eta_0', \tau_0)(x)\rightarrow(0,0,0,0),\quad {\rm as}\,\, |x|\rightarrow \infty.
$$
Denote the scaled parameters and constants by
$$r_1=\sqrt{P'(\tilde{\rho})},\;\; r_2=\frac{k(L-1)+2\mathfrak{z}\tilde{\eta}}{\sqrt{P'(\tilde{\rho})}},\;\; r_3=\frac{1}{\sqrt{P'(\tilde{\rho})}},\;\; \beta=\frac{\sqrt{P'(\tilde{\rho})}}{\tilde{\rho}}$$
and define the nonlinear functions of $\rho'$ by
$$h(\rho')=\big(\frac{P'(\tilde{\rho})}{\tilde{\rho}}-\frac{P'(\rho'+\tilde{\rho})} {\rho'+\tilde{\rho}}\big)\frac{1}{\beta},\;\;
g(\rho')=\big(\frac{1}{\tilde{\rho}}-\frac{1}{\rho'+\tilde{\rho}}\big)\frac{1}{\beta}.$$

%
Finally, (we remove all $'$ in the following system for brevity) we rewrite the system $(\ref{1})$-$(\ref{initial-condition})$ with linearized part on the left as
\beq\label{17}
\left\{
\begin{split}
&\rho_t+r_1\mathrm{div}u=S_1,\\
&u_t+r_1\nabla\rho+r_2\nabla\eta-r_3\mathrm{div}\tau=S_2,\\
&\eta_t+\beta\tilde{\eta}\,\mathrm{div}u-\varepsilon\Delta\eta=S_3,\\
&\tau_t+\frac{A_0}{2\lambda}\tau-\varepsilon\Delta\tau-\beta k\tilde{\eta}(\nabla u+\nabla^Tu)=S_4,
\end{split}
\right.
\eeq
and
\beq\label{18}
\left\{
\begin{split}
&S_1=-\beta\mathrm{div}(\rho u),\\
&S_2=-\beta u\cdot\nabla u+h(\rho)\nabla \rho+g(\rho)\big[\big(k(L-1)+2\mathfrak{z}\tilde{\eta}\big)\nabla\eta-\mathrm{div}\tau\big]
-\frac{2\mathfrak{z}}{\beta(\rho+\tilde{\rho})}\eta\nabla\eta,\\
&S_3=-\beta\mathrm{div}(\eta u),\\
&S_4=-\beta\mathrm{div} (u\tau)+\beta(\nabla u
\tau+\tau\nabla^Tu)+\beta k\eta(\nabla u+\nabla^Tu),
\end{split}
\right.
\eeq
with initial data
\begin{equation}\label{initial-condition1}
(\rho, u, \eta, \tau)(x,t)|_{t=0}=(\rho_0, u_0, \eta_0, \tau_0)(x)\rightarrow(0,0,0,0),\quad \text{as}\,\, |x|\rightarrow \infty.
\end{equation}

It is worth noticing that the proof of Theorems \ref{Theorem1} and \ref{Theorem2} can be translated into that for the solution to $(\ref{17})$-$(\ref{initial-condition1})$.

\section{Proof of Theorem \ref{Theorem1}}\label{sec3}
In this section, we will prove Theorem \ref{Theorem1} via taking vanishing viscosity limit of (\ref{clc-1}). In fact, the global existence and uniqueness of solutions to the corresponding viscous case has been achieved by the second author and the third author in \cite{W1W} where the smallness of initial data depends on the viscosity coefficients. Therefore the aim in this section is to derive some a {\it priori} estimates globally in time subject to some smallness of data independent of the viscosity coefficients. We assume that $\mu,\nu\leq1$ in the section for simplicity.

After conducted a reformulation similar to (\ref{17}), \eqref{clc-1} can be converted to the following form.
\beq\label{-17}
\left\{
\begin{split}
&\rho_t+r_1\mathrm{div}u=S_1,\\
&u_t+r_1\nabla\rho+r_2\nabla\eta-r_3\mathrm{div}\tau-\mu_1\Delta u-\mu_2\nabla \mathrm{div}u=\tilde{S}_2,\\
&\eta_t+\beta\tilde{\eta}\,\mathrm{div}u-\varepsilon\Delta\eta=S_3,\\
&\tau_t+\frac{A_0}{2\lambda}\tau-\varepsilon\Delta\tau-\beta k\tilde{\eta}(\nabla u+\nabla^Tu)=S_4,
\end{split}
\right.
\eeq
where
\[\mu_1=\frac{\mu}{\tilde{\rho}}, \, \mu_2=\frac{\mu+\nu}{\tilde{\rho}},\]
and
\begin{align}
\tilde{S}_2=&-\beta u\cdot\nabla u+h(\rho)\nabla \rho+g(\rho)\big[\big(k(L-1)+2\mathfrak{z}\tilde{\eta}\big)\nabla\eta-\mathrm{div}\tau\big]
-\frac{2\mathfrak{z}}{\beta(\rho+\tilde{\rho})}\eta\nabla\eta-\mu\beta g(\rho)\Delta u\label{-18}\\
&-(\mu+\nu)\beta g(\rho)\nabla\mathrm{div} u,\notag
\end{align}


We begin with a local existence and uniqueness result of the initial-value problem $(\ref{-17})$ and $(\ref{initial-condition1})$.
\begin{proposition}\label{Proposition0}(local existence and uniqueness) Assume that
\begin{equation*}
(\rho_0, u_0, \eta_0, \tau_0)\in H^3(\mathbb{R}^3),\,\,\inf\limits_{x\in\mathbb{R}^3}\{\rho_0(x)+\tilde{\rho}, \eta_0(x)+\tilde{\eta}\}>0.
\end{equation*}
Then, there exists a constant $T_0 > 0$ depending on $\mu$, $\nu$ and $\|(\rho_0, u_0, \eta_0, \tau_0)\|_{H^3(\mathbb{R}^3)}$, such that the initial-value problem $(\ref{-17})$ and $(\ref{initial-condition1})$ has a unique strong solution $(\rho^{\mu,\nu}, u^{\mu,\nu}, \eta^{\mu,\nu}, \tau^{\mu,\nu})$ over $\mathbb{R}^3\times [0,T_0]$, which satisfies
\begin{equation*}
\begin{split}
&\rho^{\mu,\nu}\in
\mathcal{C}([0,T_0];H^3(\mathbb{R}^3)),\ \rho^{\mu,\nu}_t\in \mathcal{C}([0,T_0]; H^{2}(\mathbb{R}^3)),\ \inf\limits_{Q_{T_0}}(\rho^{\mu,\nu}+\tilde{\rho},\eta^{\mu,\nu}+\tilde{\eta})>0,\\
& (u^{\mu,\nu},\eta^{\mu,\nu},\tau^{\mu,\nu})\in \mathcal{C}([0,T_0];H^3(\mathbb{R}^3))\cap L^2(0,T_0;H^{4}(\mathbb{R}^3)),\\
&(u^{\mu,\nu}_t,\eta^{\mu,\nu}_t, \tau^{\mu,\nu}_t)\in  \mathcal{C}([0,T_0];H^1(\mathbb{R}^3))\cap L^2(0,T_0;H^{2}(\mathbb{R}^3)),
\end{split}
\end{equation*} where $Q_{T_0}=\mathbb{R}^3\times(0,T_0)$.
\end{proposition}
\pf
The proof can be achieved by using some standard iteration arguments, please refer for instance to $\cite{FZ,HDW,LZ}$. We omit the details for brevity.
\endpf

\bigskip
\begin{proposition}\label{proposition4} ({\it A\, priori\, estimate})  Under the assumptions of Theorem \ref{Theorem1}, there exists a positive constant $\delta$ independent of $\mu$ and $\nu$ and at least bigger than $\frac{3\theta}{2}$ $($determined by $(\ref{b18})$ for some $\theta$ given by $(\ref{T1})$$)$, such that if the strong solution of the initial-value problem $(\ref{-17})$ and $(\ref{initial-condition1})$ satisfies
\beq\label{clc9}
\|(\rho^{\mu,\nu}, u^{\mu,\nu}, \eta^{\mu,\nu},\tau^{\mu,\nu})(t)\|_{H^3}\leq \delta,\\
\eeq
for any $t\in[0,T]$, where $0 <T^*\le+\infty$ is the maximum existence time for the solution and $T\in(0,T^*)$, then the following estimates
\begin{align}
&\|(\rho^{\mu,\nu}, u^{\mu,\nu}, \eta^{\mu,\nu}, \tau^{\mu,\nu})(t)\|_{H^3}^2
+\int_0^t\big(\|\nabla(\rho^{\mu,\nu},u^{\mu,\nu})\|_{H^2}^2
+\|\nabla(\eta^{\mu,\nu},\tau^{\mu,\nu})\|_{H^3}^2
\big)\mathrm{d}s\label{clc9-}\\
&+\int_0^t\big(\mu_1\|\nabla u^{\mu,\nu}\|_{H^3}^2+\mu_2\|\mathrm{div}u^{\mu,\nu}\|_{H^3}^2\big)\mathrm{d}s
\leq C\|(\rho_0, u_0, \eta_0, \tau_0)\|_{H^3}^2\leq \frac{2\delta}{3},\notag
\end{align} holds for any $t\in[0,T]$.
\end{proposition}

\begin{remark} A similar result has been obtained by the second author and the third author in \cite{W1W} (Proposition 3.2) when $\delta$ depends on $\mu$ and $\nu$. Proposition \ref{proposition4} removes the dependence between $\delta$ and the viscosity coefficients, which gives the possibilities to take the vanishing viscosity limit.
\end{remark}

\medskip

Based on the Propositions \ref{Proposition0} and \ref{proposition4}, the global existence of solutions to the initial-value
problem $(\ref{-17})$ and $(\ref{initial-condition1})$ will be established with the help of the standard continuity arguments. Then, with the aid of the uniform
estimates \eqref{clc9-} and some compactness arguments, we conclude that a subsequence of solutions $(\rho^{\mu,\nu}, u^{\mu,\nu}, \eta^{\mu,\nu},\tau^{\mu,\nu})$
converges to a limit $(\rho,u,\eta,\tau)$ (in some strong sense) which is a strong solution
to the original problem $(\ref{17})$-$(\ref{initial-condition1})$. Therefore to prove Theorem \ref{Theorem1}, it suffices to prove Proposition \ref{proposition4} which will be achieved step by step in the following lemmas.


Throughout the rest of the paper, we denote $L^p:=L^p(\mathbb{R}^3)$ and $\displaystyle\int f\,\mathrm{d}x:=\int_{\mathbb{R}^3} f\,\mathrm{d}x$, and let $C \geq 1$  represent a generic positive constant that depends on some known constants but is independent of  $\theta$, $\delta$, $\mu$, $\nu$, $t$ and $T^*$.

Although the solutions usually depend on $\mu$ and $\nu$, one can find that the following results and procedures are applicable to the case $\mu=\nu=0$. For brevity, we omit the superscripts throughout Lemmas \ref{lemma-n1}-\ref{lemma-n3}.

\begin{lemma}\label{lemma-n1}Under the same assumptions of Theorem \ref{Theorem1} and \eqref{clc9}, the following estimate
\beq\label{n00}
\begin{split}
&\frac{1}{2}\frac{d}{dt}\big(\|\rho\|_{H^3}^2
+\| u\|_{H^3}^2+\frac{r_2}{\beta\tilde{\eta}}\|\eta\|_{H^3}^2
+\frac{r_3}{2\beta
k\tilde{\eta}}\|\tau\|_{H^3}^2
-\int\frac{h(\rho)+\beta\rho}{r_1+\beta\rho}|\nabla^3\rho|^2\,\mathrm{d}x\big)\\
&+\frac{\mu_1}{2}\|\nabla u\|_{H^3}^2+\frac{\mu_2}{2}\|\mathrm{div}u\|_{H^3}^2+\frac{r_2\varepsilon}{2\beta\tilde{\eta}}\|\nabla\eta\|_{H^3}^2
+\frac{A_0 r_3}{4\lambda\beta
k\tilde{\eta}}\|\tau\|_{H^3}^2
+\frac{r_3\varepsilon}{4\beta
k\tilde{\eta}}\|\nabla\tau\|_{H^3}^2\\
\leq& C\delta(\|\nabla\rho\|_{H^2}^2+\|\nabla u\|_{H^2}^2)
\end{split}
\eeq holds for any $t\in[0,T]$.
\end{lemma}
\pf
Applying derivatives $\nabla^\ell(\ell=0,1,2,3)$ to the system \eqref{-17}, taking inner product with $\nabla^\ell\rho$, $\nabla^\ell u$, $\frac{r_2}{\beta\tilde{\eta}}\nabla^\ell \eta$ and $\frac{r_3}{2\beta k\tilde{\eta}}\nabla^\ell \tau$ respectively, and then adding the results, we can obtain
\begin{align}
&\frac{1}{2}\frac{d}{dt}\big(\|\rho\|_{H^3}^2
+\| u\|_{H^3}^2+\frac{r_2}{\beta\tilde{\eta}}\|\eta\|_{H^3}^2
+\frac{r_3}{2\beta
k\tilde{\eta}}\|\tau\|_{H^3}^2\big)\notag\\
&+\mu_1\|\nabla u\|_{H^3}^2+\mu_2\|\mathrm{div}u\|_{H^3}^2
+\frac{r_2\varepsilon}{\beta\tilde{\eta}}\|\nabla\eta\|_{H^3}^2
+\frac{A_0}{2\lambda}\frac{r_3}{2\beta
k\tilde{\eta}}\|\tau\|_{H^3}^2
+\frac{r_3\varepsilon}{2\beta
k\tilde{\eta}}\|\nabla\tau\|_{H^3}^2\label{n0}\\
=&\sum\limits_{\ell=0}^{3}\int(\nabla^\ell S_1: \nabla^\ell \rho+\nabla^\ell \tilde{S}_2: \nabla^\ell u
+\frac{r_2}{\beta\tilde{\eta}}\nabla^\ell S_3: \nabla^\ell \eta+\frac{r_3}{2\beta k\tilde{\eta}}\nabla^\ell S_4: \nabla^\ell \tau)\,\mathrm{d}x.\notag
\end{align}
Before we estimate each term on the right-hand side of \eqref{n0}, it is worth noticing that the disappearance of the viscous terms in the momentum equation \eqref{-17}$_2$ leads to partial loss of regularity of velocity $u$. Hence, to derive some uniform estimates independent
of $\mu$ and $\nu$, those terms containing the fourth derivative of density $\rho$ or velocity $u$ can not be directly controlled. In the following proof, we will list them separately and deal with them in detail.

Firstly, for the first term on the right-hand side of \eqref{n0}, by noticing the definition of $S_1$, we have
\begin{align}
&\sum\limits_{\ell=0}^{3}\int\nabla^\ell S_1: \nabla^\ell \rho\,\mathrm{d}x
=-\beta\sum\limits_{\ell=0}^{3}\int\nabla^\ell \mathrm{div}(\rho u):\nabla^\ell \rho\,\mathrm{d}x\notag\\
=&-\beta\int \mathrm{div}(\rho u) \rho\,\mathrm{d}x
-\beta\int\nabla\mathrm{div}(\rho u)\cdot\nabla \rho\,\mathrm{d}x
-\beta\sum\limits_{\ell=2}^{3}\int\nabla^\ell \mathrm{div}(\rho u):\nabla^\ell \rho\,\mathrm{d}x\label{c1}\\
:=&\sum\limits_{i=1}^{3}I_{1i}.\notag
\end{align}
The reason why we discuss $\ell$ separately here is to make the proof more concise when proving decay estimates later in this article. We first deal with the lower derivative terms, $I_{11}$ and $I_{12}$, using H\"older inequality, Sobolev inequality, Cauchy inequality and Lemma \ref{LemmaA2}, it holds that
\beq\label{c2}
|I_{11}|
\leq C\big(\|\nabla\rho\|_{L^2}\|u\|_{L^3}
+\|\rho\|_{L^3}\|\nabla u\|_{L^2}\big)\|\rho\|_{L^6}
\leq C\delta(\|\nabla \rho\|_{L^2}^2+\|\nabla u\|_{L^2}^2),
\eeq
and
\beq\label{c3}
|I_{12}|
\leq C\big(\|\nabla^2\rho\|_{L^2}\|u\|_{L^3}
+\|\rho\|_{L^3}\|\nabla^2u\|_{L^2}\big)\|\nabla\rho\|_{L^6}
\leq C\delta(\|\nabla^2\rho\|_{L^2}^2+\|\nabla^2u\|_{L^2}^2).
\eeq
Then, for $I_{13}$, which can be divided into the following five terms.
\begin{align}
I_{13}
=&-\beta\int\nabla^2 \mathrm{div}(\rho u):\nabla^2 \rho\,\mathrm{d}x
-\beta\int\nabla^3 \mathrm{div}(\rho u):\nabla^3 \rho\,\mathrm{d}x\label{b0}\\
=&-\beta\int\nabla^2 \mathrm{div}(\rho u):\nabla^2 \rho\,\mathrm{d}x
-\beta\int\nabla^3 (u\cdot\nabla\rho):\nabla^3 \rho\,\mathrm{d}x
-\beta\int\nabla^3 (\rho\,\mathrm{div}u):\nabla^3 \rho\,\mathrm{d}x\notag\\
=&-\beta\int\nabla^2 \mathrm{div}(\rho u):\nabla^2 \rho\,\mathrm{d}x
-\beta\int u\cdot\nabla\nabla^3\rho\nabla^3 \rho\,\mathrm{d}x
-\beta\sum\limits_{\ell=1}^{3}\int\mathbb{C}_3^\ell\nabla^\ell u\cdot\nabla\nabla^{3-\ell}\rho\nabla^3 \rho\,\mathrm{d}x
\notag\\
&-\beta\sum\limits_{\ell=0}^{2}\int\mathbb{C}_3^\ell\nabla^\ell\mathrm{div}u\nabla^{3-\ell} \rho:\nabla^3 \rho\,\mathrm{d}x
-\beta\int\rho\nabla^3\mathrm{div}u:\nabla^3 \rho\,\mathrm{d}x\notag\\
:=&\sum\limits_{i=1}^{5}I_{13_{i}}.\notag
\end{align}
Next, we turn to deal with the terms. In the same way,
we first deal with the lower derivative terms, $I_{13_{1}}$, $I_{13_{3}}$ and $I_{13_{4}}$, using H\"older inequality, Sobolev inequality, Cauchy inequality and Lemma \ref{LemmaA2}, it holds that
\begin{align}
|I_{13_{1}}+I_{13_{3}}+I_{13_{4}}|
\leq& C\big(\|\nabla^3\rho\|_{L^2}\|u\|_{L^3}
+\|\rho\|_{L^3}\|\nabla^3u\|_{L^2}\big)\|\nabla^2\rho\|_{L^6}\label{b1}\\
&+C\big(\|\nabla u\|_{L^\infty}\|\nabla^3\rho\|_{L^2}+\|\nabla^2 u\|_{L^3}\|\nabla^2\rho\|_{L^6}+\|\nabla^3 u\|_{L^2}\|\nabla\rho\|_{L^\infty}\big)\|\nabla^3\rho\|_{L^2}\notag\\
\leq&C\delta(\|\nabla^3\rho\|_{L^2}^2+\|\nabla^3 u\|_{L^2}^2).\notag
\end{align}
Then, for the terms containing the fourth derivative, $I_{13_{2}}$ and $I_{13_{5}}$. On one hand, for $I_{13_{2}}$, by virtue of integration by parts, H\"older inequality and Sobolev inequality, we can directly deduce
\beq\label{b2}
I_{13_{2}}
=\frac{\beta}{2}\int \mathrm{div}u\,|\nabla^3\rho|^2\,\mathrm{d}x
\leq C\|\nabla u\|_{L^\infty}\|\nabla^3\rho\|_{L^2}^2
\leq C\delta\|\nabla^3\rho\|_{L^2}^2.
\eeq
On the other hand, for $I_{13_{5}}$, by using \eqref{-17}$_1$: $\mathrm{div}u=-\frac{\rho_t+\beta u\cdot\nabla\rho}{r_1+\beta\rho}$,
we can get
\begin{align}
I_{13_{5}}
=&\beta\int\rho\nabla^3\left(\frac{\rho_t+\beta u\cdot\nabla\rho}{r_1+\beta\rho}\right):\nabla^3 \rho\,\mathrm{d}x\notag\\
=&\beta\sum\limits_{\ell=0}^{2}\int\mathbb{C}_3^\ell\rho\nabla^\ell(\rho_t+\beta u\cdot\nabla\rho)\nabla^{3-\ell}\big(\frac{1}{r_1+\beta\rho}\big):\nabla^3 \rho\,\mathrm{d}x
+\beta\int\frac{\rho}{r_1+\beta\rho}\nabla^3\rho_t:\nabla^3 \rho\,\mathrm{d}x\label{b3}\\
&+\beta^2\int\frac{\rho }{r_1+\beta\rho}\nabla^3( u\cdot\nabla\rho):\nabla^3 \rho\,\mathrm{d}x\notag\\
=&\beta\sum\limits_{\ell=0}^{2}\int\mathbb{C}_3^\ell\rho\nabla^\ell(\rho_t+\beta u\cdot\nabla\rho)\nabla^{3-\ell}\big(\frac{1}{r_1+\beta\rho}\big):\nabla^3 \rho\,\mathrm{d}x
+\beta\int\frac{\rho}{r_1+\beta\rho}\nabla^3\rho_t:\nabla^3 \rho\,\mathrm{d}x\notag\\
&+\beta^2\sum\limits_{\ell=1}^{3}\int\frac{\rho }{r_1+\beta\rho}\mathbb{C}_3^\ell\nabla^\ell u\cdot\nabla\nabla^{3-\ell}\rho:\nabla^3 \rho\,\mathrm{d}x
+\beta^2\int\frac{\rho }{r_1+\beta\rho} u\cdot\nabla\nabla^3\rho:\nabla^3 \rho\,\mathrm{d}x.\notag
\end{align}
Further, the second term and the last term on the right-hand side of \eqref{b3} equal
\beq\label{b4}
\frac{\beta}{2}\frac{d}{dt}\int\frac{\rho}{r_1+\beta\rho}|\nabla^3\rho|^2\,\mathrm{d}x
-\frac{\beta}{2}\int\left(\frac{\rho}{r_1+\beta\rho}\right)_t|\nabla^3\rho|^2\,\mathrm{d}x
-\frac{\beta^2}{2}\int\mathrm{div}\big(\frac{\rho u}{r_1+\beta\rho}\big) |\nabla^3\rho|^2\,\mathrm{d}x,
\eeq
where we use integration by parts.\\
\eqref{b3}, combined with \eqref{b4}, H\"older inequality, Sobolev inequality and Cauchy inequality, yields
\beq\label{b5}
I_{13_{5}}\leq C\delta(\|\nabla^2\rho\|_{H^1}^2+\|\nabla^3 u\|_{L^2}^2)
+\frac{\beta}{2}\frac{d}{dt}\int\frac{\rho}{r_1+\beta\rho}|\nabla^3\rho|^2\,\mathrm{d}x.
\eeq
Putting \eqref{b1}, \eqref{b2} and \eqref{b5} into \eqref{b0} yields
\beq\label{b6}
I_{13}\leq C\delta(\|\nabla^2\rho\|_{H^1}^2+\|\nabla^3 u\|_{L^2}^2)
+\frac{\beta}{2}\frac{d}{dt}\int\frac{\rho}{r_1+\beta\rho}|\nabla^3\rho|^2\,\mathrm{d}x.
\eeq
Now, substituting \eqref{c2}, \eqref{c3} and \eqref{b6} into \eqref{c1}, it holds that
\beq\label{n1}
\sum\limits_{\ell=0}^{3}\int\nabla^\ell S_1: \nabla^\ell \rho\,\mathrm{d}x
\leq C\delta(\|\nabla \rho\|_{H^2}^2+\|\nabla u\|_{H^2}^2)
+\frac{\beta}{2}\frac{d}{dt}\int\frac{\rho}{r_1+\beta\rho}|\nabla^3\rho|^2\,\mathrm{d}x.
\eeq

Secondly, the second term on the right-hand side of \eqref{n0} equals
\begin{align}
&\sum\limits_{\ell=0}^{3}\int\nabla^\ell \tilde{S}_2: \nabla^\ell u\,\mathrm{d}x\label{b00}\\
=&\int \tilde{S}_2\cdot u\,\mathrm{d}x
+\int\nabla \tilde{S}_2:\nabla u\,\mathrm{d}x
+\sum\limits_{\ell=2}^{3}\int\nabla^\ell \tilde{S}_2:\nabla^\ell u\,\mathrm{d}x\notag\\
:=&\sum\limits_{i=1}^{3}I_{2i}.\notag
\end{align}
From \eqref{-18}, we first estimate $I_{21}$ and $I_{22}$. Using H\"older inequality, Sobolev inequality and Cauchy inequality, it holds that
\begin{align}
I_{21}=&\int\big\{-\beta u\cdot\nabla u+h(\rho)\nabla \rho+g(\rho)\big[\big(k(L-1)+2\mathfrak{z}\tilde{\eta}\big)\nabla\eta-\mathrm{div}\tau\big]
-\frac{2\mathfrak{z}}{\beta(\rho+\tilde{\rho})}\eta\nabla\eta\big\}\cdot u\,\mathrm{d}x\label{c4}\\
&-\int\big(\mu\beta g(\rho)\Delta u+(\mu+\nu)\beta g(\rho)\nabla\mathrm{div} u\big)\cdot u\,\mathrm{d}x\notag\\
\leq&C\big(\|u\|_{L^3}\|\nabla u\|_{L^2}+\|h(\rho)\|_{L^3}\|\nabla\rho\|_{L^2}
+\|g(\rho)\|_{L^3}\|\nabla (\eta,\tau)\|_{L^2}+\|\frac{\eta}{\rho+\tilde{\rho}}\|_{L^3}\|\nabla \eta\|_{L^2}\big)\|u\|_{L^6}\notag\\
&+C\|g(\rho)\|_{L^\infty}\big(\mu\|\nabla u\|_{L^2}^2+(\mu+\nu)\|\mathrm{div}u\|_{L^2}^2\big)+C\|\nabla g(\rho)\|_{L^3}\big(\mu\|\nabla u\|_{L^2}+(\mu+\nu)\|\mathrm{div}u\|_{L^2}\big)\|u\|_{L^6}\notag\\
\leq&C\delta(\|\nabla \rho\|_{L^2}^2+\|\nabla u\|_{L^2}^2+\|\nabla \eta\|_{L^2}^2+\|\nabla \tau\|_{L^2}^2)\notag,
\end{align}
and
\begin{align}
I_{22}=&\int\big\{\beta u\cdot\nabla u-h(\rho)\nabla \rho-g(\rho)\big[\big(k(L-1)+2\mathfrak{z}\tilde{\eta}\big)\nabla\eta-\mathrm{div}\tau\big]
+\frac{2\mathfrak{z}}{\beta(\rho+\tilde{\rho})}\eta\nabla\eta\big\}\cdot\Delta u\,\mathrm{d}x\label{c5}\\
&+\int\big(\mu\beta g(\rho)\Delta u+(\mu+\nu)\beta g(\rho)\nabla\mathrm{div} u\big)\cdot\Delta u\,\mathrm{d}x\notag\\
\leq&C\big(\|u\|_{L^3}\|\nabla u\|_{L^6}+\|h(\rho)\|_{L^3}\|\nabla\rho\|_{L^6}
+\|g(\rho)\|_{L^3}\|\nabla (\eta,\tau)\|_{L^6}+\|\frac{\eta}{\rho+\tilde{\rho}}\|_{L^3}\|\nabla \eta\|_{L^6}\big)\|\nabla^2 u\|_{L^2}\notag\\
&+C\|g(\rho)\|_{L^\infty}\big(\mu\|\nabla^2 u\|_{L^2}+(\mu+\nu)\|\nabla\mathrm{div}u\|_{L^2}\big)\|\nabla^2u\|_{L^2}\notag\\
\leq&C\delta(\|\nabla^2\rho\|_{L^2}^2+\|\nabla^2u\|_{L^2}^2+\|\nabla^2 \eta\|_{L^2}^2+\|\nabla^2\tau\|_{L^2}^2).\notag
\end{align}
 $I_{23}$ can be divided into the following five terms:
\begin{align}
I_{23}=&\int\nabla^2\big\{-\beta u\cdot\nabla u+h(\rho)\nabla \rho+g(\rho)\big[\big(k(L-1)+2\mathfrak{z}\tilde{\eta}\big)\nabla\eta-\mathrm{div}\tau\big]
-\frac{2\mathfrak{z}}{\beta(\rho+\tilde{\rho})}\eta\nabla\eta\big\}:\nabla^2 u\,\mathrm{d}x\label{c6}\\
&-\sum\limits_{\ell=2}^{3}\int\nabla^\ell\big(\mu\beta g(\rho)\Delta u+(\mu+\nu)\beta g(\rho)\nabla\mathrm{div} u\big):\nabla^\ell u\,\mathrm{d}x\notag\\
&+\int\nabla^3\big\{g(\rho)\big[\big(k(L-1)+2\mathfrak{z}\tilde{\eta}\big)\nabla\eta-\mathrm{div}\tau\big]
-\frac{2\mathfrak{z}}{\beta(\rho+\tilde{\rho})}\eta\nabla\eta\big\}:\nabla^3 u\,\mathrm{d}x\notag\\
&+\int\nabla^3(-\beta u\cdot\nabla u):\nabla^3 u\,\mathrm{d}x
+\int\nabla^3\big(h(\rho)\nabla \rho\big):\nabla^3 u\,\mathrm{d}x\notag\\
:=&\sum\limits_{i=1}^{5}I_{23_{i}}.\notag
\end{align}
Then, we turn to deal with the terms $I_{23_{1}}$-$I_{23_{5}}$. For the terms $I_{23_{1}}$-$I_{23_{3}}$, using H\"older inequality, Sobolev inequality, Cauchy inequality and Lemma \ref{LemmaA2}, we have
\begin{align}
|I_{23_{1}}|
\leq&C\big[\|\nabla^2u\|_{L^6}\|\nabla u\|_{L^3}
+\|u\|_{L^\infty}\|\nabla^3u\|_{L^2}+\|\nabla^2h(\rho)\|_{L^6}\|\nabla \rho\|_{L^3}
+\|h(\rho)\|_{L^\infty}\|\nabla^3\rho\|_{L^2}\big]\|\nabla^2u\|_{L^2}\label{c7}\\
&+C\big[\|\nabla^2g(\rho)\|_{L^6}\|\nabla(\eta,\tau)\|_{L^3}
+\|g(\rho)\|_{L^\infty}\|\nabla^3(\eta,\tau)\|_{L^2}
\big]\|\nabla^2u\|_{L^2}\notag\\
&+C\big[\|\nabla^2(\frac{\eta }{\rho+\tilde{\rho}})\|_{L^6}\|\nabla\eta\|_{L^3}
+\|\frac{\eta }{\rho+\tilde{\rho}}\|_{L^\infty}\|\nabla^3\eta\|_{L^2}\big]\|\nabla^2u\|_{L^2}\notag\\
\leq&C\delta\big(\|\nabla^2u\|_{H^1}^2+\|\nabla^2\rho\|_{H^1}^2+\|\nabla^3\eta\|_{L^2}^2
+\|\nabla^3\tau\|_{L^2}^2\big)\notag,
\end{align}
\begin{align}
|I_{23_{2}}|
\leq&C\int\big(|\nabla^2g(\rho)||\nabla^2u|^2
+|\nabla g(\rho)||\nabla^3u||\nabla^2u|
\big)\,\mathrm{d}x
+C\int|g(\rho)|\big(\mu|\nabla^2\Delta u|
+(\mu+\nu)|\nabla^3\mathrm{div}u|\big)|\nabla^2u|\,\mathrm{d}x\notag\\
&+C\int[ |\nabla^3g(\rho)|\big(\mu|\Delta u|
+(\mu+\nu)|\nabla\mathrm{div}u|\big)|\nabla^3u| +|\nabla^2 g(\rho)|\big(\mu|\nabla\Delta u|
+(\mu+\nu)|\nabla^2\mathrm{div}u|\big)|\nabla^3u|]\,\mathrm{d}x\notag\\
&+\int C|\nabla g(\rho)|\big(\mu|\nabla^2\Delta u|
+(\mu+\nu)|\nabla^3\mathrm{div}u|\big)|\nabla^3u|-\beta g(\rho)\big(\mu\nabla^3\Delta u+(\mu+\nu)\nabla^4\mathrm{div} u\big):\nabla^3 u\,\mathrm{d}x\notag\\
\leq&C\|\nabla g(\rho)\|_{H^2}\|\nabla^2 u\|_{H^1}^2
+C\|g(\rho)\|_{H^3}\big(\mu\|\nabla^3 u\|_{H^1}
+(\mu+\nu)\|\nabla^2\mathrm{div}u\|_{H^1}\big)\|\nabla^3 u\|_{L^2}\label{c7-}\\
&+C\|g(\rho)\|_{L^\infty}\big(\mu\|\nabla^4 u\|_{L^2}^2
+(\mu+\nu)\|\nabla^3\mathrm{div}u\|_{L^2}^2\big)\notag\\
\leq&C\delta\big(\|\nabla^2 u\|_{H^1}^2+\mu\|\nabla^4 u\|_{L^2}^2
+(\mu+\nu)\|\nabla^3\mathrm{div}u\|_{L^2}^2\big),\notag
\end{align}
and
\begin{align}
|I_{23_{3}}|
\leq&C\big[\|\nabla^3g(\rho)\|_{L^2}\|\nabla(\eta,\tau)\|_{L^\infty}
+\|g(\rho)\|_{L^\infty}\|\nabla^4(\eta,\tau)\|_{L^2}
\big]\|\nabla^3u\|_{L^2}\label{c8}\\
&+C\big[\|\nabla^3(\frac{\eta }{\rho+\tilde{\rho}})\|_{L^2}\|\nabla\eta\|_{L^\infty}
+\|\frac{\eta }{\rho+\tilde{\rho}}\|_{L^\infty}\|\nabla^4\eta\|_{L^2}\big]\|\nabla^3u\|_{L^2}\notag\\
\leq&C\delta(\|\nabla^2\rho\|_{H^1}^2+\|\nabla^3 u\|_{L^2}^2+\|\nabla^3 \eta\|_{H^1}^2+\|\nabla^4 \tau\|_{L^2}^2).\notag
\end{align}
For $I_{23_{4}}$, using integration by parts, we can easily get
\begin{align}
I_{23_{4}}&=\int\nabla^3(-\beta u\cdot\nabla u):\nabla^3 u\,\mathrm{d}x\label{c9}\\
&=-\beta\sum\limits_{\ell=1}^{3}\int\mathbb{C}_3^\ell\nabla^\ell u\cdot\nabla\nabla^{3-\ell} u:\nabla^3 u\,\mathrm{d}x
+\frac{\beta}{2}\int \mathrm{div}u|\nabla^3 u|^2\,\mathrm{d}x\notag\\
&\leq C \|\nabla u\|_{L^\infty}\|\nabla^3u\|_{L^2}^2+C\|\nabla^2 u\|_{L^3}\|\nabla^2u\|_{L^6}\|\nabla^3u\|_{L^2}\notag\\
&\leq C \delta\|\nabla^3u\|_{L^2}^2,\notag
\end{align}
where we have used H\"older inequality, Sobolev inequality and Cauchy inequality.\\
Then, we are going to deal with $I_{23_{5}}$ which can be split into the following two terms.
\begin{align}
I_{23_{5}}=&\int\nabla^3\big(h(\rho)\nabla \rho\big):\nabla^3 u\,\mathrm{d}x\label{b01}\\
=&\sum\limits_{\ell=1}^{3}\int\mathbb{C}_3^\ell\nabla^\ell h(\rho)\nabla\nabla^{3-\ell} \rho:\nabla^3 u\,\mathrm{d}x
+\int h(\rho)\nabla\nabla^3 \rho:\nabla^3 u\,\mathrm{d}x.\notag
\end{align}
For the last term in \eqref{b01}, using integration by parts and \eqref{-17}$_1$ again, we have
\begin{align}
&\int h(\rho)\nabla\nabla^3 \rho:\nabla^3 u\,\mathrm{d}x\notag\\
=&-\int \nabla h(\rho)\nabla^3 \rho:\nabla^3 u\,\mathrm{d}x
-\int h(\rho)\nabla^3 \rho:\nabla^3 \mathrm{div}u\,\mathrm{d}x\label{b02}\\
=&-\int \nabla h(\rho)\nabla^3 \rho:\nabla^3 u\,\mathrm{d}x
+\int h(\rho)\nabla^3 \rho:\nabla^3\big(\frac{\rho_t+\beta u\cdot\nabla\rho}{r_1+\beta\rho}\big)\,\mathrm{d}x,\notag
\end{align} where the observation that $\mathrm{div}u=\frac{\rho_t+\beta u\cdot\nabla\rho}{r_1+\beta\rho}$ is crucial due to the loss of regularity of velocity, which can be seen for instance in \cite{Zhu}.

Further, processing method similar to $I_{13_{5}}$, the last term on the right-hand side of \eqref{b02} can be deal with like
\begin{align}
&\int h(\rho)\nabla^3\left(\frac{\rho_t+\beta u\cdot\nabla\rho}{r_1+\beta\rho}\right):\nabla^3 \rho\,\mathrm{d}x\notag\\
=&\sum\limits_{\ell=0}^{2}\int\mathbb{C}_3^\ell h(\rho)\nabla^\ell(\rho_t+\beta u\cdot\nabla\rho)\nabla^{3-\ell}\big(\frac{1}{r_1+\beta\rho}\big):\nabla^3 \rho\,\mathrm{d}x\label{b03}\\
&+\beta\sum\limits_{\ell=1}^{3}\int\frac{h(\rho) }{r_1+\beta\rho}\mathbb{C}_3^\ell\nabla^\ell u\cdot\nabla\nabla^{3-\ell}\rho:\nabla^3 \rho\,\mathrm{d}x
-\frac{\beta}{2}\int\mathrm{div}\big(\frac{h(\rho) u}{r_1+\beta\rho}\big) |\nabla^3\rho|^2\,\mathrm{d}x\notag\\
&+\frac{1}{2}\frac{d}{dt}\int\frac{h(\rho)}{r_1+\beta\rho}|\nabla^3\rho|^2\,\mathrm{d}x
-\frac{1}{2}\int\left(\frac{h(\rho)}{r_1+\beta\rho}\right)_t|\nabla^3\rho|^2\,\mathrm{d}x.\notag
\end{align}
Together with \eqref{b01}, \eqref{b02} and \eqref{b03}, using H\"older inequality, Sobolev inequality and Cauchy inequality, we get
\beq\label{c10}
|I_{23_{5}}|\leq C\delta(\|\nabla^2\rho\|_{H^1}^2+\|\nabla^3 u\|_{L^2}^2)+\frac{1}{2}\frac{d}{dt}\int\frac{h(\rho)}{r_1+\beta\rho}|\nabla^3\rho|^2\,\mathrm{d}x.
\eeq
Hence, substituting \eqref{c7}, \eqref{c7-}, \eqref{c8}, \eqref{c9} and \eqref{c10} into \eqref{c6}, it holds that
\begin{align}
I_{23}\leq& C\delta\big(\|\nabla^2u\|_{H^1}^2+\|\nabla^2\rho\|_{H^1}^2+\|\nabla^3\eta\|_{H^1}^2
+\|\nabla^3\tau\|_{H^1}^2\big)
+\frac{1}{2}\frac{d}{dt}\int\frac{h(\rho)}{r_1+\beta\rho}|\nabla^3\rho|^2\,\mathrm{d}x\label{b7}\\
&+ C\delta\big(\mu\|\nabla^4 u\|_{L^2}^2
+(\mu+\nu)\|\nabla^3\mathrm{div}u\|_{L^2}^2\big)\notag.
\end{align}
Now, putting \eqref{c4}, \eqref{c5} and \eqref{b7} into \eqref{b00}, we have
\begin{align}
\sum\limits_{\ell=0}^{3}\int\nabla^\ell \tilde{S}_2: \nabla^\ell u\,\mathrm{d}x
\leq& C\delta\big(\|\nabla u\|_{H^2}^2+\|\nabla\rho\|_{H^2}^2+\|\nabla\eta\|_{H^3}^2
+\|\nabla\tau\|_{H^3}^2\big)
+\frac{1}{2}\frac{d}{dt}\int\frac{h(\rho)}{r_1+\beta\rho}|\nabla^3\rho|^2\,\mathrm{d}x\label{n2}\\
&+ C\delta\big(\mu\|\nabla^4 u\|_{L^2}^2
+(\mu+\nu)\|\nabla^3\mathrm{div}u\|_{L^2}^2\big)\notag.
\end{align}

Next, for the third term on the right-hand side of \eqref{n0}, using \eqref{18}$_3$, we have
\begin{align*}
&\sum\limits_{\ell=0}^{3}\int\nabla^\ell S_3: \nabla^\ell \eta\,\mathrm{d}x=-\beta\sum\limits_{\ell=0}^{3}\int\nabla^\ell \mathrm{div}(\eta u): \nabla^\ell \eta\,\mathrm{d}x\\
=&-\beta\int\mathrm{div}(\eta u)\,\eta\,\mathrm{d}x
-\beta\int\nabla \mathrm{div}(\eta u): \nabla \eta\,\mathrm{d}x
-\beta\sum\limits_{\ell=2}^{3}\int\nabla^\ell \mathrm{div}(\eta u): \nabla^\ell \eta\,\mathrm{d}x\\
:=&\sum\limits_{i=1}^{3}I_{3i}.
\end{align*}
For $I_{31}$ and $I_{32}$, similar to $I_{11}$ and $I_{12}$, it is not hard to get
\beq\label{c11}
|I_{31}|
\leq C\big(\|\nabla\eta\|_{L^2}\|u\|_{L^3}
+\|\eta\|_{L^3}\|\nabla u\|_{L^2}\big)\|\eta\|_{L^6}
\leq C\delta(\|\nabla \eta\|_{L^2}^2+\|\nabla u\|_{L^2}^2),
\eeq
and
\beq\label{c12}
|I_{32}|
\leq C\big(\|\nabla^2\eta\|_{L^2}\|u\|_{L^3}
+\|\eta\|_{L^3}\|\nabla^2u\|_{L^2}\big)\|\nabla\eta\|_{L^6}
\leq C\delta(\|\nabla^2\eta\|_{L^2}^2+\|\nabla^2u\|_{L^2}^2).
\eeq
Thanks to H\"older inequality, Sobolev inequality, Cauchy inequality and Lemma \ref{LemmaA2}, $I_{33}$ can be controlled like
\beq\label{c13}
I_{33}=\beta\sum\limits_{\ell=2}^{3}\int\nabla^{\ell-1} \mathrm{div}(\eta u): \nabla^{\ell+1} \eta\,\mathrm{d}x
\leq C\delta(\|\nabla^3\eta\|_{H^1}^2+\|\nabla^3u\|_{L^2}^2).
\eeq
Note that we have used integration by parts in \eqref{c13} to reduce the order of spatial derivative of velocity $u$.\\
Hence, combining with \eqref{c11}, \eqref{c12} and \eqref{c13} yields
\beq\label{n3}
\sum\limits_{\ell=0}^{3}\int\nabla^\ell S_3: \nabla^\ell \eta\,\mathrm{d}x
\leq C\delta(\|\nabla\eta\|_{H^3}^2+\|\nabla u\|_{H^2}^2).
\eeq

Finally, for the last term on the right-hand side of \eqref{n0}, we have
\begin{align*}
\sum\limits_{\ell=0}^{3}\int\nabla^\ell S_4: \nabla^\ell \tau\,\mathrm{d}x
=&\int S_4:\tau\,\mathrm{d}x+\int\nabla S_4: \nabla \tau\,\mathrm{d}x
+\sum\limits_{\ell=2}^{3}\int\nabla^\ell S_4: \nabla^\ell \tau\,\mathrm{d}x\\
:=&\sum\limits_{i=1}^{3}I_{4i}.
\end{align*}
Thanks to \eqref{18}$_4$, and using H\"older inequality, Sobolev inequality and Cauchy inequality again, $I_{41}$ and $I_{42}$ can be controlled like
\begin{align}
I_{41}=&\beta\int \big(-\mathrm{div} (u\tau)+(\nabla u
\tau+\tau\nabla^Tu)+k\eta(\nabla u+\nabla^Tu)\big): \tau\,\mathrm{d}x\notag\\
\leq&C\big(\|\nabla u\|_{L^2}\|(\eta,\tau)\|_{L^3}
+\|\nabla \tau\|_{L^2}\| u\|_{L^3}\big)\|\tau\|_{L^6}\label{c14}\\
\leq& C\delta(\|\nabla \tau\|_{L^2}^2+\|\nabla u\|_{L^2}^2),\notag
\end{align}
and
\begin{align}
I_{42}=&\beta\int\nabla \big(-\mathrm{div} (u\tau)+(\nabla u
\tau+\tau\nabla^Tu)+k\eta(\nabla u+\nabla^Tu)\big): \nabla \tau\,\mathrm{d}x\notag\\
=&\beta\int \big(\mathrm{div} (u\tau)-(\nabla u
\tau+\tau\nabla^Tu)-k\eta(\nabla u+\nabla^Tu)\big): \Delta \tau\,\mathrm{d}x\label{c15}\\
\leq&C\big(\|\nabla u\|_{L^6}\|(\eta,\tau)\|_{L^3}
+\|\nabla\tau\|_{L^6}\|u\|_{L^3}\big)\|\nabla^2\tau\|_{L^2}\notag\\
\leq& C\delta(\|\nabla^2\tau\|_{L^2}^2+\|\nabla^2u\|_{L^2}^2).\notag
\end{align}
For $I_{43}$, we will take advantage of the higher integrability of $\tau$ to deal with each term on the right-hand side through integration by parts. In other words, thanks to \eqref{18}$_4$ again, $I_{43}$ can be estimated as
\begin{align}
I_{43}&=\sum\limits_{\ell=2}^{3}\int\nabla^\ell S_4: \nabla^\ell \tau\,\mathrm{d}x
=-\sum\limits_{\ell=2}^{3}\int\nabla^{\ell-1} S_4: \nabla^{\ell-1}\Delta \tau\,\mathrm{d}x\notag\\
&=\beta\sum\limits_{\ell=2}^{3}\int\nabla^{\ell-1} \big(\mathrm{div} (u\tau)-(\nabla u
\tau+\tau\nabla^Tu)-k\eta(\nabla u+\nabla^Tu)\big): \nabla^{\ell-1}\Delta \tau\,\mathrm{d}x\label{c16}\\
\leq&C\big(\|\nabla^2u\|_{L^6}\|(\eta, \tau)\|_{L^3}+\|\nabla^2\tau\|_{L^6}\|u\|_{L^3}
+\|\nabla u\|_{L^6}\|\nabla(\eta,\tau)\|_{L^3}\big)\|\nabla^3\tau\|_{L^2}\notag\\
&+C\big(\|\nabla^3u\|_{L^2}\|(\eta, \tau)\|_{L^\infty}+\|\nabla^3\tau\|_{L^2}\|u\|_{L^\infty}
+\|\nabla u\|_{L^6}\|\nabla^2(\eta,\tau)\|_{L^3}\big)\|\nabla^4\tau\|_{L^2}\notag\\
\leq&C\delta(\|\nabla^2u\|_{H^1}^2+\|\nabla^3\tau\|_{H^1}^2),\notag
\end{align}
where we have used H\"older inequality, Sobolev inequality, Cauchy inequality and Lemma \ref{LemmaA2}.\\
Owing to \eqref{c14}, \eqref{c15} and \eqref{c16}, we get
\beq\label{n4}
\sum\limits_{\ell=0}^{3}\int\nabla^\ell S_4: \nabla^\ell \tau\,\mathrm{d}x
\leq C\delta(\|\nabla u\|_{H^2}^2+\|\nabla\tau\|_{H^3}^2).
\eeq

Finally, plugging \eqref{n1}, \eqref{n2}, \eqref{n3} and \eqref{n4} into \eqref{n0}, we then obtain the following inequality:
\begin{align*}
&\frac{1}{2}\frac{d}{dt}\big(\|\rho\|_{H^3}^2
+\| u\|_{H^3}^2+\frac{r_2}{\beta\tilde{\eta}}\|\eta\|_{H^3}^2
+\frac{r_3}{2\beta
k\tilde{\eta}}\|\tau\|_{H^3}^2
-\int\frac{h(\rho)+\beta\rho}{r_1+\beta\rho}|\nabla^3\rho|^2\,\mathrm{d}x\big)\\
&+\mu_1\|\nabla u\|_{H^3}^2+\mu_2\|\mathrm{div}u\|_{H^3}^2
+\frac{r_2\varepsilon}{\beta\tilde{\eta}}\|\nabla\eta\|_{H^3}^2
+\frac{A_0}{2\lambda}\frac{r_3}{2\beta
k\tilde{\eta}}\|\tau\|_{H^3}^2
+\frac{r_3\varepsilon}{2\beta
k\tilde{\eta}}\|\nabla\tau\|_{H^3}^2\\
\leq& C\delta\big(\|\nabla\rho\|_{H^2}^2+\|\nabla u\|_{H^2}^2+\|\nabla\eta\|_{H^3}^2+\|\nabla\tau\|_{H^3}^2+\mu\|\nabla^4 u\|_{L^2}^2
+(\mu+\nu)\|\nabla^3\mathrm{div}u\|_{L^2}^2\big).
\end{align*}
Choosing $\delta$ sufficiently small in the above inequality, $(\ref{n00})$ will be established. Thus, we complete the proof of this lemma.
\endpf

\medskip

In the following lemmas, we obtain some dissipation estimates of velocity $u$ and density $\rho$ which are independent of the viscosity coefficients.
\begin{lemma}\label{lemma-n2}Under the same assumptions of Theorem \ref{Theorem1} and \eqref{clc9}, we then have the following estimate:
\begin{align}
&\frac{d}{dt}\sum\limits_{\ell=0}^{2}\int\nabla^\ell\mathrm{div}\tau:\nabla^\ell u\,\mathrm{d}x+\frac{\beta k \tilde{\eta}}{2}(\|\nabla u\|_{H^2}^2+\|\mathrm{div}u\|_{H^2}^2)\label{n6}\\
\leq&C(\epsilon+\delta)\|\nabla \rho\|_{H^2}^2+C(\epsilon+\delta)\|\nabla \eta\|_{H^2}^2+C_\epsilon\|\nabla\tau\|_{H^3}^2
+\frac{A_0}{2\lambda}C\|\tau\|_{H^2}^2.\notag
\end{align}
\end{lemma}
\pf
Let $\ell= 0,1,2$. Applying operator $\nabla^\ell\mathrm{div}$ to the equation of \eqref{-17}$_4$ and $\nabla^\ell$ to the equation of \eqref{-17}$_2$, multiplying the results by $\nabla^\ell u$ and $\nabla^\ell\mathrm{div}\tau$, respectively, summing them up and then integrating it over $\mathbb{R}^3$, we get the equality:
\begin{align}
&\frac{d}{dt}\sum\limits_{\ell=0}^{2}\int\nabla^\ell\mathrm{div}\tau:\nabla^\ell u\,\mathrm{d}x+\beta k \tilde{\eta}(\|\nabla u\|_{H^2}^2+\|\mathrm{div}u\|_{H^2}^2)\notag\\
=&\sum\limits_{\ell=0}^{2}\int(\nabla^\ell\mathrm{div}S_4-\frac{A_0}{2\lambda}\nabla^\ell\mathrm{div}\tau
+\varepsilon\nabla^\ell\Delta\mathrm{div}\tau):\nabla^\ell u\,\mathrm{d}x\label{n5}\\
&+\sum\limits_{\ell=0}^{2}\int(\nabla^\ell \tilde{S}_2-r_1\nabla^{\ell+1}\rho-r_2\nabla^{\ell+1}\eta
+r_3\nabla^\ell\mathrm{div}\tau
+\mu_1\nabla^\ell\Delta u+\mu_2\nabla^{\ell+1} \mathrm{div}u):\nabla^\ell\mathrm{div}\tau\,\mathrm{d}x.\notag
\end{align}
First, by the definition of $S_4$ and $\tilde{S}_2$, integration by parts, H\"older inequality, Sobolev inequality, Cauchy inequality and Lemma \ref{LemmaA2}, it holds that
\begin{align}
&\sum\limits_{\ell=0}^{2}\int\nabla^\ell\mathrm{div}S_4:\nabla^\ell u\,\mathrm{d}x=-\sum\limits_{\ell=0}^{2}\int\nabla^\ell S_4:\nabla^{\ell+1} u\,\mathrm{d}x\notag\\
=&\beta\sum\limits_{\ell=0}^{2}\int\nabla^\ell\big(\mathrm{div} (u\tau)-(\nabla u
\tau+\tau\nabla^Tu)-k\eta(\nabla u+\nabla^Tu)\big):\nabla^{\ell+1} u\,\mathrm{d}x\label{b11}\\
\leq&C\big(\|\nabla u\|_{L^2}\|(\eta,\tau)\|_{L^\infty}+\|\nabla \tau\|_{L^2}\|u\|_{L^\infty}\big)\|\nabla u\|_{L^2}\notag\\
&+C\big(\|\nabla^2 u\|_{L^2}\|(\eta,\tau)\|_{L^\infty}+\|\nabla^2 \tau\|_{L^2}\|u\|_{L^\infty}+\|\nabla u\|_{L^6}\|\nabla(\eta,\tau)\|_{L^3}\big)\|\nabla^2 u\|_{L^2}\notag\\
&+C\big(\|\nabla^3 u\|_{L^2}\|(\eta,\tau)\|_{L^\infty}+\|\nabla^3 \tau\|_{L^2}\|u\|_{L^\infty}+\|\nabla u\|_{L^6}\|\nabla^2(\eta,\tau)\|_{L^3}\big)\|\nabla^3 u\|_{L^2}\notag\\
\leq&C\delta(\|\nabla u\|_{H^2}^2+\|\nabla\tau\|_{H^2}^2)\notag
\end{align}
and
\begin{align*}
\sum\limits_{\ell=0}^{2}\int\nabla^\ell \tilde{S}_2:\nabla^\ell\mathrm{div}\tau\,\mathrm{d}x
=\int \tilde{S}_2\cdot\mathrm{div}\tau\,\mathrm{d}x
-\int \tilde{S}_2\cdot\Delta\mathrm{div}\tau\,\mathrm{d}x
+\int\nabla^2 \tilde{S}_2:\nabla^2\mathrm{div}\tau\,\mathrm{d}x,
\end{align*}
where the terms on the right-hand side of the above equality can be controlled as
\begin{align}
&\int \tilde{S}_2\cdot\mathrm{div}\tau\,\mathrm{d}x\label{b12}\\
=&\int\left(-\beta u\cdot\nabla u+h(\rho)\nabla \rho+g(\rho)\big[\big(k(L-1)+2\mathfrak{z}\tilde{\eta}\big)\nabla\eta-\mathrm{div}\tau\big]
-\frac{2\mathfrak{z}}{\beta(\rho+\tilde{\rho})}\eta\nabla\eta\right) \cdot\mathrm{div}\tau\,\mathrm{d}x\notag\\
&-\int\big(\mu\beta g(\rho)\Delta u+(\mu+\nu)\beta g(\rho)\nabla\mathrm{div} u\big)\cdot\mathrm{div}\tau\,\mathrm{d}x\notag\\
\leq&C\big[\|u\|_{L^\infty}\|\nabla u\|_{L^2}+\|h(\rho)\|_{L^\infty}\|\nabla\rho\|_{L^2}
+\|g(\rho)\|_{L^\infty}\|\nabla(\eta,\tau)\|_{L^2}
+\|\frac{\eta}{\rho+\tilde{\rho}}\|_{L^\infty}\|\nabla \eta\|_{L^2}\big]\|\nabla\tau\|_{L^2}\notag\\
&+C\big(\mu\|\nabla u\|_{L^2}+(\mu+\nu)\|\mathrm{div}u\|_{L^2}\big)
\big(\|g(\rho)\|_{L^\infty}\|\nabla\mathrm{div}\tau\|_{L^2}
+\|\nabla g(\rho)\|_{L^3}\|\mathrm{div}\tau\|_{L^6}\big)\notag\\
\leq&C\delta(\|\nabla\rho\|_{L^2}^2+\|\nabla u\|_{L^2}^2+\|\nabla\eta\|_{L^2}^2+\|\nabla\tau\|_{L^2}^2
+\|\nabla^2\tau\|_{L^2}^2)\notag,
\end{align}
\begin{align}
&-\int \tilde{S}_2\cdot\Delta\mathrm{div}\tau\,\mathrm{d}x\label{c17}\\
=&\int\left(\beta u\cdot\nabla u-h(\rho)\nabla \rho-g(\rho)\big[\big(k(L-1)+2\mathfrak{z}\tilde{\eta}\big)\nabla\eta-\mathrm{div}\tau\big]
+\frac{2\mathfrak{z}}{\beta(\rho+\tilde{\rho})}\eta\nabla\eta\right) \cdot\Delta\mathrm{div}\tau\,\mathrm{d}x\notag\\
&+\int\big(\mu\beta g(\rho)\Delta u+(\mu+\nu)\beta g(\rho)\nabla\mathrm{div} u\big)\cdot\Delta\mathrm{div}\tau\,\mathrm{d}x\notag\\
\leq&C\big[\|u\|_{L^3}\|\nabla u\|_{L^6}+\|h(\rho)\|_{L^3}\|\nabla\rho\|_{L^6}
+\|g(\rho)\|_{L^3}\|\nabla(\eta,\tau)\|_{L^6}
+\|\frac{\eta}{\rho+\tilde{\rho}}\|_{L^3}\|\nabla \eta\|_{L^6}\big]\|\nabla^3\tau\|_{L^2}\notag\\
&+C\|g(\rho)\|_{L^\infty}\big(\mu\|\Delta u\|_{L^2}+(\mu+\nu)\|\nabla\mathrm{div}u\|_{L^2}\big)\|\Delta\mathrm{div}\tau\|_{L^2}
\notag\\
\leq&C\delta(\|\nabla^2\rho\|_{L^2}^2+\|\nabla^2 u\|_{L^2}^2+\|\nabla^2\eta\|_{L^2}^2+\|\nabla^2\tau\|_{H^1}^2)\notag,
\end{align}
and
\begin{align}
&\int\nabla^2 \tilde{S}_2:\nabla^2\mathrm{div}\tau\,\mathrm{d}x\label{c18}\\
=&\int\nabla^2\left(-\beta u\cdot\nabla u+h(\rho)\nabla \rho+g(\rho)\big[\big(k(L-1)+2\mathfrak{z}\tilde{\eta}\big)\nabla\eta-\mathrm{div}\tau\big]
-\frac{2\mathfrak{z}}{\beta(\rho+\tilde{\rho})}\eta\nabla\eta\right) :\nabla^2\mathrm{div}\tau\,\mathrm{d}x\notag\\
&-\int\nabla^2\big(\mu\beta g(\rho)\Delta u+(\mu+\nu)\beta g(\rho)\nabla\mathrm{div} u\big)\cdot\nabla^2\mathrm{div}\tau\,\mathrm{d}x\notag\\
\leq&C\big[\|u\|_{L^\infty}\|\nabla^3 u\|_{L^2}+\|\nabla^2u\|_{L^6}\|\nabla u\|_{L^3}
+\|h(\rho)\|_{L^\infty}\|\nabla^3\rho\|_{L^2}+\|\nabla^2h(\rho)\|_{L^6}\|\nabla\rho\|_{L^3}
\big]\|\nabla^3\tau\|_{L^2}\notag\\
&+C\big[\|g(\rho)\|_{L^\infty}\|\nabla^3(\eta,\tau)\|_{L^2}
+\|\nabla^2g(\rho)\|_{L^6}\|\nabla(\eta,\tau)\|_{L^3}
+\|\frac{\eta}{\rho+\tilde{\rho}}\|_{L^\infty}\|\nabla^3 \eta\|_{L^2}\big]\|\nabla^3\tau\|_{L^2}\notag\\
&+C\|\nabla^2\big(\frac{\eta}{\rho+\tilde{\rho}}\big)\|_{L^6}\|\nabla \eta\|_{L^3}\|\nabla^3\tau\|_{L^2}
+C\|g(\rho)\|_{L^\infty}\big(\mu\|\nabla^3 u\|_{L^2}+(\mu+\nu)\|\nabla^2\mathrm{div}u\|_{L^2}\big)\|\nabla^4\tau\|_{L^2}
\notag\\
&+C\|\nabla g(\rho)\|_{L^3}\big(\mu\|\nabla^2 u\|_{L^6}+(\mu+\nu)\|\nabla\mathrm{div}u\|_{L^6}\big)\|\nabla^4\tau\|_{L^2}
\notag\\
\leq&C\delta(\|\nabla^2\rho\|_{H^1}^2+\|\nabla^3 u\|_{L^2}^2+\|\nabla^3\eta\|_{L^2}^2+\|\nabla^3\tau\|_{H^1}^2)\notag.
\end{align}
For the remaining two terms on the second line of \eqref{n5}, by using integration by parts, H\"oder inequality and Cauchy inequality, it is straightforward to show that
\begin{align}
&\sum\limits_{\ell=0}^{2}\int\big(-\frac{A_0}{2\lambda}\nabla^\ell\mathrm{div}\tau
+\varepsilon\nabla^\ell\Delta\mathrm{div}\tau\big):\nabla^\ell u\,\mathrm{d}x\notag\\
=&\sum\limits_{\ell=0}^{2}\int\big(\frac{A_0}{2\lambda}\nabla^\ell\tau
-\varepsilon\nabla^\ell\Delta\tau\big):\nabla^{\ell+1} u\,\mathrm{d}x\label{b13}\\
\leq&C \|\nabla u\|_{H^2}(\|\nabla^2\tau\|_{H^2}+\frac{A_0}{2\lambda}\|\tau\|_{H^2})\notag\\
\leq& \frac{\beta k \tilde{\eta}}{4} \|\nabla u\|_{H^2}^2+C(\|\nabla^2\tau\|_{H^2}^2+\frac{A_0}{2\lambda}\|\tau\|_{H^2}^2).\notag
\end{align}
Finally, the remaining terms in \eqref{n5} can be estimated by using  integration by parts, H\"older inequality and Cauchy inequality as follows
\begin{align}
&\sum\limits_{\ell=0}^{2}\int(-r_1\nabla^{\ell+1}\rho-r_2\nabla^{\ell+1}\eta
+r_3\nabla^\ell\mathrm{div}\tau+\mu_1\nabla^\ell\Delta u+\mu_2\nabla^{\ell+1} \mathrm{div}u):\nabla^\ell\mathrm{div}\tau\,\mathrm{d}x\label{b14}\\
\leq&C(\|\nabla \rho\|_{H^2}+\|\nabla \eta\|_{H^2}+\|\nabla\tau\|_{H^2})\|\nabla\tau\|_{H^2}
+\epsilon\mu_1\|\nabla u\|_{H^2}^2+\epsilon\mu_2\|\mathrm{div} u\|_{H^2}^2
+C_\epsilon\|\nabla\mathrm{div}\tau\|_{H^2}^2
\notag\\
\leq&C\epsilon\|\nabla \rho\|_{H^2}^2+C\epsilon\|\nabla \eta\|_{H^2}^2
+C\epsilon\|\nabla u\|_{H^2}^2+C\epsilon\|\mathrm{div} u\|_{H^2}^2
+C_\epsilon\|\nabla\tau\|_{H^3}^2.\notag
\end{align}
Plugging \eqref{b11}-\eqref{b14} into \eqref{n5}, then choosing $\delta$ sufficiently small yields
\begin{align*}
&\frac{d}{dt}\sum\limits_{\ell=0}^{2}\int\nabla^\ell\mathrm{div}\tau:\nabla^\ell u\,\mathrm{d}x+\frac{\beta k \tilde{\eta}}{2}(\|\nabla u\|_{H^2}^2+\|\mathrm{div}u\|_{H^2}^2)\\
\leq&C(\epsilon+\delta)\|\nabla \rho\|_{H^2}^2+C(\epsilon+\delta)\|\nabla \eta\|_{H^2}^2+C_\epsilon\|\nabla\tau\|_{H^3}^2
+\frac{A_0}{2\lambda}C\|\tau\|_{H^2}^2,
\end{align*}
which is \eqref{n6}. Thus, the proof of Lemma \ref{lemma-n2} is complete.
\endpf

\begin{lemma}\label{lemma-n3} Under the same assumptions of Theorem \ref{Theorem1} and \eqref{clc9}, we then have the following estimate:
\begin{align}
&\frac{d}{dt}\sum\limits_{\ell=0}^{2}\int\nabla^\ell u:\nabla^{\ell+1} \rho\,\mathrm{d}x+\frac{r_1}{2}\|\nabla\rho\|_{H^2}^2\label{n7}\\
\leq& C(\|\nabla u\|_{H^2}^{2}+\|\nabla\eta\|_{H^2}^{2}+\|\nabla \tau\|_{H^2}^{2}+\mu_1\|\nabla^4 u\|_{L^2}^2
+\mu_2\|\nabla^3\mathrm{div}u\|_{L^2}^2).\notag
\end{align}
\end{lemma}
\pf
Let $\ell= 0,1,2$. Applying operator $\nabla^\ell$ to the equation of \eqref{-17}$_2$ and \eqref{-17}$_1$, multiplying the results by $\nabla^{\ell+1} \rho$ and $-\nabla^\ell\mathrm{div}u$ respectively, then summing them up and integrating it over $\mathbb{R}^3$, we get the equality:
\begin{align}
&\frac{d}{dt}\sum\limits_{\ell=0}^{2}\int\nabla^\ell u:\nabla^{\ell+1} \rho\,\mathrm{d}x+r_1\|\nabla\rho\|_{H^2}^2\notag\\
=&\sum\limits_{\ell=0}^{2}\int(\nabla^\ell \tilde{S}_2-r_2\nabla^{\ell+1}\eta
+r_3\nabla^\ell\mathrm{div}\tau+\mu_1\nabla^\ell\Delta u+\mu_2\nabla^{\ell+1} \mathrm{div}u):\nabla^{\ell+1} \rho\,\mathrm{d}x\label{n8}\\
&+\sum\limits_{\ell=0}^{2}\int(-\nabla^\ell S_1
+r_1\nabla^\ell\mathrm{div}u):\nabla^\ell \mathrm{div}u\,\mathrm{d}x.\notag
\end{align}
Likewise for \eqref{b12}-\eqref{c18}, from \eqref{-18}, H\"oder inequality, Sobolev inequality, Cauchy inequality and Lemma \ref{LemmaA2}, the first line on the right-hand side of \eqref{n8} can be controlled like
\begin{align}
&\sum\limits_{\ell=0}^{2}\int\nabla^\ell \tilde{S}_2:\nabla^{\ell+1} \rho\,\mathrm{d}x\label{b15}\\
=&\sum\limits_{\ell=0}^{2}\int\nabla^\ell\left(-\beta u\cdot\nabla u+h(\rho)\nabla \rho+g(\rho)\big[\big(k(L-1)+2\mathfrak{z}\tilde{\eta}\big)\nabla\eta-\mathrm{div}\tau\big]
-\frac{2\mathfrak{z}}{\beta(\rho+\tilde{\rho})}\eta\nabla\eta\right) :\nabla^{\ell+1} \rho\,\mathrm{d}x\notag\\
&-\sum\limits_{\ell=0}^{2}\int\nabla^\ell\left(\mu\beta g(\rho)\Delta u+(\mu+\nu)\beta g(\rho)\nabla\mathrm{div} u\right) :\nabla^{\ell+1} \rho\,\mathrm{d}x\notag\\
\leq&C\delta\left(\|\nabla u\|_{H^2}^2+\|\nabla \rho\|_{H^2}^2+\|\nabla\eta\|_{H^2}^2+\|\nabla \tau\|_{H^2}^2+\mu\|\nabla^4 u\|_{L^2}^2
+(\mu+\nu)\|\nabla^3\mathrm{div}u\|_{L^2}^2\right),\notag
\end{align}
and
\begin{align}
&\sum\limits_{\ell=0}^{2}\int(-r_2\nabla^{\ell+1}\eta
+r_3\nabla^\ell\mathrm{div}\tau+\mu_1\nabla^\ell\Delta u+\mu_2\nabla^{\ell+1} \mathrm{div}u):\nabla^{\ell+1} \rho\,\mathrm{d}x\label{b16}\\
\leq& C\left(\|\nabla\eta\|_{H^2}+\|\nabla \tau\|_{H^2}+\mu_1\|\nabla^2u\|_{H^2}+\mu_2\|\nabla\mathrm{div}u\|_{H^2}\right)\|\nabla\rho\|_{H^{2}}\notag\\
\leq&\frac{r_1}{4}\|\nabla\rho\|_{H^{2}}^{2}
+C\big(\|\nabla\eta\|_{H^2}^{2}+\|\nabla \tau\|_{H^2}^{2}
+\mu_1\|\nabla^2u\|_{H^2}^2+\mu_2\|\nabla\mathrm{div}u\|_{H^2}^2\big).\notag
\end{align}
Similarly, for the last line in \eqref{n8}, recalling the definition of $S_1$, using H\"oder inequality, Sobolev inequality, Cauchy inequality and Lemma \ref{LemmaA2} again, we have
\begin{align}
&\sum\limits_{\ell=0}^{2}\int(-\nabla^\ell S_1
+r_1\nabla^\ell\mathrm{div}u):\nabla^\ell \mathrm{div}u\,\mathrm{d}x\notag\\
=&\sum\limits_{\ell=0}^{2}\int(\beta\nabla^\ell\mathrm{div}(\rho u)
+r_1\nabla^\ell\mathrm{div}u):\nabla^\ell\mathrm{div} u\,\mathrm{d}x\label{b17}\\
\leq&C(\|\nabla u\|_{H^{2}}\|\rho\|_{L^\infty}+\|\nabla \rho\|_{H^{2}}\|u\|_{L^\infty})\|\nabla u\|_{H^{2}}
+C\|\nabla u\|_{H^{2}}^2\notag\\
\leq&C\delta\|\nabla\rho\|_{H^{2}}^{2}
+C\|\nabla u\|_{H^2}^{2}.\notag
\end{align}

Now, plugging \eqref{b15}-\eqref{b17} into \eqref{n8}, then choosing $\delta$ sufficiently small yields
$$
\frac{d}{dt}\sum\limits_{\ell=0}^{2}\int\nabla^\ell u:\nabla^{\ell+1} \rho\,\mathrm{d}x+\frac{r_1}{2}\|\nabla\rho\|_{H^2}^2
\leq C(\|\nabla u\|_{H^2}^{2}+\|\nabla\eta\|_{H^2}^{2}+\|\nabla \tau\|_{H^2}^{2}+\mu_1\|\nabla^4 u\|_{L^2}^2
+\mu_2\|\nabla^3\mathrm{div}u\|_{L^2}^2),
$$
which is \eqref{n7}. Thus, we finish the proof of Lemma \ref{lemma-n3}.
\endpf

\bigskip

In what follows, based on Lemmas \ref{lemma-n1}-\ref{lemma-n3}, we are ready to prove Proposition \ref{proposition4}.
\paragraph{Proof of Proposition \ref{proposition4}:}

\medskip

Combined with \eqref{n00}, $\epsilon_2$\eqref{n6} and $\epsilon_1$\eqref{n7}, it holds that
\begin{align*}
&\frac{1}{2}\frac{d}{dt}\big(\|\rho\|_{H^3}^2
+\| u\|_{H^3}^2+\frac{r_2}{\beta\tilde{\eta}}\|\eta\|_{H^3}^2
+\frac{r_3}{2\beta
k\tilde{\eta}}\|\tau\|_{H^3}^2\big)\\
&+\frac{d}{dt}\int\big(\epsilon_1\sum\limits_{\ell=0}^{2}\nabla^\ell u:\nabla^{\ell+1} \rho+\epsilon_2\sum\limits_{\ell=0}^{2}\nabla^\ell\mathrm{div}\tau:\nabla^\ell u-\frac{1}{2}\frac{h(\rho)+\beta\rho}{r_1+\beta\rho}|\nabla^3\rho|^2\big)\,\mathrm{d}x
+\frac{\mu_1}{2}\|\nabla u\|_{H^3}^2+\frac{\mu_2}{2}\|\mathrm{div}u\|_{H^3}^2\\
&+\epsilon_1\frac{r_1}{2}\|\nabla\rho\|_{H^2}^2
+\epsilon_2\frac{\beta k \tilde{\eta}}{2}\|\nabla u\|_{H^2}^2
+\frac{r_2\varepsilon}{2\beta\tilde{\eta}}\|\nabla\eta\|_{H^3}^2
+\frac{A_0 r_3}{4\lambda\beta
k\tilde{\eta}}\|\tau\|_{H^3}^2
+\frac{r_3\varepsilon}{4\beta
k\tilde{\eta}}\|\nabla\tau\|_{H^3}^2\\
\leq&C(\delta+\epsilon_2\epsilon)\|\nabla\rho\|_{H^2}^2
+C(\delta+\epsilon_1)\|\nabla u\|_{H^2}^2+C(\epsilon_2\epsilon+\epsilon_1+\delta)\|\nabla \eta\|_{H^2}^2+C_\epsilon(\epsilon_2+\epsilon_1)\|\nabla\tau\|_{H^3}^2
+\epsilon_2\frac{A_0}{2\lambda}C\|\tau\|_{H^2}^2\\
&+\mu_1\epsilon_1C\|\nabla^4 u\|_{L^2}^2
+\mu_2\epsilon_1C\|\nabla^3\mathrm{div}u\|_{L^2}^2.
\end{align*}
Firstly, choosing a fixed positive constant $\epsilon\leq \frac{\beta k \tilde{\eta}r_1}{64C^2}$, and then taking
\[\epsilon_1\leq\min\left\{\frac{\epsilon_2\beta k \tilde{\eta}}{8C}, \frac{r_2\varepsilon}{16C\beta\tilde{\eta}},\frac{r_3\varepsilon}{16C_\epsilon\beta k \tilde{\eta}}, \frac{1}{4C}\right\}\]
and
\[\epsilon_2\leq\min\left\{\frac{\epsilon_1r_1}{8C\epsilon}, \frac{r_2\varepsilon}{16C\epsilon\beta\tilde{\eta}},\frac{r_3\varepsilon}{16C_\epsilon\beta k \tilde{\eta}},\frac{r_3}{4C\beta k \tilde{\eta}}\right\},\]
and finally choosing $\delta$ sufficiently small, we obtain
\begin{align}
&\frac{1}{2}\frac{d}{dt}\big(\|\rho\|_{H^3}^2
+\| u\|_{H^3}^2+\frac{r_2}{\beta\tilde{\eta}}\|\eta\|_{H^3}^2
+\frac{r_3}{2\beta
k\tilde{\eta}}\|\tau\|_{H^3}^2\big)\label{b19}\\
&+\frac{d}{dt}\int\big(\epsilon_1\sum\limits_{\ell=0}^{2}\nabla^\ell u:\nabla^{\ell+1} \rho+\epsilon_2\sum\limits_{\ell=0}^{2}\nabla^\ell\mathrm{div}\tau:\nabla^\ell u-\frac{1}{2}\frac{h(\rho)+\beta\rho}{r_1+\beta\rho}|\nabla^3\rho|^2\big)\,\mathrm{d}x
+\frac{\mu_1}{4}\|\nabla u\|_{H^3}^2+\frac{\mu_2}{4}\|\mathrm{div}u\|_{H^3}^2\notag\\
&+\epsilon_1\frac{r_1}{4}\|\nabla\rho\|_{H^2}^2
+\epsilon_2\frac{\beta k \tilde{\eta}}{4}\|\nabla u\|_{H^2}^2
+\frac{r_2\varepsilon}{4\beta\tilde{\eta}}\|\nabla\eta\|_{H^3}^2
+\frac{A_0 r_3}{8\lambda\beta
k\tilde{\eta}}\|\tau\|_{H^3}^2
+\frac{r_3\varepsilon}{8\beta
k\tilde{\eta}}\|\nabla\tau\|_{H^3}^2
\leq0.\notag
\end{align}
Next, integrating \eqref{b19} over $(0,t)$, we get
\begin{align}
&\mathcal{J}(t)
+\int_0^t\big(\epsilon_1\frac{r_1}{4}\|\nabla\rho\|_{H^2}^2
+\epsilon_2\frac{\beta k \tilde{\eta}}{4}\|\nabla u\|_{H^2}^2
+\frac{r_2\varepsilon}{4\beta\tilde{\eta}}\|\nabla\eta\|_{H^3}^2
+\frac{A_0 r_3}{8\lambda\beta
k\tilde{\eta}}\|\tau\|_{H^3}^2
+\frac{r_3\varepsilon}{8\beta
k\tilde{\eta}}\|\nabla\tau\|_{H^3}^2\big)\mathrm{d}s\label{b20}\\
&+\int_0^t\big(\frac{\mu_1}{4}\|\nabla u\|_{H^3}^2+\frac{\mu_2}{4}\|\mathrm{div}u\|_{H^3}^2\big)\mathrm{d}s
\leq\mathcal{J}(0),\notag
\end{align}
where
\begin{align*}
\mathcal{J}(t)=&\frac{1}{2}\big(\|\rho\|_{H^3}^2
+\| u\|_{H^3}^2+\frac{r_2}{\beta\tilde{\eta}}\|\eta\|_{H^3}^2
+\frac{r_3}{2\beta
k\tilde{\eta}}\|\tau\|_{H^3}^2\big)\notag\\
&+\int\big(\epsilon_1\sum\limits_{\ell=0}^{2}\nabla^\ell u:\nabla^{\ell+1} \rho+\epsilon_2\sum\limits_{\ell=0}^{2}\nabla^\ell\mathrm{div}\tau:\nabla^\ell u-\frac{1}{2}\frac{h(\rho)+\beta\rho}{r_1+\beta\rho}|\nabla^3\rho|^2\big)\,\mathrm{d}x.
\end{align*}

Since \eqref{clc9} and the smallness of $\delta$, $\epsilon_1$ and $\epsilon_2$, it is easy to check that $\mathcal{J}(t)$ is equivalent to $$\|\rho\|_{H^3}^2+\|u\|_{H^3}^2+\|\eta\|_{H^3}^2+\| \tau\|_{H^3}^2.$$
Moreover, by the virtue of \eqref{T1} and \eqref{b20}, there exists a constant $C_1$ independent of $\mu,\nu,\theta,\delta,\epsilon,\epsilon_1,\epsilon_2,t$ and $T^*$, such that
$$
\|(\rho, u, \eta,\tau)(t)\|_{H^3}\leq C_1\theta.
$$
Letting
\begin{equation}\label{b18}
C_1\theta\leq\frac{2}{3}\delta,
\end{equation}
and then we can get \eqref{clc9-} and complete the proof of Proposition \ref{proposition4}.
\endpf

\bigskip

Due to the \emph{priori} estimate stated in Proposition \ref{proposition4} and the standard continuity arguments, the following estimate
\beq\label{n10}
\|(\rho, u, \eta,\tau)(t)\|_{H^3}\leq \delta,\;\;\; \forall\; t\in[0,T^*),
\eeq
holds. Next, we only need to show $T^*=\infty$. In fact, owing to Proposition \ref{Proposition0} and the time-uniform estimates stated in \eqref{n10}, it concludes that $T^*=\infty$. Thus, we get the global existence of solutions to the initial-value
problem $(\ref{-17})$ and $(\ref{initial-condition1})$.

Thanks to \eqref{clc9-}, we deduce that $(\rho^{\mu,\nu}, u^{\mu,\nu},\eta^{\mu,\nu},\tau^{\mu,\nu})$, $(\rho^{\mu,\nu}_t, u^{\mu,\nu}_t)$ and $(\eta^{\mu,\nu}_t,\tau^{\mu,\nu}_t)$ are uniformly bounded in $L^\infty([0,\infty);H^3)$, $L^\infty([0,\infty);H^2)$ and $L^\infty([0,\infty);H^1)$, respectively. Moreover, $(\eta^{\mu,\nu},\tau^{\mu,\nu})$ is uniformly bounded in $L^2([0,\infty);H^4)$. Hence, there exists a subsequence $(\rho^{\mu,\nu}, u^{\mu,\nu}, \eta^{\mu,\nu},\tau^{\mu,\nu})$ such that
\[(\rho^{\mu,\nu}, u^{\mu,\nu}, \eta^{\mu,\nu}, \tau^{\mu,\nu})\xrightarrow{*}(\rho,u,\eta,\tau) \;\;\text{in}\;\;L^\infty([0,\infty);H^3),\]
\[(\eta^{\mu,\nu}, \tau^{\mu,\nu})\xrightarrow{w}(\eta,\tau) \;\;\text{in}\;\;L^2([0,\infty);H^4),\]
\[(\rho^{\mu,\nu}, u^{\mu,\nu}, \eta^{\mu,\nu}, \tau^{\mu,\nu})\rightarrow(\rho,u,\eta,\tau) \;\;\text{in}\;\;\mathcal{C}_{loc}([0,\infty);H_{loc}^2).\]
The regularity of the limit is good enough to ensure that $(\rho, u, \eta,\tau)$ is a strong solution to the original problem $(\ref{17})$-$(\ref{initial-condition1})$. Therefore the proof for the first part of Theorem \ref{Theorem1}, i.e., global existence, is complete. By the standard energy method, we can prove that the solution in Theorem \ref{Theorem1} is unique, provided that $\|(\rho, u, \eta,\tau)(t)\|_{H^3}$ is sufficiently small. Therefore, we finish the proof of Theorem \ref{Theorem1}.


\section{Proof of Theorem \ref{Theorem2}}\label{sec4}
In this part, we are going to obtain the decay estimates of $(\rho,u,\eta,\tau)$ to the original problem $(\ref{17})$-$(\ref{initial-condition1})$. To do this, the strategy is to combine the energy method with the spectral analysis of the corresponding linear system. The connection between the two aspects is the Duhamel's principle.
\begin{proposition}\label{prodecay} Under the assumptions of Theorem \ref{Theorem2}, there is a constant $C$ independent of $t$ such that the solution $(\rho,u,\eta,\tau)$ of initial-value problem $(\ref{17})$-$(\ref{initial-condition1})$ satisfies the following estimates:
\begin{align*}
&\|\nabla^m\tau(t)\|_{L^2(\mathbb{R}^3)}\leq C(1+t)^{-\frac{5}{4}-\frac{m}{2}},\;\;\;m=0,1,2,\\
&\|\nabla^m(\rho, u, \eta)(t)\|_{L^2(\mathbb{R}^3)}\leq C(1+t)^{-\frac{3}{4}-\frac{m}{2}},\;\;\;m=0,1,2,3,\\
&\|\nabla^3\tau(t)\|_{L^2(\mathbb{R}^3)}\leq C(1+t)^{-\frac{9}{4}},
\end{align*}for any $t\geq0$.
\end{proposition}
\pf
The proof of Proposition \ref{prodecay} consists of Propositions \ref{lemmas1}-\ref{lemmas04}.
\endpf

\subsection{Energy estimates}
First of all, we are going to get the optimal time-decay estimate of $\|\nabla(\rho,u,\eta)(t)\|_{L^2}$.
\begin{lemma}\label{lemmaa1} Under the same assumptions of Theorem \ref{Theorem2} and \eqref{clc9}, there exist two small positive constants $\epsilon_3$ and $\epsilon_4$ which will be determined in the proof of the lemma and Lemma \ref{lemmas01}, such that
\beq\label{clc-31}
\frac{1}{2}\frac{d}{dt}\mathcal{H}_1(t)
+\epsilon_3\frac{r_1}{4}\|\nabla^2\rho\|_{H^1}^2
+\epsilon_4\frac{\beta k\tilde{\eta}}{4}\|\nabla^2u\|_{H^1}^2
+\frac{r_2\varepsilon}{4\beta\tilde{\eta}}\|\nabla^2\eta\|_{H^2}^2
+\frac{r_3A_0}{8\lambda\beta k\tilde{\eta}}\|\nabla\tau\|_{H^2}^2
\leq 0,
\eeq
where
\begin{align*}
\mathcal{H}_1(t)=&\|\nabla \rho\|_{H^2}^2+\|\nabla u\|_{H^2}^2+\frac{r_2}{\beta\tilde{\eta}}\|\nabla\eta\|_{H^2}^2
+\frac{r_3}{2\beta k\tilde{\eta}}\|\nabla \tau\|_{H^2}^2\\
&+\int\big(2\epsilon_3\sum\limits_{\ell=1}^{2}\nabla^\ell u:\nabla^{\ell+1} \rho+2\epsilon_4\sum\limits_{\ell=1}^{2}\nabla^\ell\mathrm{div}\tau:\nabla^\ell u-\frac{h(\rho)+\beta\rho}{r_1+\beta\rho}|\nabla^3\rho|^2\big)\,\mathrm{d}x.
\end{align*}
\end{lemma}
\pf
Following arguments similar to the proof of Lemma \ref{lemma-n1} for the case $\ell=1,2,3$ and $\mu=\nu=0$, one has,
\beq\label{b21}
\begin{split}
&\frac{1}{2}\frac{d}{dt}\big(\|\nabla\rho\|_{H^2}^2
+\|\nabla u\|_{H^2}^2+\frac{r_2}{\beta\tilde{\eta}}\|\nabla\eta\|_{H^2}^2
+\frac{r_3}{2\beta
k\tilde{\eta}}\|\nabla\tau\|_{H^2}^2
-\int\frac{h(\rho)+\beta\rho}{r_1+\beta\rho}|\nabla^3\rho|^2\,\mathrm{d}x\big)\\
&+\frac{r_2\varepsilon}{2\beta\tilde{\eta}}\|\nabla^2\eta\|_{H^2}^2
+\frac{A_0 r_3}{4\lambda\beta
k\tilde{\eta}}\|\nabla\tau\|_{H^2}^2
+\frac{r_3\varepsilon}{4\beta
k\tilde{\eta}}\|\nabla^2\tau\|_{H^2}^2\\
\leq& C\delta(\|\nabla^2\rho\|_{H^1}^2+\|\nabla^2 u\|_{H^1}^2).
\end{split}
\eeq
In addition, for the case $\ell=1,2$ and $\mu=\nu=0$ in \eqref{n6} and \eqref{n7}, we can get
\begin{align}
&\frac{d}{dt}\sum\limits_{\ell=1}^{2}\int\nabla^\ell\mathrm{div}\tau:\nabla^\ell u\,\mathrm{d}x+\frac{\beta k \tilde{\eta}}{2}(\|\nabla^2 u\|_{H^1}^2+\|\nabla\mathrm{div}u\|_{H^1}^2)\label{b22}\\
\leq&C(\epsilon+\delta)\|\nabla^2 \rho\|_{H^1}^2+C(\epsilon+\delta)\|\nabla^2 \eta\|_{H^1}^2+C_\epsilon\|\nabla^2\tau\|_{H^2}^2
+\frac{A_0}{2\lambda}C\|\nabla\tau\|_{H^1}^2,\notag
\end{align}
and
\beq\label{b23}
\frac{d}{dt}\sum\limits_{\ell=1}^{2}\int\nabla^\ell u:\nabla^{\ell+1} \rho\,\mathrm{d}x+\frac{r_1}{2}\|\nabla^2\rho\|_{H^1}^2
\leq C\|\nabla^2 u\|_{H^1}^{2}+C\|\nabla^2\eta\|_{H^1}^{2}+C\|\nabla^2 \tau\|_{H^1}^{2}.
\eeq
Hence, $\epsilon_3$\eqref{b23} together with \eqref{b21} and $\epsilon_4$\eqref{b22} yields
\begin{align*}
&\frac{1}{2}\frac{d}{dt}\mathcal{H}_1(t)
+\epsilon_3\frac{r_1}{2}\|\nabla^2\rho\|_{H^1}^2
+\epsilon_4\frac{\beta k\tilde{\eta}}{2}\|\nabla^2u\|_{H^1}^2
+\frac{r_2\varepsilon}{2\beta\tilde{\eta}}\|\nabla^2\eta\|_{H^2}^2
+\frac{A_0 r_3}{4\lambda\beta k\tilde{\eta}}\|\nabla\tau\|_{H^2}^2
+\frac{r_3\varepsilon}{4\beta
k\tilde{\eta}}\|\nabla^2\tau\|_{H^2}^2\\
\leq&C(\epsilon_4\epsilon+\delta)\|\nabla^2 \rho\|_{H^1}^2 +C(\epsilon_3+\delta)\|\nabla^2 u\|_{H^1}^{2}+C(\epsilon_3+\epsilon_4\epsilon+\delta)\|\nabla^2 \eta\|_{H^1}^2\\
&+C_\epsilon(\epsilon_4+\epsilon_3)\|\nabla^2\tau\|_{H^2}^2
+\epsilon_4\frac{A_0}{2\lambda}C\|\nabla\tau\|_{H^1}^2.
\end{align*}
Firstly, choosing a fixed positive constant $\epsilon\leq \frac{\beta k \tilde{\eta}r_1}{64C^2}$, and taking
\[\epsilon_3\leq\min\left\{\frac{\epsilon_4\beta k \tilde{\eta}}{8C}, \frac{r_2\varepsilon}{16C\beta\tilde{\eta}},\frac{r_3\varepsilon}{16C_\epsilon\beta k \tilde{\eta}}\right\}\]
and
\[\epsilon_4\leq\min\left\{\frac{\epsilon_3r_1}{8C\epsilon}, \frac{r_2\varepsilon}{16C\epsilon\beta\tilde{\eta}},\frac{r_3\varepsilon}{16C_\epsilon\beta k \tilde{\eta}},\frac{r_3}{4C\beta k \tilde{\eta}}\right\},\]
and finally choosing $\delta$ sufficiently small, we get \eqref{clc-31}.
\endpf

\bigskip

Moreover, with Lemmas \ref{lemmaa1} and $\ref{lemma}$, the following result holds.
\begin{lemma}\label{lemmas01}Under the same assumptions of Theorem \ref{Theorem2} and $(\ref{clc9})$, we have
\begin{align}
\mathcal{H}_1(t)\leq e^{-C_2t}\mathcal{H}_1(0)+C\int_0^te^{-C_2(t-s)}\left(\|\nabla \rho^L\|_{L^2}^2+\|\nabla u^L\|_{L^2}^2+\|\nabla \eta^L\|_{L^2}^2\right)(s)\mathrm{d}s,\label{e12}
\end{align}
for some positive constant $C_2$ independent of $\delta$.
\end{lemma}
\pf
By Lemma $\ref{lemma}$, we have
\begin{equation*}
c_0\|\nabla \rho^h\|_{L^2}\leq \|\nabla^2 \rho\|_{L^2},\;\;c_0\|\nabla u^h\|_{L^2}\leq \|\nabla^2 u\|_{L^2},\;\;c_0\|\nabla \eta^h\|_{L^2}\leq \|\nabla^2 \eta\|_{L^2}.
\end{equation*}
Thus, $(\ref{clc-31})$ leads to
\begin{align}
&\frac{1}{2}\frac{d}{dt}\mathcal{H}_1(t)
+\epsilon_3\frac{r_1c_0^2}{8}\|\nabla\rho^h\|_{L^2}^2
+\epsilon_3\frac{r_1}{8}\|\nabla^2 \rho\|_{H^1}^2+\epsilon_4\frac{\beta k\tilde{\eta}c_0^2}{8}\|\nabla u^h\|_{L^2}^2+\epsilon_4\frac{\beta k\tilde{\eta}}{8}\|\nabla^2 u\|_{H^1}^2
\label{clc-s4}\\
&+\frac{r_2\varepsilon c_0^2}{8\beta\tilde{\eta}}\|\nabla\eta^h\|_{L^2}^2
+\frac{r_2\varepsilon}{8\beta\tilde{\eta}}\|\nabla^2\eta\|_{H^2}^2
+\frac{r_3A_0}{8\lambda\beta k\tilde{\eta}}\|\nabla\tau\|_{H^2}^2
\leq 0.\notag
\end{align}
By adding $\epsilon_3\frac{r_1c_0^2}{8}\|\nabla \rho^L\|_{L^2}^2+ \epsilon_4\frac{\beta k\tilde{\eta}c_0^2}{8}\|\nabla u^L\|_{L^2}^2+\frac{r_2\varepsilon c_0^2}{8\beta\tilde{\eta}}\|\nabla\eta^L\|_{L^2}^2$ to both sides of inequality $(\ref{clc-s4})$, we have
\begin{align*}
&\frac{1}{2}\frac{d}{dt}\mathcal{H}_1(t)
+\epsilon_3\frac{r_1c_0^2}{8}\|\nabla\rho\|_{H^2}^2
+\epsilon_4\frac{\beta k\tilde{\eta}c_0^2}{8}\|\nabla u\|_{H^2}^2
+\frac{r_2\varepsilon c_0^2}{8\beta\tilde{\eta}}\|\nabla\eta\|_{H^2}^2
+\frac{r_3A_0}{8\lambda\beta k\tilde{\eta}}\|\nabla\tau\|_{H^2}^2\\
\leq&\epsilon_3\frac{r_1c_0^2}{8}\|\nabla \rho^L\|_{L^2}^2
+\epsilon_4\frac{\beta k\tilde{\eta}c_0^2}{8}\|\nabla u^L\|_{L^2}^2+\frac{r_2\varepsilon c_0^2}{8\beta\tilde{\eta}}\|\nabla\eta^L\|_{L^2}^2,
\end{align*} where we let $c_0\in (0,1]$.

Note that, by virtue of  \eqref{clc9} and the smallness of $\delta$, $\epsilon_3$ and $\epsilon_4$,  it is easy to check that $\mathcal{H}_1(t)$ is equivalent to $$\|\nabla \rho\|_{H^2}^2+\|\nabla u\|_{H^2}^2+\|\nabla \eta\|_{H^2}^2+\|\nabla \tau\|_{H^2}^2.$$
Then there exists a positive constant $C_2>0$ such that
\begin{align*}
\frac{d}{dt}\mathcal{H}_1(t)+C_2\,\mathcal{H}_1(t)\leq C \|\nabla \rho^L\|_{L^2}^2+C\|\nabla u^L\|_{L^2}^2+C\|\nabla \eta^L\|_{L^2}^2.
\end{align*}
By using Gronwall's inequality, we get \eqref{e12}.
\endpf

\bigskip

In the same way, we show the following estimates of $\|\nabla^2(\rho,u,\eta,\tau)(t)\|_{H^1}$ which are the basis for getting the optimal decay estimate of $\|\nabla^2(\rho,u,\eta)(t)\|_{L^2}$.
\begin{lemma}\label{proposition2} Under the same assumptions of Theorem \ref{Theorem2} and \eqref{clc9}, there exist two small positive constants $\epsilon_5$ and $\epsilon_6$ which will be determined in the proof of the lemma and Lemma \ref{lemmas02}, such that
\begin{align}
&\frac{1}{2}\frac{d}{dt}\mathcal{H}_2(t)
+\epsilon_5\frac{r_1}{4}\|\nabla^3\rho\|_{L^2}^2
+\epsilon_6\frac{\beta k\tilde{\eta}}{4}\|\nabla^3u\|_{L^2}^2
+\frac{r_2\varepsilon}{4\beta\tilde{\eta}}\|\nabla^3\eta\|_{H^1}^2
+\frac{r_3A_0}{8\lambda\beta k\tilde{\eta}}\|\nabla^2\tau\|_{H^1}^2\label{b24}\\
\leq& C\delta\|\nabla^2 u\|_{L^2}^{2}+C\delta\|\nabla^2 \rho\|_{L^2}^{2},\notag
\end{align}
where
\begin{align*}
\mathcal{H}_2(t)=&\|\nabla^2 \rho\|_{H^1}^2+\|\nabla^2 u\|_{H^1}^2+\frac{r_2}{\beta\tilde{\eta}}\|\nabla^2\eta\|_{H^1}^2
+\frac{r_3}{2\beta k\tilde{\eta}}\|\nabla^2 \tau\|_{H^1}^2\\
&+\int\big(2\epsilon_5\nabla^2 u:\nabla^3 \rho+2\epsilon_6\nabla^2\mathrm{div}\tau:\nabla^2 u-\frac{h(\rho)+\beta\rho}{r_1+\beta\rho}|\nabla^3\rho|^2\big)\,\mathrm{d}x.
\end{align*}
\end{lemma}
\pf
Following arguments similar to the proof of Lemma \ref{lemma-n1} for the case $\ell=2,3$ and $\mu=\nu=0$, one has,
\beq\label{b25}
\begin{split}
&\frac{1}{2}\frac{d}{dt}\big(\|\nabla^2\rho\|_{H^1}^2
+\|\nabla^2 u\|_{H^1}^2+\frac{r_2}{\beta\tilde{\eta}}\|\nabla^2\eta\|_{H^1}^2
+\frac{r_3}{2\beta
k\tilde{\eta}}\|\nabla^2\tau\|_{H^1}^2
-\int\frac{h(\rho)+\beta\rho}{r_1+\beta\rho}|\nabla^3\rho|^2\,\mathrm{d}x\big)\\
&+\frac{r_2\varepsilon}{2\beta\tilde{\eta}}\|\nabla^3\eta\|_{H^1}^2
+\frac{A_0 r_3}{4\lambda\beta
k\tilde{\eta}}\|\nabla^2\tau\|_{H^1}^2
+\frac{r_3\varepsilon}{4\beta
k\tilde{\eta}}\|\nabla^3\tau\|_{H^1}^2\\
\leq& C\delta\big(\|\nabla^2u\|_{H^1}^2+\|\nabla^2\rho\|_{H^1}^2\big).
\end{split}
\eeq
In addition, for the case $\ell=2$ and $\mu=\nu=0$ in \eqref{n6} and \eqref{n7}, we can deduce
\begin{align}
&\frac{d}{dt}\int\nabla^2\mathrm{div}\tau:\nabla^2 u\,\mathrm{d}x+\frac{\beta k \tilde{\eta}}{2}(\|\nabla^3 u\|_{L^2}^2+\|\nabla^2\mathrm{div}u\|_{L^2}^2)\label{b26}\\
\leq&C(\epsilon+\delta)\|\nabla^3 \rho\|_{L^2}^2+C\delta\|\nabla^2 \rho\|_{L^2}^2+C\delta\|\nabla^2 u\|_{L^2}^2+C(\epsilon+\delta)\|\nabla^3 \eta\|_{L^2}^2+C_\epsilon\|\nabla^3\tau\|_{H^1}^2
+\frac{A_0}{2\lambda}C\|\nabla^2\tau\|_{L^2}^2,\notag
\end{align}
and
\beq\label{b27}
\frac{d}{dt}\int\nabla^2 u:\nabla^3 \rho\,\mathrm{d}x+\frac{r_1}{2}\|\nabla^3\rho\|_{L^2}^2
\leq C\delta\|\nabla^2 \rho\|_{L^2}^{2}+C\|\nabla^3 u\|_{L^2}^{2}+C\|\nabla^3\eta\|_{L^2}^{2}+C\|\nabla^3 \tau\|_{L^2}^{2}.
\eeq

Together with $\epsilon_5$\eqref{b27}, \eqref{b25} and $\epsilon_6$\eqref{b26} yields
\begin{align*}
&\frac{1}{2}\frac{d}{dt}\mathcal{H}_2(t)
+\epsilon_5\frac{r_1}{2}\|\nabla^3\rho\|_{L^2}^2
+\epsilon_6\frac{\beta k\tilde{\eta}}{2}\|\nabla^3u\|_{L^2}^2
+\frac{r_2\varepsilon}{2\beta\tilde{\eta}}\|\nabla^3\eta\|_{H^1}^2
+\frac{A_0 r_3}{4\lambda\beta k\tilde{\eta}}\|\nabla^2\tau\|_{H^1}^2
+\frac{r_3\varepsilon}{4\beta
k\tilde{\eta}}\|\nabla^3\tau\|_{H^1}^2\\
\leq&C\epsilon_6\epsilon\|\nabla^3 \rho\|_{L^2}^2
+C\delta\|\nabla^2 \rho\|_{H^1}^2 +C\delta\|\nabla^2 u\|_{H^1}^{2}+C\epsilon_5\|\nabla^3 u\|_{L^2}^{2}
+C(\epsilon_5+\epsilon_6\epsilon+\delta)\|\nabla^3 \eta\|_{L^2}^2\\
&+C_\epsilon\epsilon_6\|\nabla^3\tau\|_{H^1}^2
+C\epsilon_5\|\nabla^3\tau\|_{L^2}^2
+\epsilon_6\frac{A_0}{2\lambda}C\|\nabla^2\tau\|_{L^2}^2.
\end{align*}
Firstly, choosing a fixed positive constant $\epsilon\leq \frac{\beta k \tilde{\eta}r_1}{64C^2}$, and taking
\[\epsilon_5\leq\min\left\{\frac{\epsilon_6\beta k \tilde{\eta}}{8C}, \frac{r_2\varepsilon}{16C\beta\tilde{\eta}},\frac{r_3\varepsilon}{16C\beta k \tilde{\eta}}\right\}\]
and
\[\epsilon_6\leq\min\left\{\frac{\epsilon_5r_1}{8C\epsilon}, \frac{r_2\varepsilon}{16C\epsilon\beta\tilde{\eta}},\frac{r_3\varepsilon}{16C_\epsilon\beta k \tilde{\eta}},\frac{r_3}{4C\beta k \tilde{\eta}}\right\},\]
and finally choosing $\delta$ sufficiently small, we get \eqref{b24}.
\endpf

\bigskip

Based on Lemmas \ref{proposition2} and $\ref{lemma}$, the following result holds.
\begin{lemma}\label{lemmas02}Under the same assumptions of Theorem \ref{Theorem2} and $(\ref{clc9})$, we have
\begin{align}
\mathcal{H}_2(t)\leq e^{-\tilde{C_2}t}\mathcal{H}_2(0)+C\int_0^te^{-\tilde{C_2}(t-s)}\left(\|\nabla^2 \rho^L\|_{L^2}^2+\|\nabla^2 u^L\|_{L^2}^2+\|\nabla^2 \eta^L\|_{L^2}^2\right)(s)\mathrm{d}s,\label{b28}
\end{align}
for some positive constant $\tilde{C_2}$ independent of $\delta$.
\end{lemma}
\pf
By Lemma $\ref{lemma}$, we have
\begin{equation*}
c_0\|\nabla^2 \rho^h\|_{L^2}\leq \|\nabla^3 \rho\|_{L^2},\;\;c_0\|\nabla^2 u^h\|_{L^2}\leq \|\nabla^3 u\|_{L^2},\;\;c_0\|\nabla^2 \eta^h\|_{L^2}\leq \|\nabla^3 \eta\|_{L^2}.
\end{equation*}
Thus, \eqref{b24} leads to
\begin{align}
&\frac{1}{2}\frac{d}{dt}\mathcal{H}_2(t)
+\epsilon_5\frac{r_1c_0^2}{8}\|\nabla^2\rho^h\|_{L^2}^2
+\epsilon_5\frac{r_1}{8}\|\nabla^3 \rho\|_{L^2}^2+\epsilon_6\frac{\beta k\tilde{\eta}c_0^2}{8}\|\nabla^2 u^h\|_{L^2}^2+\epsilon_6\frac{\beta k\tilde{\eta}}{8}\|\nabla^3 u\|_{L^2}^2
\label{b29}\\
&+\frac{r_2\varepsilon c_0^2}{8\beta\tilde{\eta}}\|\nabla^2\eta^h\|_{L^2}^2
+\frac{r_2\varepsilon}{8\beta\tilde{\eta}}\|\nabla^3\eta\|_{H^1}^2
+\frac{r_3A_0}{8\lambda\beta k\tilde{\eta}}\|\nabla^2\tau\|_{H^1}^2
\leq C\delta\|\nabla^2 u\|_{L^2}^{2}+C\delta\|\nabla^2 \rho\|_{L^2}^{2}.\notag
\end{align}
By adding $\epsilon_5\frac{r_1c_0^2}{8}\|\nabla^2 \rho^L\|_{L^2}^2+ \epsilon_6\frac{\beta k\tilde{\eta}c_0^2}{8}\|\nabla^2 u^L\|_{L^2}^2+\frac{r_2\varepsilon c_0^2}{8\beta\tilde{\eta}}\|\nabla^2\eta^L\|_{L^2}^2$ to both sides of inequality \eqref{b29}, and choosing $\delta$ sufficiently small, we have
\begin{align*}
&\frac{1}{2}\frac{d}{dt}\mathcal{H}_2(t)
+\epsilon_5\frac{r_1c_0^2}{16}\|\nabla^2\rho\|_{H^1}^2
+\epsilon_6\frac{\beta k\tilde{\eta}c_0^2}{16}\|\nabla^2 u\|_{H^1}^2
+\frac{r_2\varepsilon c_0^2}{8\beta\tilde{\eta}}\|\nabla^2\eta\|_{H^1}^2
+\frac{r_3A_0}{8\lambda\beta k\tilde{\eta}}\|\nabla^2\tau\|_{H^1}^2\\
\leq&\epsilon_5\frac{r_1c_0^2}{8}\|\nabla^2 \rho^L\|_{L^2}^2
+\epsilon_6\frac{\beta k\tilde{\eta}c_0^2}{8}\|\nabla^2 u^L\|_{L^2}^2+\frac{r_2\varepsilon c_0^2}{8\beta\tilde{\eta}}\|\nabla^2\eta^L\|_{L^2}^2,
\end{align*} where we let $c_0\in (0,1]$.

Moreover, by virtue of  \eqref{clc9} and the smallness of $\delta$, $\epsilon_5$ and $\epsilon_6$,  it is easy to check that $\mathcal{H}_2(t)$ is equivalent to $$\|\nabla^2 \rho\|_{H^1}^2+\|\nabla^2 u\|_{H^1}^2+\|\nabla^2 \eta\|_{H^1}^2+\|\nabla^2 \tau\|_{H^1}^2.$$
Then there exists a positive constant $\tilde{C_2}>0$ such that
\begin{align*}
\frac{d}{dt}\mathcal{H}_2(t)+\tilde{C_2}\mathcal{H}_2(t)\leq C \|\nabla^2 \rho^L\|_{L^2}^2+C\|\nabla^2 u^L\|_{L^2}^2+C\|\nabla^2 \eta^L\|_{L^2}^2.
\end{align*}
By using Gronwall's inequality again, we get \eqref{b28} directly.
\endpf

\subsection{Decay estimates of the low-frequency parts}
Next, with the help of Lemmas \ref{lemmas01} and \ref{lemmas02}, we will study the decay rates of solution $(\rho,u,\eta,\tau)$. As it can be seen from \eqref{e12} and \eqref{b28}, we only need to analyze the low-frequency part ($|\xi|\leq c_0$) of $(\rho,u,\eta)$.

Letting $\mathbb{A}$ be the following matrix of differential operators of the form
\begin{equation*}       
\mathbb{A}=
\left(                 
  \begin{array}{cccc}   
    0 & r_1\mathrm{div} & 0 & 0\\  
    r_1\nabla & 0 & r_2\nabla & -r_3\\  
    0 & \beta\tilde{\eta}\mathrm{div} & -\varepsilon\Delta & 0\\
    0 & -\beta k\tilde{\eta}(\Delta+\nabla\mathrm{div}) & 0 & \frac{A_0}{2\lambda}-\varepsilon\Delta\\
  \end{array}
\right),                 
\end{equation*}
and setting
$$\bar{\mathbb{U}}(t):=(\bar{\rho}(t),\bar{u}(t),\bar{\eta}(t),\mathrm{div}\bar{\tau}(t))^T\;
\;\text{and}\;\;\mathbb{U}(0):=(\rho_0,u_0,\eta_0,\mathrm{div}\tau_0)^T,$$
we obtain from the linearized problem of \eqref{A1} as below:
\begin{equation}\label{clc-41}
\left\{
\begin{split}
&\partial_t\bar{\mathbb{U}}+\mathbb{A}\bar{\mathbb{U}}= 0,\;\text{for}\;\;t>0,\\
&\bar{\mathbb{U}}\big|_{t=0}=\mathbb{U}(0).
\end{split}
\right.
\end{equation}
Applying the Fourier transform to \eqref{clc-41} with respect to the $x$-variable and solving the ordinary differential equation with respect to $t$, we have
$$\bar{\mathbb{U}}(t)=\mathcal{A}(t)\mathbb{U}(0),$$
where $\mathcal{A}(t)=e^{-t\mathbb{A}}(t\geq0)$  is the semigroup generated by the linear operator $\mathbb{A}$ and
$\mathcal{A}(t)f:=\mathcal{F}^{-1}(e^{-t\mathbb{A}_{\xi}}\hat{f}(\xi))$ with
\begin{equation*}       
\mathbb{A}_{\xi}=
\left(                 
  \begin{array}{cccc}   
    0 & \sqrt{-1}r_1\xi^T & 0 & 0\\  
    \sqrt{-1}r_1\xi & 0 & \sqrt{-1}r_2\xi & -r_3\\  
    0 & \sqrt{-1}\beta\tilde{\eta}\xi^T & \varepsilon|\xi|^2 & 0\\
    0 & \beta k\tilde{\eta}(|\xi|^2\delta_{ij}+\xi_i\xi_j) & 0 & (\frac{A_0}{2\lambda}+\varepsilon|\xi|^2)\delta_{ij}\\
  \end{array}
\right).                 
\end{equation*}
Then, from Proposition \ref{proA2}, we have the following result.
\begin{lemma}\label{lemmal-m} For any integer $m\geq 0$, the following time-decay estimates for the low-frequency part, i.e.,
\beq \label{a1}
\|\nabla^m(\mathcal{A}(t)\mathbb{U}^L(0))\|_{L^2}
\leq C (1+t)^{-\frac{3}{4}-\frac{m}{2}}\|\mathbb{U}(0)\|_{L^1},
\eeq
and
\beq \label{a2}
\|\nabla^m(\mathcal{A}(t)\mathbb{U}^L(0))\|_{L^2}
\leq C (1+t)^{-\frac{m}{2}}\|\mathbb{U}(0)\|_{L^2}
\eeq hold for any $t\in(0,\infty)$.
\end{lemma}

\begin{remark} \eqref{a1} and \eqref{a2} are used to obtain the optimal time-decay estimates of $\|(\rho,u,\eta,\tau)(t)\|_{H^2}$ and $\|\nabla^3(\rho,u,\eta,\tau)(t)\|_{L^2}$, respectively.
\end{remark}

In what follows, based on the estimates in Lemma \ref{lemmal-m}, we establish time-decay estimates for the low-frequency part of solutions to the nonlinear problem \eqref{17}-\eqref{initial-condition1}. Denoting
$$\mathbb{U}(t):=(\rho(t),u(t),\eta(t),\mathrm{div}\tau(t))^T,$$
then from \eqref{18}, we have
\begin{equation}\label{clc-42}
\left\{
\begin{split}
&\partial_t\mathbb{U}+\mathbb{A}\mathbb{U}=S(\mathbb{U}) ,\;\text{for}\;\;t>0,\\
&\mathbb{U}\big|_{t=0}=\mathbb{U}(0),
\end{split}
\right.
\end{equation}
where
$$S(\mathbb{U})=(S_1,S_2,S_3,\mathrm{div}S_4)^T.$$
Using the Duhamel's principle, the solution of \eqref{clc-42} can be stated as follows:
\beq\label{clc-43}
\mathbb{U}(t)=\mathcal{A}(t)\mathbb{U}(0)
+\int_0^t\mathcal{A}(t-s)S(\mathbb{U})(s)\mathrm{d}s.
\eeq

\begin{lemma}\label{prol-m} Under the assumptions of Theorem \ref{Theorem2} and \eqref{clc9}, for any integer $m\geq1$, there exists a positive constant $C$ such that
\begin{align}
\|\nabla^m\mathbb{U}^L(t)\|_{L^2}
\leq & C(1+t)^{-\frac{3}{4}-\frac{m}{2}}\|\mathbb{U}(0)\|_{L^1}\label{clc-44}\\
&+C\delta\int_0^t(1+t-s)^{-\frac{3}{4}-\frac{m}{2}}\big(\|\nabla(\rho,\eta)(s)\|_{L^2}
+\|\nabla(u,\tau)(s)\|_{H^1}\big)\mathrm{d}s\notag
\end{align}
and
\begin{align}
\|\nabla^m\mathbb{U}^L(t)\|_{L^2}
\leq & C(1+t)^{-\frac{3}{4}-\frac{m}{2}}\|\mathbb{U}(0)\|_{L^1}\label{b33}\\
&+C\int_0^\frac{t}{2}(1+t-s)^{-\frac{3}{4}-\frac{m}{2}}\big(\|(\rho,u,\eta)\|_{L^2}\|\nabla(\rho,u,\eta,\tau)\|_{L^2}
+\|\nabla(\eta,\tau)\|_{L^2}\|\nabla u\|_{L^2}\big)(s)\mathrm{d}s\notag\\
&+C\int_0^\frac{t}{2}(1+t-s)^{-\frac{3}{4}-\frac{m}{2}}\big(\|(\eta,\tau)\|_{L^2}\|\nabla^2u\|_{L^2}
+\|\nabla^2\tau\|_{L^2}\|u\|_{L^2}\big)(s)\mathrm{d}s\notag\\
&+C\int_\frac{t}{2}^t(1+t-s)^{-\frac{m}{2}} \big(\|(\rho,u,\eta,\tau)\|_{H^2}\|\nabla^2(\rho,u,\eta,\tau)\|_{L^2}\big)(s)\mathrm{d}s.\notag
\end{align}
\end{lemma}
\begin{remark} With the aid of \eqref{clc-44} and \eqref{b33}, the optimal time-decay estimates of $\|(\rho,u,\eta,\tau)(t)\|_{H^2}$ and $\|\nabla^3(\rho,u,\eta,\tau)(t)\|_{L^2}$ will be obtained, respectively.
\end{remark}
\pf
From \eqref{clc-43}, using Lemma \ref{lemmal-m}, we have
\begin{align}
&\|\nabla^m\mathbb{U}^L(t)\|_{L^2}\notag\\
\leq& C\big\|\nabla^m \big(\mathcal{A}(t)\mathbb{U}^L(0)\big)\big\|_{L^2}
+C\big\|\nabla^m\int_0^t\mathcal{A}(t-s)S^L(\mathbb{U})(s)\mathrm{d}s\big\|_{L^2}\label{b31}\\
\leq & C(1+t)^{-\frac{3}{4}-\frac{m}{2}}\|\mathbb{U}(0)\|_{L^1}
+C\int_0^t(1+t-s)^{-\frac{3}{4}-\frac{m}{2}}\|S(\mathbb{U})(s)\|_{L^1}\mathrm{d}s\notag\\
\leq & C(1+t)^{-\frac{3}{4}-\frac{m}{2}}\|\mathbb{U}(0)\|_{L^1}
+C\delta\int_0^t(1+t-s)^{-\frac{3}{4}-\frac{m}{2}}\big(\|\nabla(\rho,\eta)(s)\|_{L^2}+\|\nabla(u,\tau)(s)\|_{H^1}\big)\mathrm{d}s,\notag
\end{align}
where we have used the fact that
\begin{align}
\|S(\mathbb{U})\|_{L^1}\leq&C\|(S_1,S_2,S_3,\mathrm{div}S_4)(\mathbb{U})\|_{L^1}\label{b32}\\
\leq&C\big(\|(\rho,u,\eta)\|_{L^2}\|\nabla(\rho,u,\eta,\tau)\|_{L^2}
+\|\nabla(\eta,\tau)\|_{L^2}\|\nabla u\|_{L^2}
+\|(\eta,\tau)\|_{L^2}\|\nabla^2u\|_{L^2}
+\|\nabla^2\tau\|_{L^2}\|u\|_{L^2}\big).\notag
\end{align}
Hence, \eqref{clc-44} is obtained.

Moreover, for \eqref{b33}, using \eqref{clc-43} and Lemma \ref{lemmal-m} again, we have
\begin{align}
&\|\nabla^m\mathbb{U}^L(t)\|_{L^2}\label{b36}\\
\leq& C\left\|\nabla^m (\mathcal{A}(t)\mathbb{U}^L(0))\right\|_{L^2}
+C\big\|\nabla^m\int_0^\frac{t}{2}\mathcal{A}(t-s)S^L(\mathbb{U})(s)\mathrm{d}s\big\|_{L^2}
+C\big\|\nabla^m\int_{\frac{t}{2}}^t\mathcal{A}(t-s)S^L(\mathbb{U})(s)\mathrm{d}s\big\|_{L^2}\notag\\
\leq & C(1+t)^{-\frac{3}{4}-\frac{m}{2}}\|\mathbb{U}(0)\|_{L^1}
+C\int_0^\frac{t}{2}(1+t-s)^{-\frac{3}{4}-\frac{m}{2}}\|S(\mathbb{U})(s)\|_{L^1}\mathrm{d}s
+C\int_\frac{t}{2}^t(1+t-s)^{-\frac{m}{2}}\|S(\mathbb{U})(s)\|_{L^2}\mathrm{d}s,\notag
\end{align}
where
\begin{align}
\|S(\mathbb{U})\|_{L^2}\leq&C\|(S_1,S_2,S_3,\mathrm{div}S_4)(\mathbb{U})\|_{L^2}\label{b37}\\
\leq&C\big(\|(\rho,u,\eta)\|_{L^3}\|\nabla(\rho,u,\eta,\tau)\|_{L^6}
+\|\nabla(\eta,\tau)\|_{L^6}\|\nabla u\|_{L^3}
+\|(\eta,\tau)\|_{L^\infty}\|\nabla^2u\|_{L^2}
+\|\nabla^2\tau\|_{L^2}\|u\|_{L^\infty}\big)\notag\\
\leq& C\|(\rho,u,\eta,\tau)\|_{H^2}\|\nabla^2(\rho,u,\eta,\tau)\|_{L^2}.\notag
\end{align}
Thus, together with \eqref{b32}, \eqref{b36} and \eqref{b37},   \eqref{b33} can be obtained.
\endpf

\subsection{Optimal decay rates of $\|(\rho,u,\eta)(t)\|_{H^2}$ and $\|\tau(t)\|_{H^1}$}
In this subsection, we will obtain the time-decay estimates of $(\rho, u,\eta,\tau)$ with the aid of Lemmas \ref{lemmas01}, \ref{lemmas02} and \ref{prol-m}. Firstly, we consider the decay estimate of $\|(\nabla\rho, \nabla u,\nabla\eta,\nabla\tau)(t)\|_{H^2}$.
\begin{proposition}\label{lemmas1}Under the same assumptions of Theorem \ref{Theorem2} and $(\ref{clc9})$, we have
\beq\label{clc-ss1}
\|\nabla \rho(t)\|_{H^2}+\|\nabla u(t)\|_{H^2}+\|\nabla \eta(t)\|_{H^2}+\|\nabla \tau(t)\|_{H^2}\leq C(1+t)^{-\frac{5}{4}}.
\eeq
\end{proposition}
\pf
Owing to $(\ref{e12})$ and \eqref{clc-44}, we can obtain
\begin{align}
\mathcal{H}_1(t)\leq & e^{-C_2t}\mathcal{H}_1(0)+C\int_0^te^{-C_2(t-s)}\left(\|\nabla \rho^L\|_{L^2}^2+\|\nabla u^L\|_{L^2}^2+\|\nabla \eta^L\|_{L^2}^2\right)(s)\mathrm{d}s\notag\\
\leq&e^{-C_2t}\mathcal{H}_1(0)
+C\int_0^te^{-C_2(t-s)} (1+s)^{-\frac{5}{2}}\mathrm{d}s
+\delta\int_0^te^{-C_2(t-s)} \big(\int_0^s(1+s-s')^{-\frac{5}{2}}\mathcal{H}_1(s')\mathrm{d}s'\big)\mathrm{d}s\notag\\
\leq&C(1+t)^{-\frac{5}{2}}
+C\delta\mathcal{I}(t)\int_0^te^{-C_2(t-s)}\big(\int_0^s(1+s-s')^{-\frac{5}{2}}(1+s')^{-\frac{5}{2}}\mathrm{d}s'\big)\mathrm{d}s\notag\\
\leq&C(1+t)^{-\frac{5}{2}}+C\delta\mathcal{I}(t)\int_0^te^{-C_2(t-s)}(1+s)^{-\frac{5}{2}}\mathrm{d}s\notag\\
\leq& C(1+t)^{-\frac{5}{2}}+C\delta\mathcal{I}(t)(1+t)^{-\frac{5}{2}},\notag
\end{align}
where $\mathcal{I}(t)=\sup\limits_{0\leq s\leq t}(1+s)^{\frac{5}{2}}\mathcal{H}_1(s)$.

Further, by virtue of the definition of $\mathcal{I}(t)$ and the smallness of $\delta$, we can obtain
\begin{align*}
\mathcal{I}(t)\leq C,
\end{align*}
which is $(\ref{clc-ss1})$.  We complete the proof of the proposition.\endpf

\bigskip

Then, based on Proposition \ref{lemmas1}, we can obtain the next proposition.

\begin{proposition}\label{lemmas2}Under the same assumptions of Theorem \ref{Theorem2} and $(\ref{clc9})$, we have
\begin{align}
&\|\tau(t)\|_{L^2}\leq C(1+t)^{-\frac{5}{4}},\label{clc-ss01}\\
&\|( \rho, u, \eta)(t)\|_{L^2}\leq C(1+t)^{-\frac{3}{4}}.\label{clc-ss02}
\end{align}
\end{proposition}
\pf
Firstly, for $\|\tau(t)\|_{L^2}$, we multiply $(\ref{17})_4$ by $2\tau$ and then integrate the resulting equation over $\mathbb{R}^3$ to obtain
\begin{align}
\frac{d}{dt}\int|\tau|^2\,\mathrm{d}x+\frac{A_0}{\lambda}\int|\tau|^2\,\mathrm{d}x
+2\varepsilon\int|\nabla\tau|^2\,\mathrm{d}x
=&2\int S_4:\tau\,\mathrm{d}x+2\beta k\tilde{\eta}\int(\nabla u+\nabla^Tu):\tau\,\mathrm{d}x\label{clc-46}\\
\leq&(\frac{A_0}{4\lambda}+\delta)\|\tau\|_{L^2}^2+C\|\nabla u\|_{L^2}^2,\notag
\end{align}
where we have used Cauchy inequality and the fact that
$$\|S_4\|_{L^2}\leq \|u\|_{L^6}\|\nabla\tau\|_{L^3}+\|\nabla u\|_{L^2}(\|\tau\|_{L^\infty}+\|\eta\|_{L^\infty}).$$
Choosing $\delta$ sufficiently small in \eqref{clc-46}, it holds that
$$
\frac{d}{dt}\int|\tau|^2\,\mathrm{d}x+\frac{A_0}{2\lambda}\int|\tau|^2\,\mathrm{d}x
\leq C\|\nabla u\|_{L^2}^2.
$$
Further, using Gronwall's inequality, the above inequality gives
\begin{align*}
\|\tau(t)\|_{L^2}^2&\leq C e^{-\frac{A_0}{2\lambda}t}\|\tau(0)\|_{L^2}^2+C\int_0^te^{-\frac{A_0}{2\lambda}(t-s)}\|\nabla u(s)\|_{L^2}^2\mathrm{d}s\\
&\leq C e^{-\frac{A_0}{2\lambda}t}\|\tau(0)\|_{L^2}^2+C\int_0^te^{-\frac{A_0}{2\lambda}(t-s)}(1+s)^{-\frac{5}{2}}\mathrm{d}s\\
&\leq C(1+t)^{-\frac{5}{2}},
\end{align*}
where we have used $(\ref{clc-ss1})$. Then we obtain \eqref{clc-ss01}.

Next, for $\|(\rho, u, \eta)(t)\|_{L^2}$, thanks to \eqref{clc-44}, let $m=0$, we find that
\begin{align}
\|(\rho,u,\eta,\mathrm{div}\tau)^L(t)\|_{L^2}
\leq & C(1+t)^{-\frac{3}{4}}\|(\rho,u,\eta,\mathrm{div}\tau)(0)\|_{L^1}\label{clc-47}\\
&+C\delta\int_0^t(1+t-s)^{-\frac{3}{4}}\big(\|\nabla(\rho,\eta)(s)\|_{L^2}+\|\nabla(u,\tau)(s)\|_{H^1}\big)\mathrm{d}s\notag\\
\leq&C(1+t)^{-\frac{3}{4}}
+C\delta\int_0^t(1+t-s)^{-\frac{3}{4}}(1+s)^{-\frac{5}{4}}\mathrm{d}s\notag\\
\leq&C(1+t)^{-\frac{3}{4}},\notag
\end{align}
where we have used \eqref{clc-ss1}.

In addition, by using Lemma \ref{lemma} and \eqref{clc-ss1} again, we have
\beq\label{clc-48}
\|(\rho,u,\eta,\mathrm{div}\tau)^h(t)\|_{L^2}
\leq  \frac{1}{c_0}\|\nabla(\rho,u,\eta,\mathrm{div}\tau)(t)\|_{L^2}
\leq C(1+t)^{-\frac{5}{4}}.
\eeq
Combining with \eqref{clc-47} and \eqref{clc-48}, and owing to \eqref{5.28}, we can get \eqref{clc-ss02} directly.
\endpf

\bigskip

Further, the optimal decay estimates of second order for $(\rho,u,\eta)$ and first order for $\tau$ in the sense of $L^2$ norm are obtained as below.
\begin{proposition}\label{lemmas3}Under the same assumptions of Theorem \ref{Theorem2} and $(\ref{clc9})$, we have
\begin{align}
&\|\nabla^2 \rho(t)\|_{H^1}+\|\nabla^2 u(t)\|_{H^1}+\|\nabla^2 \eta(t)\|_{H^1}+\|\nabla^2 \tau(t)\|_{H^1}\leq C(1+t)^{-\frac{7}{4}},\label{b30}\\
&\|\nabla\tau(t)\|_{L^2}\leq C(1+t)^{-\frac{7}{4}}.\label{b34}
\end{align}
\end{proposition}
\pf
Owing to \eqref{b28}, \eqref{b31}, \eqref{b32} and Propositions \ref{lemmas1} and \ref{lemmas2}, we can obtain
\begin{align}
\mathcal{H}_2(t)\leq & e^{-\tilde{C_2}t}\mathcal{H}_2(0)+C\int_0^te^{-\tilde{C_2}(t-s)}\left(\|\nabla^2 \rho^L\|_{L^2}^2+\|\nabla^2 u^L\|_{L^2}^2+\|\nabla^2 \eta^L\|_{L^2}^2\right)(s)\mathrm{d}s\notag\\
\leq&e^{-\tilde{C_2}t}\mathcal{H}_2(0)
+C\int_0^te^{-\tilde{C_2}(t-s)} (1+s)^{-\frac{7}{2}}\mathrm{d}s
+C\int_0^te^{-\tilde{C_2}(t-s)} \big(\int_0^s(1+s-s')^{-\frac{7}{2}}(1+s)^{-\frac{8}{2}}(s')\mathrm{d}s'\big)\mathrm{d}s\notag\\
\leq&C(1+t)^{-\frac{7}{2}}+C\int_0^te^{-\tilde{C_2}(t-s)}(1+s)^{-\frac{7}{2}}\mathrm{d}s\notag\\
\leq& C(1+t)^{-\frac{7}{2}},\notag
\end{align}
which is $(\ref{b30})$.

Finally, for \eqref{b34}, multiplying $\nabla(\ref{17})_4$ by $2\nabla\tau$ and then integrating the result equation over $\mathbb{R}^3$, similar to \eqref{clc-ss01}, we get
$$
\frac{d}{dt}\int|\nabla\tau|^2\,\mathrm{d}x+\frac{A_0}{2\lambda}\int|\nabla\tau|^2\,\mathrm{d}x
\leq C\|\nabla^2 u\|_{L^2}^2.
$$
Further, using Gronwall's inequality, the above inequality gives
\begin{align*}
\|\nabla\tau(t)\|_{L^2}^2&\leq C e^{-\frac{A_0}{2\lambda}t}\|\nabla\tau(0)\|_{L^2}^2+C\int_0^te^{-\frac{A_0}{2\lambda}(t-s)}\|\nabla^2 u(s)\|_{L^2}^2\mathrm{d}s\\
&\leq C e^{-\frac{A_0}{2\lambda}t}\|\nabla\tau(0)\|_{L^2}^2+C\int_0^te^{-\frac{A_0}{2\lambda}(t-s)}(1+s)^{-\frac{7}{2}}\mathrm{d}s\\
&\leq C(1+t)^{-\frac{7}{2}},
\end{align*}
where we have used $(\ref{b30})$. Hence, we complete the proof of the proposition.
\endpf

\subsection{Optimal decay rates of $\|\nabla^3(\rho,u,\eta)(t)\|_{L^2}$ and $\|\nabla^2\tau(t)\|_{L^2}$}\label{sec4.4}
Inspired by \cite{WW, WZZ}, we are going to study the optimal decay estimates of $\|\nabla^3(\rho,u,\eta)(t)\|_{L^2}$ and $\|\nabla^2\tau(t)\|_{L^2}$. In the process, we have made full use of the benefit of frequency decomposition.
\begin{proposition}\label{lemmas04}Under the same assumptions of Theorem \ref{Theorem2} and $(\ref{clc9})$, we have
\begin{align}
&\|\nabla^3 \rho(t)\|_{L^2}+\|\nabla^3 u(t)\|_{L^2}+\|\nabla^3 \eta(t)\|_{L^2}+\|\nabla^3 \tau(t)\|_{L^2}\leq C(1+t)^{-\frac{9}{4}},\label{a30}\\
&\|\nabla^2\tau(t)\|_{L^2}\leq C(1+t)^{-\frac{9}{4}}.\label{ab34}
\end{align}
\end{proposition}
\begin{remark}
The proof of Proposition \ref{lemmas04} consists of Lemmas \ref{proposition3}-\ref{lemmas03} below.
\end{remark}
\begin{lemma}\label{proposition3} Under the same assumptions of Theorem \ref{Theorem2} and \eqref{clc9}, there exist two small positive constants $\epsilon_7$ and $\epsilon_8$ which will be determined in the proof of the lemma and Lemma \ref{lemmas03}, such that
\begin{align}
&\frac{1}{2}\frac{d}{dt}\mathcal{H}_3(t)
+\epsilon_7\frac{r_1}{4}\|\nabla^3\rho^h\|_{L^2}^2
+\epsilon_8\frac{\beta k\tilde{\eta}}{4}\|\nabla^3u^h\|_{L^2}^2
+\frac{r_2\varepsilon}{2\beta\tilde{\eta}}\|\nabla^4\eta\|_{L^2}^2
+\frac{A_0 r_3}{8\lambda\beta k\tilde{\eta}}\|\nabla^3\tau\|_{L^2}^2
+\frac{r_3\varepsilon}{8\beta
k\tilde{\eta}}\|\nabla^4\tau\|_{L^2}^2\notag\\
\leq&C(\delta+\epsilon_7)\|\nabla^3 \rho^L\|_{L^2}^2 +C(\delta+\epsilon_7+\epsilon_8)\|\nabla^3 u^L\|_{L^2}^{2}
+C(\epsilon_7+\epsilon_8\epsilon+\delta)\|\nabla^3 \eta\|_{L^2}^2,\label{ab24}
\end{align}
where
\begin{align*}
\mathcal{H}_3(t)=&\|\nabla^3 \rho\|_{L^2}^2+\|\nabla^3 u\|_{L^2}^2+\frac{r_2}{\beta\tilde{\eta}}\|\nabla^3\eta\|_{L^2}^2
+\frac{r_3}{2\beta k\tilde{\eta}}\|\nabla^3 \tau\|_{L^2}^2\\
&+\int\big(2\epsilon_7\nabla^2 u:\nabla^3 \rho^h+2\epsilon_8\nabla^2\mathrm{div}\tau:\nabla^2 u^h-\frac{h(\rho)+\beta\rho}{r_1+\beta\rho}|\nabla^3\rho|^2\big)\,\mathrm{d}x.
\end{align*}
\end{lemma}
\pf
Following some arguments similar to the proof of Lemma \ref{lemma-n1} for the case $\ell=3$ and $\mu=\nu=0$, one has
\beq\label{ab25}
\begin{split}
&\frac{1}{2}\frac{d}{dt}\big(\|\nabla^3\rho\|_{L^2}^2
+\|\nabla^3 u\|_{L^2}^2+\frac{r_2}{\beta\tilde{\eta}}\|\nabla^3\eta\|_{L^2}^2
+\frac{r_3}{2\beta
k\tilde{\eta}}\|\nabla^3\tau\|_{L^2}^2
-\int\frac{h(\rho)+\beta\rho}{r_1+\beta\rho}|\nabla^3\rho|^2\,\mathrm{d}x\big)\\
&+\frac{r_2\varepsilon}{2\beta\tilde{\eta}}\|\nabla^4\eta\|_{L^2}^2
+\frac{A_0 r_3}{4\lambda\beta
k\tilde{\eta}}\|\nabla^3\tau\|_{L^2}^2
+\frac{r_3\varepsilon}{4\beta
k\tilde{\eta}}\|\nabla^4\tau\|_{L^2}^2\\
\leq& C\delta\big(\|\nabla^3u\|_{L^2}^2+\|\nabla^3\rho\|_{L^2}^2\big),
\end{split}
\eeq
where we have used the following inequality:
\begin{align*}
\|\nabla^3\big(\frac{1}{r_1+\beta\rho}\big)\|_{L^2}
&\leq C\|\nabla\rho\nabla\rho\nabla\rho\|_{L^2}+C\|\nabla^2\rho\nabla\rho\|_{L^2}+C\|\nabla^3\rho\|_{L^2}\\
&\leq C\|\nabla\rho\|_{L^6}^3+C\|\nabla^2\rho\|_{L^6}\|\nabla\rho\|_{L^3}+C\|\nabla^3\rho\|_{L^2}\\
&\leq C\|\rho\|_{L^\infty}^2\|\nabla^3\rho\|_{L^2}+C\|\nabla^3\rho\|_{L^2}(\|\nabla\rho\|_{L^3}+1)
\end{align*} which is established by H\"older inequality, sobolev inequality and Gagliardo-Nirenberg inequality.

In addition, applying operator $\nabla^2\mathrm{div}$ to the equation of \eqref{17}$_4$ and $\nabla^2$ to the equation of \eqref{17}$_2^h$, multiplying the results by $\nabla^2 u^h$ and $\nabla^2\mathrm{div}\tau$ respectively, then summing them up and integrating it over $\mathbb{R}^3$, we can get
\begin{align}
&\frac{d}{dt}\int\nabla^2\mathrm{div}\tau:\nabla^2 u^h\,\mathrm{d}x-\beta k \tilde{\eta}\int\nabla^2\big(\Delta (u^h+u^L)+\nabla\mathrm{div}(u^h+u^L)\big):\nabla^2 u^h\,\mathrm{d}x\notag\\
=&\int(\nabla^2\mathrm{div}S_4-\frac{A_0}{2\lambda}\nabla^2\mathrm{div}\tau
+\varepsilon\nabla^2\Delta\mathrm{div}\tau):\nabla^2 u^h\,\mathrm{d}x\label{an1}\\
&+\int(\nabla^2 S_2^h-r_1\nabla^3\rho^h-r_2\nabla^3\eta^h
+r_3\nabla^2\mathrm{div}\tau^h):\nabla^2\mathrm{div}\tau\,\mathrm{d}x,\notag
\end{align}
where we have used the frequency decomposition \eqref{5.28}.
Then, similar to the case $\ell=2$ and $\mu=\nu=0$ in \eqref{n6}, we have
\begin{align}
&\int(\nabla^2\mathrm{div}S_4
+\varepsilon\nabla^2\Delta\mathrm{div}\tau):\nabla^2 u^h\,\mathrm{d}x\label{an2}\\
\leq&(\frac{\beta k \tilde{\eta}}{8}+C\delta)\|\nabla^3 u^h\|_{L^2}^2
+C\delta\|\nabla^3 u\|_{L^2}^2+C\delta\|\nabla^3 (\eta,\tau)\|_{L^2}^2+C\|\nabla^4 \tau\|_{L^2}^2,\notag
\end{align}
and
\begin{align}
&\int(\nabla^2 S_2^h-r_1\nabla^3\rho^h-r_2\nabla^3\eta^h
+r_3\nabla^2\mathrm{div}\tau^h):\nabla^2\mathrm{div}\tau\,\mathrm{d}x\label{an3}\\
\leq&\epsilon\|\nabla^3 \rho^h\|_{L^2}^2+C\delta\|\nabla^3 \rho\|_{L^2}^2+C\delta\|\nabla^3 u\|_{L^2}^2+C(\epsilon+\delta)\|\nabla^3 \eta\|_{L^2}^2+C_\epsilon\|\nabla^3\tau\|_{L^2}^2.\notag
\end{align}
Finally, combined with Lemma \ref{lemma}, the rest term of \eqref{an1} can be deal with like
\begin{align}
\int-\frac{A_0}{2\lambda}\nabla^2\mathrm{div}\tau:\nabla^2 u^h\,\mathrm{d}x
\leq&\epsilon\|\nabla^2 u^h\|_{L^2}^2
+C_\epsilon\|\nabla^3 \tau\|_{L^2}^2\label{an4}\\
\leq&C\epsilon\|\nabla^3 u\|_{L^2}^2
+C_\epsilon\|\nabla^3 \tau\|_{L^2}^2.\notag
\end{align}
Together with \eqref{an2}, \eqref{an3} and \eqref{an4}, using H\"older inequality and Cauchy inequality, and choosing $\delta$ and $\epsilon$ small enough, we can deduce from \eqref{an1}
\begin{align}
&\frac{d}{dt}\int\nabla^2\mathrm{div}\tau:\nabla^2 u^h\,\mathrm{d}x+\frac{\beta k \tilde{\eta}}{2}(\|\nabla^3 u^h\|_{L^2}^2+\|\nabla^2\mathrm{div}u^h\|_{L^2}^2)\label{ab26}\\
\leq&C(\epsilon+\delta)\|\nabla^3 \rho^h\|_{L^2}^2+C\delta\|\nabla^3 \rho^L\|_{L^2}^2+C\|\nabla^3 u^L\|_{L^2}^2+C(\epsilon+\delta)\|\nabla^3 \eta\|_{L^2}^2+C_\epsilon\|\nabla^3\tau\|_{H^1}^2,\notag
\end{align}
where we have used the fact that
\[\|f\|_{L^2}\leq \|f^L\|_{L^2}+\|f^h\|_{L^2}, ~\forall f \in L^2(\mathbb{R}^3).\]

Similarly, applying operator $\nabla^2$ to the equation of \eqref{17}$_2$ and \eqref{17}$_1^h$, multiplying the results by $\nabla^3 \rho^h$ and $-\nabla^2\mathrm{div}u$ respectively, then summing them up and integrating it over $\mathbb{R}^3$, we have
\begin{align}
&\frac{d}{dt}\int\nabla^2 u:\nabla^3 \rho^h\,\mathrm{d}x+r_1\int\nabla^3(\rho^L+\rho^h):\nabla^3\rho^h\,\mathrm{d}x\notag\\
=&\int(\nabla^2 S_2-r_2\nabla^3\eta
+r_3\nabla^2\mathrm{div}\tau):\nabla^3 \rho^h\,\mathrm{d}x\label{ab27}\\
&+\int(-\nabla^2 S_1^h
+r_1\nabla^2\mathrm{div}u^h):\nabla^2\mathrm{div}u\,\mathrm{d}x.\notag
\end{align}
Further, from \eqref{ab27} and referring to the case $\ell=2$ and $\mu=\nu=0$ in \eqref{n7}, we can deduce the following inequality:
\begin{align}
&\frac{d}{dt}\int\nabla^2 u:\nabla^3 \rho^h\,\mathrm{d}x+\frac{r_1}{2}\|\nabla^3\rho^h\|_{L^2}^2\label{ab28}\\
\leq&C\|\nabla^3\rho^L\|_{L^2}^2+C\|\nabla^3u\|_{L^2}^2
+C\|\nabla^3\eta\|_{L^2}^2+C\|\nabla^3\tau\|_{L^2}^2.\notag
\end{align}

Hence, $\epsilon_7$\eqref{ab28} together with \eqref{ab25} and $\epsilon_8$\eqref{ab26} yields
\begin{align*}
&\frac{1}{2}\frac{d}{dt}\mathcal{H}_3(t)
+\epsilon_7\frac{r_1}{2}\|\nabla^3\rho^h\|_{L^2}^2
+\epsilon_8\frac{\beta k\tilde{\eta}}{2}\|\nabla^3u^h\|_{L^2}^2
+\frac{r_2\varepsilon}{2\beta\tilde{\eta}}\|\nabla^4\eta\|_{L^2}^2
+\frac{A_0 r_3}{4\lambda\beta k\tilde{\eta}}\|\nabla^3\tau\|_{L^2}^2
+\frac{r_3\varepsilon}{4\beta
k\tilde{\eta}}\|\nabla^4\tau\|_{L^2}^2\\
\leq&C(\epsilon_8\epsilon+\delta)\|\nabla^3 \rho^h\|_{L^2}^2
+C(\delta+\epsilon_7)\|\nabla^3 \rho^L\|_{L^2}^2 +C(\delta+\epsilon_7)\|\nabla^3 u\|_{L^2}^{2}+C\epsilon_8\|\nabla^3 u^L\|_{L^2}^{2}\\
&+C(\epsilon_7+\epsilon_8\epsilon+\delta)\|\nabla^3 \eta\|_{L^2}^2
+C_\epsilon\epsilon_8\|\nabla^3\tau\|_{H^1}^2
+C\epsilon_7\|\nabla^3\tau\|_{L^2}^2.
\end{align*}
Firstly, choosing a fixed positive constant $\epsilon\leq \frac{\beta k \tilde{\eta}r_1}{64C^2}$, and taking
\[\epsilon_7\leq\min\left\{\frac{\epsilon_8\beta k \tilde{\eta}}{8C}, \frac{A_0r_3}{16C\lambda\beta k\tilde{\eta}}\right\}\]
and
\[\epsilon_8\leq\min\left\{\frac{\epsilon_7r_1}{8C\epsilon}, \frac{A_0r_3}{16C_\epsilon\lambda\beta k\tilde{\eta}},\frac{r_3\varepsilon}{8C_\epsilon\beta k \tilde{\eta}}\right\},\]
and finally choosing $\delta$ sufficiently small, we get \eqref{ab24}.
\endpf

\bigskip

Moreover, with Lemmas \ref{proposition3} and $\ref{lemma}$, the following result holds.
\begin{lemma}\label{lemmas03}Under the same assumptions of Theorem \ref{Theorem2} and $(\ref{clc9})$, we have
\begin{align}
\mathcal{H}_3(t)\leq e^{-\tilde{C_3}t}\mathcal{H}_3(0)+C\int_0^te^{-\tilde{C_3}(t-s)}\left(\|\nabla^3 \rho^L\|_{L^2}^2+\|\nabla^3 u^L\|_{L^2}^2+\|\nabla^3 \eta^L\|_{L^2}^2\right)(s)\mathrm{d}s,\label{ab29}
\end{align}
for some positive constant $\tilde{C_3}$ independent of $\delta$.
\end{lemma}
\pf
By Lemma $\ref{lemma}$, we have
\begin{equation*}
c_0\|\nabla^3 \eta^h\|_{L^2}\leq \|\nabla^4 \eta\|_{L^2}.
\end{equation*}
Thus, \eqref{ab24} leads to
\begin{align}
&\frac{1}{2}\frac{d}{dt}\mathcal{H}_3(t)
+\epsilon_7\frac{r_1}{4}\|\nabla^3\rho^h\|_{L^2}^2
+\epsilon_8\frac{\beta k\tilde{\eta}}{4}\|\nabla^3u^h\|_{L^2}^2
+\frac{r_2\varepsilon c_0^2}{2\beta\tilde{\eta}}\|\nabla^3\eta^h\|_{L^2}^2
+\frac{A_0 r_3}{8\lambda\beta k\tilde{\eta}}\|\nabla^3\tau\|_{L^2}^2\label{ab30}\\
\leq&C(\delta+\epsilon_7)\|\nabla^3 \rho^L\|_{L^2}^2 +C(\delta+\epsilon_7+\epsilon_8)\|\nabla^3 u^L\|_{L^2}^{2}
+C(\epsilon_7+\epsilon_8\epsilon+\delta)(\|\nabla^3 \eta^L\|_{L^2}^2+\|\nabla^3 \eta^h\|_{L^2}^2),\notag
\end{align}
By adding $\epsilon_7\frac{r_1}{4}\|\nabla^3 \rho^L\|_{L^2}^2+ \epsilon_8\frac{\beta k\tilde{\eta}}{4}\|\nabla^3 u^L\|_{L^2}^2+\frac{r_2\varepsilon c_0^2}{2\beta\tilde{\eta}}\|\nabla^3\eta^L\|_{L^2}^2$ to both sides of inequality \eqref{ab30}, taking
\[\epsilon_7\leq\frac{r_2\varepsilon c_0^2}{8C\beta\tilde{\eta}}
\;\; \text{and}\;\;
\epsilon_8\leq\frac{r_2\varepsilon c_0^2}{8C\epsilon\beta\tilde{\eta}},\]
and choosing $\delta$ sufficiently small, we have
\begin{align*}
&\frac{1}{2}\frac{d}{dt}\mathcal{H}_3(t)
+\epsilon_7\frac{r_1}{4}\|\nabla^3\rho\|_{L^2}^2
+\epsilon_8\frac{\beta k\tilde{\eta}}{4}\|\nabla^3u\|_{L^2}^2
+\frac{r_2\varepsilon c_0^2}{4\beta\tilde{\eta}}\|\nabla^3\eta\|_{L^2}^2
+\frac{A_0 r_3}{8\lambda\beta k\tilde{\eta}}\|\nabla^3\tau\|_{L^2}^2\\
\leq&C\|\nabla^3 \rho^L\|_{L^2}^2 +C\|\nabla^3 u^L\|_{L^2}^{2}+C\|\nabla^3 \eta^L\|_{L^2}^{2}.
\end{align*}
Moreover, it follows from integration by parts, the Young inequality and Lemma \ref{lemma} that
\begin{align*}
&\int\big(\nabla^2 u:\nabla^3 \rho^h+\nabla^2\mathrm{div}\tau:\nabla^2 u^h\big)\,\mathrm{d}x\\
=&\int\big(-\nabla^2\mathrm{div}u:\nabla^2 \rho^h+\nabla^2\mathrm{div}\tau:\nabla^2 u^h\big)\,\mathrm{d}x\\
\leq&\frac{1}{2}\|\nabla^2\mathrm{div}u\|_{L^2}
+\frac{1}{2}\|\nabla^2 \rho^h\|_{L^2}
+\frac{1}{2}\|\nabla^2\mathrm{div}\tau\|_{L^2}
+\frac{1}{2}\|\nabla^2 u^h\|_{L^2}\\
\leq&\frac{1}{2}\|\nabla^3\,u\|_{L^2}
+\frac{1}{2}\|\nabla^3 \rho\|_{L^2}
+\frac{1}{2}\|\nabla^3\,\tau\|_{L^2}.
\end{align*}
Hence, by virtue of \eqref{clc9} and the smallness of $\delta$, $\epsilon_7$ and $\epsilon_8$, it is easy to check that $\mathcal{H}_3(t)$ is equivalent to $$\|\nabla^3 \rho\|_{L^2}^2+\|\nabla^3 u\|_{L^2}^2+\|\nabla^3 \eta\|_{L^2}^2+\|\nabla^3 \tau\|_{L^2}^2.$$
Then there exists a positive constant $\tilde{C_3}>0$ such that
\begin{align*}
\frac{d}{dt}\mathcal{H}_3(t)+\tilde{C_3}\mathcal{H}_3(t)\leq C \|\nabla^3 \rho^L\|_{L^2}^2+C\|\nabla^3 u^L\|_{L^2}^2+C\|\nabla^3 \eta^L\|_{L^2}^2.
\end{align*}
By using Gronwall's inequality, we get \eqref{ab29}.
\endpf

\bigskip

With the help of Lemmas \ref{proposition3}-\ref{lemmas03}, we are ready to prove Proposition \ref{lemmas04}.
\paragraph{Proof of Proposition \ref{lemmas04}:}
Thanks to the case $m=3$ in \eqref{b33}, \eqref{ab29}, Propositions \ref{lemmas1}, \ref{lemmas2} and \ref{lemmas3}, we can obtain
\begin{align}
\mathcal{H}_3(t)\leq & e^{-\tilde{C_3}t}\mathcal{H}_3(0)+C\int_0^te^{-\tilde{C_3}(t-s)}\left(\|\nabla^3 \rho^L\|_{L^2}^2+\|\nabla^3 u^L\|_{L^2}^2+\|\nabla^3 \eta^L\|_{L^2}^2\right)(s)\mathrm{d}s\notag\\
\leq&e^{-\tilde{C_3}t}\mathcal{H}_3(0)
+C\int_0^te^{-\tilde{C_3}(t-s)} (1+s)^{-\frac{9}{2}}\mathrm{d}s
+C\int_0^te^{-\tilde{C_3}(t-s)} \big(\int_0^\frac{s}{2}(1+s-s')^{-\frac{9}{2}}(1+s)^{-\frac{8}{2}}(s')\mathrm{d}s'\big)\mathrm{d}s\notag\\
&+C\int_0^te^{-\tilde{C_3}(t-s)} \big(\int_\frac{s}{2}^s(1+s-s')^{-3}(1+s)^{-\frac{10}{2}}(s')\mathrm{d}s'\big)\mathrm{d}s\notag\\
\leq&C(1+t)^{-\frac{9}{2}}+C\int_0^te^{-\tilde{C_3}(t-s)}(1+s)^{-\frac{9}{2}}\mathrm{d}s\notag\\
\leq& C(1+t)^{-\frac{9}{2}},\notag
\end{align}
which is $(\ref{a30})$.

Then, for \eqref{ab34}, multiplying $\nabla(\ref{17})_4$ by $2\nabla^2\tau$ and then integrating the result equation over $\mathbb{R}^3$, similar to \eqref{clc-ss01} and \eqref{b34}, we get
$$
\frac{d}{dt}\int|\nabla^2\tau|^2\,\mathrm{d}x+\frac{A_0}{2\lambda}\int|\nabla^2\tau|^2\,\mathrm{d}x
\leq C\|\nabla^3 u\|_{L^2}^2.
$$
Using Gronwall's inequality, the above inequality yields
\begin{align*}
\|\nabla^2\tau(t)\|_{L^2}^2&\leq C e^{-\frac{A_0}{2\lambda}t}\|\nabla^2\tau(0)\|_{L^2}^2+C\int_0^te^{-\frac{A_0}{2\lambda}(t-s)}\|\nabla^3 u(s)\|_{L^2}^2\mathrm{d}s\\
&\leq C e^{-\frac{A_0}{2\lambda}t}\|\nabla^2\tau(0)\|_{L^2}^2+C\int_0^te^{-\frac{A_0}{2\lambda}(t-s)}(1+s)^{-\frac{9}{2}}\mathrm{d}s\\
&\leq C(1+t)^{-\frac{9}{2}},
\end{align*}
where we have used $(\ref{a30})$. Hence, we complete the proof of the proposition.
\endpf

Finally, based on Propositions \ref{lemmas1}-\ref{lemmas04}, the decay rates of the solution stated in Proposition \ref{prodecay} is obtained. Thus, We finish the proof of Theorem \ref{Theorem2}.

\appendix
\section{Appendix}
\subsection{Estimates on the linearized system}

\medskip
Let us consider the following linear system for $(\bar{\rho}, \bar{u}, \bar{\eta}, \mathrm{div}\bar{\tau})$:
\beq\label{A1}
\left\{
\begin{split}
&\bar{\rho}_t+r_1\mathrm{div}\bar{u}=0,\\
&\bar{u}_t+r_1\nabla\bar{\rho}+r_2\nabla\bar{\eta}-r_3\mathrm{div}\bar{\tau}=0,\\
&\bar{\eta}_t+\beta\tilde{\eta}\,\mathrm{div}\bar{u}-\varepsilon\Delta\bar{\eta}=0,\\
&\mathrm{div}\bar{\tau}_t+\frac{A_0}{2\lambda}\mathrm{div}\bar{\tau}-\varepsilon\Delta\mathrm{div}\bar{\tau}-\beta k\tilde{\eta}(\Delta \bar{u}+\nabla\mathrm{div}\bar{u})=0.
\end{split}
\right.
\eeq
As it can be seen from \eqref{e12}, \eqref{b28} and \eqref{ab29}, to study the decay estimates of $(\rho,u,\eta,\tau)$, we only need to analyze the low frequency part ($|\xi|\leq c_0$) of $(\bar{\rho},\bar{u},\bar{\eta},\bar{\tau})$.

If we adopt $\Lambda^s:=(-\Delta)^\frac{s}{2}$ as the notation for the pseudo-differential operator defined by $\Lambda^sf:=\mathcal{F}^{-1}(|\xi|^s\hat{f}(\xi))$, we only need to study $\bar{\rho}$, $d:=\Lambda^{-1}\mathrm{div}\bar{u}$ and $\mathbb{P}\bar{u}:=\Lambda^{-1}\text{curl}\bar{u}$, where $\text{curl}_i^j\bar{u}=\partial_j\bar{u}^i-\partial_i\bar{u}^j$; $\bar{\eta}$, $q:=\Lambda^{-1}\mathrm{div}\mathrm{div}\bar{\tau}$ and $\mathbb{P}\mathrm{div}\bar{\tau}:=\Lambda^{-1}\text{curl}\mathrm{div}\bar{\tau}$. Indeed, by the definition of $\mathbb{P}$, we have
$$\bar{u}=-\Lambda^{-1}\nabla d-\Lambda^{-1}\mathrm{div}\mathbb{P}\bar{u},$$
$$\mathrm{div}\bar{\tau}=-\Lambda^{-1}\nabla q-\Lambda^{-1}\mathrm{div}\mathbb{P}\mathrm{div}\bar{\tau}.$$
We see that $(\bar{\rho},d, \bar{\eta},q)$ and $(\mathbb{P}\bar{u},\mathbb{P}\mathrm{div}\bar{\tau})$ satisfy
\beq\label{A2}
\left\{
\begin{split}
&\bar{\rho}_t+r_1\Lambda d=0,\\
&d_t-r_1\Lambda\bar{\rho}-r_2\Lambda\bar{\eta}-r_3q=0,\\
&\bar{\eta}_t+\beta\tilde{\eta}\Lambda d-\varepsilon\Delta\bar{\eta}=0,\\
&q_t+\frac{A_0}{2\lambda}q-\varepsilon\Delta q-2\beta k\tilde{\eta}\Delta d=0,
\end{split}
\right.
\eeq
and
\beq\label{A3}
\left\{
\begin{split}
&\mathbb{P}\bar{u}_t-r_3\mathbb{P}\mathrm{div}\bar{\tau}=0,\\
&\mathbb{P}\mathrm{div}\bar{\tau}_t+\frac{A_0}{2\lambda}\mathbb{P}\mathrm{div}\bar{\tau}
-\varepsilon\Delta\mathbb{P}\mathrm{div}\bar{\tau}-\beta k\tilde{\eta}\Delta\mathbb{P}\bar{u}=0.
\end{split}
\right.
\eeq
Applying Fourier transform to the linearized system \eqref{A2} and \eqref{A3}, we arrive at
\beq\label{A4}
\left\{
\begin{split}
&\hat{\bar{\rho}}_t+r_1|\xi| \hat{d}=0,\\
&\hat{d}_t-r_1|\xi|\hat{\bar{\rho}}-r_2|\xi|\hat{\bar{\eta}}-r_3\hat{q}=0,\\
&\hat{\bar{\eta}}_t+\beta\tilde{\eta}|\xi| \hat{d}+\varepsilon|\xi|^2\hat{\bar{\eta}}=0,\\
&\hat{q}_t+\frac{A_0}{2\lambda}\hat{q}+\varepsilon|\xi|^2 \hat{q}+2\beta k\tilde{\eta}|\xi|^2 \hat{d}=0,
\end{split}
\right.
\eeq
and
\beq\label{A5}
\left\{
\begin{split}
&\widehat{\mathbb{P}\bar{u}}_t-r_3\widehat{\mathbb{P}\mathrm{div}\bar{\tau}}=0,\\
&\widehat{\mathbb{P}\mathrm{div}\bar{\tau}}_t+\frac{A_0}{2\lambda}\widehat{\mathbb{P}\mathrm{div}\bar{\tau}}
+\varepsilon|\xi|^2\widehat{\mathbb{P}\mathrm{div}\bar{\tau}}+\beta k\tilde{\eta}|\xi|^2\widehat{\mathbb{P}\bar{u}}=0.
\end{split}
\right.
\eeq

\subsubsection{Estimates on $(\hat{\bar{\rho}},\hat{d},\hat{\bar{\eta}},\hat{q})$.}
We introduce the following corrected modes different from those in \cite{W1W}:
\begin{align}
&\hat{\mathbf{a}}=\hat{\bar{\rho}}
+\frac{2\lambda}{A_0}r_3r_1|\xi|[\frac{A_0}{2\lambda}+(\varepsilon-2\beta k\tilde{\eta}\frac{2\lambda}{A_0}r_3)|\xi|^2]^{-1}\hat{q},\label{A-1}\\
&\hat{\mathbf{o}}=\hat{d}+\frac{2\lambda}{A_0}r_3\hat{q},\label{A-2}\\
&\hat{\mathbf{z}}=\hat{\bar{\eta}}
+\frac{2\lambda}{A_0}r_3\beta\tilde{\eta}|\xi|[\frac{A_0}{2\lambda}+(\varepsilon-2\beta k\tilde{\eta}\frac{2\lambda}{A_0}r_3)|\xi|^2]^{-1}\hat{q},\label{A-3}\\
&\hat{\mathbf{q}}=\hat{q}\label{A-4}.
\end{align}
Then the system \eqref{A4} can be rewritten as
\beq\label{A6}
\left\{
\begin{split}
&\hat{\mathbf{a}}_t+r_1|\xi| \hat{\mathbf{o}}
=-\frac{4\lambda}{A_0}r_3r_1\beta k\tilde{\eta}A_1|\xi|^3\hat{\mathbf{o}},\\
&\hat{\mathbf{o}}_t+\frac{4\lambda}{A_0}r_3\beta k\tilde{\eta}|\xi|^2\hat{\mathbf{o}}-r_1|\xi|\hat{\mathbf{a}}
-r_2|\xi|\hat{\mathbf{z}}
=\frac{2\lambda}{A_0}r_3\left(\frac{4\lambda}{A_0}r_3\beta k\tilde{\eta}-\varepsilon-(r_2\beta\tilde{\eta}+r_1^2)
A_1\right)|\xi|^2\hat{\mathbf{q}},\\
&\hat{\mathbf{z}}_t+\beta\tilde{\eta}|\xi|\hat{\mathbf{o}} +\varepsilon|\xi|^2\hat{\mathbf{z}}
=\varepsilon\frac{2\lambda}{A_0}r_3\beta\tilde{\eta}
A_1|\xi|^3\hat{\mathbf{q}}
-\frac{4\lambda}{A_0}r_3\beta^2 k\tilde{\eta}^2
A_1|\xi|^3\hat{\mathbf{o}},\\
&\hat{\mathbf{q}}_t+[\frac{A_0}{2\lambda}+(\varepsilon-2\beta k\tilde{\eta}\frac{2\lambda}{A_0}r_3)|\xi|^2]\hat{\mathbf{q}}
=-2\beta k\tilde{\eta}|\xi|^2\hat{\mathbf{o}},
\end{split}
\right.
\eeq
where the coefficient $A_1$ is defined by
$$A_1=\left[\frac{A_0}{2\lambda}+(\varepsilon-2\beta k\tilde{\eta}\frac{2\lambda}{A_0}r_3)|\xi|^2\right]^{-1}.$$
From the corrected modes, it is not hard to find that the estimates of $(\hat{\mathbf{a}},\hat{\mathbf{o}},\hat{\mathbf{z}},\hat{\mathbf{q}})$ can be easily translated into the estimates of $(\hat{\bar{\rho}},\hat{d},\hat{\bar{\eta}},\hat{q})$ for small $|\xi|$. Next, let us turn to study the estimates of $(\hat{\mathbf{a}},\hat{\mathbf{o}},\hat{\mathbf{z}},\hat{\mathbf{q}})$. From \eqref{A6}, we easily obtain
\begin{align}
&\frac{1}{2}\frac{d}{dt}(|\hat{\mathbf{a}}|^2+|\hat{\mathbf{o}}|^2
+\frac{r_2}{\beta \tilde{\eta}}|\hat{\mathbf{z}}|^2)
+\frac{r_2\varepsilon}{\beta \tilde{\eta}}|\xi|^2|\hat{\mathbf{z}}|^2+\frac{4\lambda}{A_0}r_3\beta k\tilde{\eta}|\xi|^2|\hat{\mathbf{o}}|^2\notag\\
=&-\frac{4\lambda}{A_0}r_3r_1\beta k\tilde{\eta}A_1|\xi|^3\text{Re}(\hat{\mathbf{o}}\,\bar{\hat{\mathbf{a}}})
+\frac{2\lambda}{A_0}r_3\left(\frac{4\lambda}{A_0}r_3\beta k\tilde{\eta}-\varepsilon-(r_2\beta\tilde{\eta}+r_1^2)
A_1\right)|\xi|^2\text{Re}(\hat{\mathbf{q}}\,\bar{\hat{\mathbf{o}}})\label{A7}\\
&+\frac{r_2\varepsilon}{\beta \tilde{\eta}}\frac{2\lambda}{A_0}r_3\beta\tilde{\eta}
A_1|\xi|^3\text{Re}(\hat{\mathbf{q}}\,\bar{\hat{\mathbf{z}}})
-\frac{r_2}{\beta \tilde{\eta}}\frac{4\lambda}{A_0}r_3\beta^2 k\tilde{\eta}^2
A_1|\xi|^3\text{Re}(\hat{\mathbf{o}}\,\bar{\hat{\mathbf{z}}}).\notag
\end{align}
Multiplying \eqref{A6}$_1$ and \eqref{A6}$_2$ by $\bar{\hat{\mathbf{o}}}$ and $\bar{\hat{\mathbf{a}}}$, respectively, yields
\begin{align}
&\frac{d}{dt}\text{Re}(\bar{\hat{\mathbf{a}}}\,\hat{\mathbf{o}}) +r_1|\xi||\hat{\mathbf{o}}|^2-r_1|\xi||\hat{\mathbf{a}}|^2
-r_2|\xi|\text{Re}(\hat{\mathbf{z}}\,\bar{\hat{\mathbf{a}}})\label{A8}\\
=&-\frac{4\lambda}{A_0}r_3r_1\beta k\tilde{\eta}A_1|\xi|^3|\hat{\mathbf{o}}|^2
+\frac{2\lambda}{A_0}r_3\left(\frac{4\lambda}{A_0}r_3\beta k\tilde{\eta}-\varepsilon-(r_2\beta\tilde{\eta}+r_1^2)
A_1\right)|\xi|^2\text{Re}(\hat{\mathbf{q}}\,\bar{\hat{\mathbf{a}}})
-\frac{4\lambda}{A_0}r_3\beta k\tilde{\eta}|\xi|^2\text{Re}(\hat{\mathbf{o}}\,\bar{\hat{\mathbf{a}}}).\notag
\end{align}
Combined with \eqref{A7} and $-\tilde{\epsilon}|\xi|\times$\eqref{A8}, it holds that
\begin{align}
&\frac{1}{2}\frac{d}{dt}\left(|\hat{\mathbf{a}}|^2+|\hat{\mathbf{o}}|^2
+\frac{r_2}{\beta \tilde{\eta}}|\hat{\mathbf{z}}|^2
-2\tilde{\epsilon}|\xi|\text{Re}(\bar{\hat{\mathbf{a}}}\,\hat{\mathbf{o}})\right)\notag
+r_1\tilde{\epsilon}|\xi|^2|\hat{\mathbf{a}}|^2+\frac{r_2\varepsilon}{\beta \tilde{\eta}}|\xi|^2|\hat{\mathbf{z}}|^2
+(\frac{4\lambda}{A_0}r_3\beta k\tilde{\eta}-\tilde{\epsilon}r_1)|\xi|^2|\hat{\mathbf{o}}|^2\notag\\
=&\frac{2\lambda}{A_0}r_3\left(\frac{4\lambda}{A_0}r_3\beta k\tilde{\eta}-\varepsilon-(r_2\beta\tilde{\eta}+r_1^2)
A_1\right)|\xi|^2\text{Re}(\hat{\mathbf{q}}\,\bar{\hat{\mathbf{o}}})
-\tilde{\epsilon}r_2|\xi|^2\text{Re}(\hat{\mathbf{z}}\,\bar{\hat{\mathbf{a}}})\label{A9}\\
&+\tilde{\epsilon}\frac{4\lambda}{A_0}r_3r_1\beta k\tilde{\eta}A_1|\xi|^4|\hat{\mathbf{o}}|^2
-\tilde{\epsilon}\frac{2\lambda}{A_0}r_3\left(\frac{4\lambda}{A_0}r_3\beta k\tilde{\eta}-\varepsilon-(r_2\beta\tilde{\eta}+r_1^2)
A_1\right)|\xi|^3\text{Re}(\hat{\mathbf{q}}\,\bar{\hat{\mathbf{a}}})\notag\\
&+\tilde{\epsilon}\frac{4\lambda}{A_0}r_3\beta k\tilde{\eta}|\xi|^3\text{Re}(\hat{\mathbf{o}}\,\bar{\hat{\mathbf{a}}})
-\frac{4\lambda}{A_0}r_3r_1\beta k\tilde{\eta}A_1|\xi|^3\text{Re}(\hat{\mathbf{o}}\,\bar{\hat{\mathbf{a}}})
+\frac{r_2\varepsilon}{\beta \tilde{\eta}}\frac{2\lambda}{A_0}r_3\beta\tilde{\eta}
A_1|\xi|^3\text{Re}(\hat{\mathbf{q}}\,\bar{\hat{\mathbf{z}}})\notag\\
&-\frac{r_2}{\beta \tilde{\eta}}\frac{4\lambda}{A_0}r_3\beta^2 k\tilde{\eta}^2
A_1|\xi|^3\text{Re}(\hat{\mathbf{o}}\,\bar{\hat{\mathbf{z}}})\notag\\
=&\frac{2\lambda}{A_0}r_3\left(\frac{4\lambda}{A_0}r_3\beta k\tilde{\eta}-\varepsilon-(r_2\beta\tilde{\eta}+r_1^2)
A_1\right)|\xi|^2\text{Re}(\hat{\mathbf{q}}\,\bar{\hat{\mathbf{o}}})
-\tilde{\epsilon}r_2|\xi|^2\text{Re}(\hat{\mathbf{z}}\,\bar{\hat{\mathbf{a}}})
+II_1.\notag
\end{align}
It is natural to derive the estimates for those terms on the right-hand side of \eqref{A9}. First, the first two terms can be controlled by
\beq\label{A10}
\frac{\lambda}{A_0}r_3\beta k\tilde{\eta}|\xi|^2|\hat{\mathbf{o}}|^2
+\frac{ A_0}{4\lambda r_3\beta k\tilde{\eta}}A_2^2|\xi|^2||\hat{\mathbf{q}}|^2
+\frac{\tilde{\epsilon}r_1}{2}|\xi|^2|\hat{\mathbf{a}}|^2
+\frac{\tilde{\epsilon}r_2^2}{2r_1}|\xi|^2||\hat{\mathbf{z}}|^2,
\eeq
where
$$A_2=\frac{2\lambda}{A_0}r_3\left(\frac{4\lambda}{A_0}r_3\beta k\tilde{\eta}-\varepsilon-(r_2\beta\tilde{\eta}+r_1^2)
A_1\right).$$
Similarly, we can drive the bound of the last term.
\beq\label{A12}
|II_1|\leq C|\xi|^3|(\hat{\mathbf{a}},\hat{\mathbf{z}})||(\hat{\mathbf{q}},\hat{\mathbf{o}})|
+C|\xi|^4|\hat{\mathbf{o}}|^2.
\eeq
Substituting \eqref{A10} and \eqref{A12} into \eqref{A9} yields
\begin{align}
&\frac{1}{2}\frac{d}{dt}\left(|\hat{\mathbf{a}}|^2+|\hat{\mathbf{o}}|^2
+\frac{r_2}{\beta \tilde{\eta}}|\hat{\mathbf{z}}|^2
-2\tilde{\epsilon}|\xi|\text{Re}(\bar{\hat{\mathbf{a}}}\,\hat{\mathbf{o}})\right)\notag\\
&+\frac{r_1\tilde{\epsilon}}{2}|\xi|^2|\hat{\mathbf{a}}|^2
+(\frac{r_2\varepsilon}{\beta \tilde{\eta}}-\frac{\tilde{\epsilon}r_2^2}{2r_1})|\xi|^2|\hat{\mathbf{z}}|^2
+(\frac{3\lambda}{A_0}r_3\beta k\tilde{\eta}-\tilde{\epsilon}r_1)|\xi|^2|\hat{\mathbf{o}}|^2\label{A13}\\
\leq&\frac{ A_0}{4\lambda r_3\beta k\tilde{\eta}}A_2^2|\xi|^2||\hat{\mathbf{q}}|^2
+C|\xi|^3|(\hat{\mathbf{a}},\hat{\mathbf{z}})||(\hat{\mathbf{q}},\hat{\mathbf{o}})|
+C|\xi|^4|\hat{\mathbf{o}}|^2.\notag
\end{align}
Now, we move on and derive the estimates of $\hat{\mathbf{q}}$ as follows. Multiplying \eqref{A6}$_4$ by $\bar{\hat{\mathbf{q}}}$ and using Cauchy inequality, we then obtain the following inequality:
\begin{align*}
\frac{1}{2}\frac{d}{dt}|\hat{\mathbf{q}}|^2
+\left[\frac{A_0}{2\lambda}+(\varepsilon-2\beta k\tilde{\eta}\frac{2\lambda}{A_0}r_3)|\xi|^2\right]|\hat{\mathbf{q}}|^2
=&-2\beta k\tilde{\eta}|\xi|^2\text{Re}(\hat{\mathbf{o}}\,\bar{\hat{\mathbf{q}}})\\
\leq&\frac{\lambda}{2A_0}r_3\beta k\tilde{\eta}|\xi|^2||\hat{\mathbf{o}}|^2
+\frac{2A_0}{\lambda r_3}\beta k\tilde{\eta}|\xi|^2|\hat{\mathbf{q}}|^2,
\end{align*}
which implies
\begin{align}
\frac{1}{2}\frac{d}{dt}|\hat{\mathbf{q}}|^2
+\left[\frac{A_0}{2\lambda}+(\varepsilon-2\beta k\tilde{\eta}\frac{2\lambda}{A_0}r_3-\frac{2A_0}{\lambda r_3}\beta k\tilde{\eta})|\xi|^2\right]|\hat{\mathbf{q}}|^2
\leq\frac{\lambda}{2A_0}r_3\beta k\tilde{\eta}|\xi|^2||\hat{\mathbf{o}}|^2.\label{A14}
\end{align}
Finally, combining with \eqref{A13} and \eqref{A14} yields
\begin{align}
&\frac{1}{2}\frac{d}{dt}\left(|\hat{\mathbf{a}}|^2+|\hat{\mathbf{o}}|^2
+\frac{r_2}{\beta \tilde{\eta}}|\hat{\mathbf{z}}|^2+|\hat{\mathbf{q}}|^2
-2\tilde{\epsilon}|\xi|\text{Re}(\bar{\hat{\mathbf{a}}}\,\hat{\mathbf{o}})\right)
+\frac{r_1\tilde{\epsilon}}{2}|\xi|^2|\hat{\mathbf{a}}|^2
+(\frac{r_2\varepsilon}{\beta \tilde{\eta}}-\frac{\tilde{\epsilon}r_2^2}{2r_1})|\xi|^2|\hat{\mathbf{z}}|^2\notag\\
&+(\frac{2\lambda}{A_0}r_3\beta k\tilde{\eta}-\tilde{\epsilon}r_1)|\xi|^2|\hat{\mathbf{o}}|^2
+\left[\frac{A_0}{2\lambda}+(\varepsilon-2\beta k\tilde{\eta}\frac{2\lambda}{A_0}r_3-\frac{2A_0}{\lambda r_3}\beta k\tilde{\eta}-\frac{ A_0}{4\lambda r_3\beta k\tilde{\eta}}A_2^2)|\xi|^2\right]|\hat{\mathbf{q}}|^2\label{A15}\\
\leq&C|\xi|^3|(\hat{\mathbf{a}},\hat{\mathbf{z}})||(\hat{\mathbf{q}},\hat{\mathbf{o}})|
+C|\xi|^4|\hat{\mathbf{o}}|^2.\notag
\end{align}
Taking $\tilde{\epsilon}:=\min\{\frac{r_1\varepsilon}{2r_2\beta \tilde{\eta}},\frac{\lambda r_3\beta k\tilde{\eta}}{A_0r_1}\}$, and introducing the Lyapunov functional
$$\mathcal{L}_{com}^2:=|\hat{\mathbf{a}}|^2+|\hat{\mathbf{o}}|^2
+\frac{r_2}{\beta \tilde{\eta}}|\hat{\mathbf{z}}|^2+|\hat{\mathbf{q}}|^2
-2\tilde{\epsilon}|\xi|\text{Re}(\bar{\hat{\mathbf{a}}}\,\hat{\mathbf{o}}).$$
It is clear that, for $|\xi|\leq \frac{1}{\tilde{\epsilon}}$, we have $\mathcal{L}_{com}^2$ is equivalent to $|\hat{\mathbf{a}}|^2+|\hat{\mathbf{o}}|^2
+\frac{r_2}{\beta \tilde{\eta}}|\hat{\mathbf{z}}|^2+|\hat{\mathbf{q}}|^2$. Hence, from \eqref{A15}, using Cauchy inequality, there exists a small positive constant $c_1\leq\frac{1}{\tilde{\epsilon}}$ depending only on the parameters $\varepsilon, k, A_0, \beta, \tilde{\eta}, \lambda, \tilde{\epsilon}$ and $r_i\;(i=1,2,3)$ such that
\begin{equation*}
\frac{1}{2}\frac{d}{dt}\mathcal{L}_{com}^2
+\frac{r_1\tilde{\epsilon}}{4}|\xi|^2|\hat{\mathbf{a}}|^2
+\frac{r_2\varepsilon}{4\beta \tilde{\eta}}|\xi|^2|\hat{\mathbf{z}}|^2
+\frac{\lambda}{2A_0}r_3\beta k\tilde{\eta}|\xi|^2|\hat{\mathbf{o}}|^2
+\frac{A_0}{4\lambda}|\xi|^2|\hat{\mathbf{q}}|^2
\leq0.
\end{equation*}
Namely, we have
\beq\label{A16}
\frac{d}{dt}\mathcal{L}_{com}^2+C_3|\xi|^2\mathcal{L}_{com}^2\leq0,
\eeq
for some positive constant $C_3$ independent of $|\xi|$. Then it follows form \eqref{A16} that for $\xi\leq c_1$
\beq\label{A17}
|(\hat{\mathbf{a}},\hat{\mathbf{o}},\hat{\mathbf{z}},\hat{\mathbf{q}})(t)|^2
\leq Ce^{-C_3|\xi|^2t}|(\hat{\mathbf{a}},\hat{\mathbf{o}},\hat{\mathbf{z}},\hat{\mathbf{q}})(0)|^2.
\eeq
Recalling the relation between $(\hat{\mathbf{a}},\hat{\mathbf{o}},\hat{\mathbf{z}},\hat{\mathbf{q}})$ with $(\hat{\bar{\rho}},\hat{d},\hat{\bar{\eta}},\hat{q})$ in \eqref{A-1}-\eqref{A-4}, and thanks to \eqref{A17}, we can easily deduce the following inequality holds.
\beq\label{A18}
|(\hat{\bar{\rho}},\hat{d},\hat{\bar{\eta}},\hat{q})(t)|^2
\leq Ce^{-C_3|\xi|^2t}|(\hat{\bar{\rho}},\hat{d},\hat{\bar{\eta}},\hat{q})(0)|^2,
\;\;\text{for}\;|\xi|\leq c_1.
\eeq

\subsubsection{Estimates on $(\widehat{\mathbb{P}\bar{u}},\widehat{\mathbb{P}\mathrm{div}\bar{\tau}})$}
We introduce the following corrected modes:
\begin{align}
&\hat{\mathbf{v}}=\widehat{\mathbb{P}\bar{u}}
+\frac{2\lambda}{A_0}r_3\widehat{\mathbb{P}\mathrm{div}\bar{\tau}},\label{A-5}\\
&\hat{\mathbf{w}}=\widehat{\mathbb{P}\mathrm{div}\bar{\tau}},\label{A-6}
\end{align}
Then the system \eqref{A5} can be rewritten as
\beq\label{A5-}
\left\{
\begin{split}
&\hat{\mathbf{v}}_t+\frac{2\lambda}{A_0}r_3\beta k\tilde{\eta}|\xi|^2 \hat{\mathbf{v}}
+\frac{2\lambda}{A_0}r_3\big(\varepsilon-\frac{2\lambda}{A_0}r_3\beta k\tilde{\eta}\big)|\xi|^2\hat{\mathbf{w}}=0,\\
&\hat{\mathbf{w}}_t+\frac{A_0}{2\lambda}\hat{\mathbf{w}}
+\big(\varepsilon-\frac{2\lambda}{A_0}r_3\beta k\tilde{\eta}\big)|\xi|^2\hat{\mathbf{w}}
+\beta k\tilde{\eta}|\xi|^2\hat{\mathbf{v}}
=0.
\end{split}
\right.
\eeq
Multiplying \eqref{A5-} with $\bar{\hat{\mathbf{v}}}$ and $\bar{\hat{\mathbf{w}}}$ respectively and use Cauchy inequality, we have
\begin{align}
\frac{1}{2}\frac{d}{dt}|\hat{\mathbf{v}}|^2
+\frac{2\lambda}{A_0}r_3\beta k\tilde{\eta}|\xi|^2|\hat{\mathbf{v}}|^2=&
-\frac{2\lambda}{A_0}r_3\big(\varepsilon-\frac{2\lambda}{A_0}r_3\beta k\tilde{\eta}\big)|\xi|^2\text{Re}(\hat{\mathbf{w}}\,\bar{\hat{\mathbf{v}}})\label{AA1}\\
\leq&\frac{\lambda}{2 A_0}r_3\beta k\tilde{\eta}|\xi|^2|\hat{\mathbf{v}}|^2
+\frac{2\lambda}{A_0\beta k\tilde{\eta}}r_3\big(\varepsilon-\frac{2\lambda}{A_0}r_3\beta k\tilde{\eta}\big)^2|\xi|^2|\hat{\mathbf{w}}|^2,\notag
\end{align}
and
\begin{align}
\frac{1}{2}\frac{d}{dt}|\hat{\mathbf{w}}|^2
+\big(\frac{A_0}{2\lambda}+(\varepsilon-\frac{2\lambda}{A_0}r_3\beta k\tilde{\eta})|\xi|^2\big)|\hat{\mathbf{w}}|^2=&
-\beta k\tilde{\eta}|\xi|^2\text{Re}(\hat{\mathbf{v}}\,\bar{\hat{\mathbf{w}}})\label{AA2}\\
\leq&\frac{ A_0}{2\lambda r_3}\beta k\tilde{\eta}|\xi|^2|\hat{\mathbf{w}}|^2
+\frac{\lambda}{2A_0}r_3\beta k\tilde{\eta}|\xi|^2|\hat{\mathbf{v}}|^2.\notag
\end{align}
Together with \eqref{AA1} and \eqref{AA2}, it holds that
\begin{align}
&\frac{1}{2}\frac{d}{dt}(|\hat{\mathbf{v}}|^2+|\hat{\mathbf{w}}|^2)
+\frac{\lambda}{A_0}r_3\beta k\tilde{\eta}|\xi|^2|\hat{\mathbf{v}}|^2\label{AA3}\\
&+\left[\frac{A_0}{2\lambda}+\big(\varepsilon-\frac{2\lambda}{A_0}r_3\beta k\tilde{\eta}-\frac{A_0}{2\lambda r_3}\beta k\tilde{\eta}-\frac{2\lambda}{A_0\beta k\tilde{\eta}}r_3(\varepsilon-\frac{2\lambda}{A_0}r_3\beta k\tilde{\eta})^2\big)|\xi|^2\right]|\hat{\mathbf{w}}|^2\leq0.\notag
\end{align}
By choosing $|\xi|\leq c_2$, which makes  $\frac{A_0}{2\lambda}+\big(\varepsilon-\frac{2\lambda}{A_0}r_3\beta k\tilde{\eta}-\frac{A_0}{2\lambda r_3}\beta k\tilde{\eta}-\frac{2\lambda}{A_0\beta k\tilde{\eta}}r_3(\varepsilon-\frac{2\lambda}{A_0}r_3\beta k\tilde{\eta})^2\big)|\xi|^2\geq\frac{A_0}{4\lambda}\geq\frac{A_0}{4\lambda}|\xi|^2$ true.  Then, from \eqref{AA3}, we have
\[\frac{1}{2}\frac{d}{dt}(|\hat{\mathbf{v}}|^2+|\hat{\mathbf{w}}|^2)
+\frac{\lambda}{A_0}r_3\beta k\tilde{\eta}|\xi|^2|\hat{\mathbf{v}}|^2
+\frac{A_0}{4\lambda}|\xi|^2|\hat{\mathbf{w}}|^2\leq0.\]
Hence, there exists a positive constant $C_4$ independent of $|\xi|$ such that
\begin{equation}\label{AA4}
\frac{d}{dt}(|\hat{\mathbf{v}}|^2+|\hat{\mathbf{w}}|^2)
+C_4|\xi|^2(|\hat{\mathbf{v}}|^2+|\hat{\mathbf{w}}|^2)\leq0.
\end{equation}
Thanks to \eqref{AA4}, we can deduce
\begin{equation}\label{AA5}
|(\hat{\mathbf{v}},\hat{\mathbf{w}})(t)|^2
\leq Ce^{-C_4|\xi|^2 t}|(\hat{\mathbf{v}},\hat{\mathbf{w}})(0)|^2,\;\;\text{for}\;\; |\xi|\leq c_2.
\end{equation}
Further, from the definitions \eqref{A-5} and \eqref{A-6}, \eqref{AA5} implies that, for $|\xi|\leq c_2$,
\beq\label{A27}
|(\widehat{\mathbb{P}\bar{u}},\widehat{\mathbb{P}\mathrm{div}\bar{\tau}})(t)|^2
\leq Ce^{-C_4|\xi|^2 t}|(\widehat{\mathbb{P}\bar{u}},\widehat{\mathbb{P}\mathrm{div}\bar{\tau}})(0)|^2.
\eeq

\subsubsection{Decay estimates of the low frequency part}
Taking $c_0=\min\{c_1, c_2\}$, and combining with \eqref{A18} and \eqref{A27}, we have the following proposition.
\begin{proposition}\label{proA1}It holds that, for $|\xi|\leq c_0$,
$$
|(\hat{\bar{\rho}},\hat{d},\hat{\bar{\eta}},\hat{q})(\xi,t)|^2
\leq Ce^{-2C_{5}|\xi|^2 t}|(\hat{\bar{\rho}},\hat{d},\hat{\bar{\eta}},\hat{q})(\xi,0)|^2
$$
and
$$
|(\widehat{\mathbb{P}\bar{u}},\widehat{\mathbb{P}\mathrm{div}\bar{\tau}})(\xi,t)|^2
\leq Ce^{-2C_{5}|\xi|^2 t}|(\widehat{\mathbb{P}\bar{u}},\widehat{\mathbb{P}\mathrm{div}\bar{\tau}})(\xi,0)|^2,
$$
for some positive constant $C_5$ independent of $\xi$.
\end{proposition}

For the low frequency part, we have the following decay estimates for solutions to the linearized problem \eqref{A1}, \eqref{initial-condition1}.

\begin{proposition}\label{proA2}It holds that
$$
\|\partial_x^m(\bar{\rho},\bar{u},\bar{\eta},\mathrm{div}\bar{\tau} )^L(t)\|_{L^2}
\leq C(1+t)^{-\frac{3}{4}-\frac{|m|}{2}}\|(\bar{\rho},\bar{u},\bar{\eta},\mathrm{div}\bar{\tau} )^L(0)\|_{L^1},
$$
and
$$
\|\partial_x^m(\bar{\rho},\bar{u},\bar{\eta},\mathrm{div}\bar{\tau} )^L(t)\|_{L^2}
\leq C(1+t)^{-\frac{|m|}{2}}\|(\bar{\rho},\bar{u},\bar{\eta},\mathrm{div}\bar{\tau} )^L(0)\|_{L^2}
$$
for any $|m|>0$.
\end{proposition}
\pf
By using Plancherel's theorem and Proposition \ref{proA1}, we have
\begin{align}
\|\partial_x^m (\bar{\rho},d,\bar{\eta},q)^L(t)\|_{L^2}
&=\|(i\xi)^m(\hat{\bar{\rho}},\hat{d},\hat{\bar{\eta}},\hat{q})
\|_{L_{\xi}^2(|\xi|\leq c_0)}\notag\\
&\leq\left(\int_{|\xi|\leq c_0}|\xi|^{2|m|}|(\hat{\bar{\rho}},\hat{d},\hat{\bar{\eta}},\hat{q})(\xi, t)|^2\mathrm{d}\xi\right)^\frac{1}{2}\label{A52}\\
&\leq C\left(\int_{|\xi|\leq c_0}e^{2C_5|\xi|^2}|\xi|^{2|m|}e^{-2C_5|\xi|^2(1+t)}|(\hat{\bar{\rho}},\hat{d},\hat{\bar{\eta}},\hat{q})(\xi, 0)|^2\mathrm{d}\xi\right)^\frac{1}{2}\notag\\
&\leq C\left(\int_{|\beta|\leq c_0\sqrt{1+t}}|\beta|^{2|m|}(1+t)^{-|m|-\frac{3}{2}}e^{-2C_5|\beta|^2}|(\hat{\bar{\rho}},\hat{d},\hat{\bar{\eta}},\hat{q})(\xi, 0)|^2\mathrm{d}\beta\right)^\frac{1}{2}\notag\\
&\leq C(1+t)^{-\frac{3}{4}-\frac{|m|}{2}}\|(\hat{\bar{\rho}},\hat{d},\hat{\bar{\eta}},\hat{q})(0)\|_{L^\infty}\notag\\
&\leq C(1+t)^{-\frac{3}{4}-\frac{|m|}{2}}\|(\bar{\rho},d,\bar{\eta},q)(0)\|_{L^1}.\notag
\end{align}
Moreover, we can also deduce
\[\|\partial_x^m (\bar{\rho},d,\bar{\eta},q)^L(t)\|_{L^2}
\leq C(1+t)^{-\frac{|m|}{2}}\|(\bar{\rho},d,\bar{\eta},q)(0)\|_{L^2}.\]

By similar calculations, we can get
$$
\|\partial_x^m (\mathbb{P}\bar{u},\mathbb{P}\mathrm{div}\bar{\tau})^L(t)\|_{L^2}
\leq C(1+t)^{-\frac{3}{4}-\frac{|m|}{2}}\|(\mathbb{P}\bar{u},\mathbb{P}\mathrm{div}\bar{\tau})(0)\|_{L^1}.
$$
and
$$
\|\partial_x^m (\mathbb{P}\bar{u},\mathbb{P}\mathrm{div}\bar{\tau})^L(t)\|_{L^2}
\leq C(1+t)^{-\frac{|m|}{2}}\|(\mathbb{P}\bar{u},\mathbb{P}\mathrm{div}\bar{\tau})(0)\|_{L^2}.
$$
Then, we finish the proof of Proposition \ref{proA2}.
\endpf
\subsection{Some useful inequalities}
Let $0\leq\varphi_0(\xi)\leq1$ be a function in $\mathcal{C}_0^\infty(\mathbb{R}^3)$ such that
\begin{eqnarray*}
\varphi_0(\xi)=
\begin{cases}
1, & {|\xi|\leq \frac{c_0}{2},}\\
0, & {|\xi|\geq c_0,}
\end{cases}
\end{eqnarray*}
where $c_0$ is a positive constant. Based on the Fourier transform, we can define a low and high frequency decomposition $(f^L(x), f^h(x))$ for a function $f(x)$ as follows
\begin{equation}\label{5.28}
f^L(x)=\mathcal{F}^{-1}(\varphi_0(\xi)\hat{f}(\xi)), \;\;and\;\;f^h(x)=f(x)-f^L(x).
\end{equation}
The following lemma can be obtained directly from the definition \eqref{5.28} and Plancherel's theorem.

\begin{lemma}\label{lemma}(\cite{LLW}) If $f\in H^m(\mathbb{R}^3)(m \geq 2)$ is divided into two parts $(f^L, f^h)$ by the low and high frequency decomposition \eqref{5.28}. It holds that
\begin{equation*}
c_0^{m_1-m_2}\|\nabla^{m_2}f^h\|_{L^2(\mathbb{R}^3)}\leq\|\nabla^{m_1}f\|_{L^2(\mathbb{R}^3)},
\end{equation*} for any integers $m_1$ and $m_2$ with $m_2 \leq m_1\leq m$.
\end{lemma}

Finally, the following elementary inequality will also be used.
\begin{lemma}(\cite{DU}) If $a>1$ and $b\in [0,a]$, then it holds that
\begin{equation*}
\int_0^t(1+t-s)^{-a}(1+s)^{-b}\mathrm{d}s\leq C(a,b) (1+t)^{-b}.
\end{equation*}
\end{lemma}

\begin{lemma}\label{LemmaA2}(\cite{DRZ,MB})
Let $m\geq 1$ be an integer, then we have
\begin{equation*}
\|\nabla^m(fg)\|_{L^p(\mathbb{R}^n)}\leq C\|f\|_{L^{p_1}(\mathbb{R}^n)}\|\nabla^m g\|_{L^{p_2}(\mathbb{R}^n)}+
C\|\nabla^m f\|_{L^{p_3}(\mathbb{R}^n)}\|g\|_{L^{p_4}(\mathbb{R}^n)},
\end{equation*}
where $1\leq p, p_i\leq +\infty,\ (i=1,2,3,4)$ and
\begin{equation*}
\frac{1}{p}=\frac{1}{p_1}+\frac{1}{p_2}=\frac{1}{p_3}+\frac{1}{p_4}.
\end{equation*}
\end{lemma}
%
%

\section*{Acknowledgement}
W. Wang was supported by the National Natural Science Foundation of China $\#$ 11871341.
H. Wen was partially supported by the National Natural Science Foundation of China $\#$ 12071152.

\end{document}